\documentclass[11pt]{article}
\usepackage{natbib,amsmath,amssymb,amsfonts,dsfont,enumitem,float,pifont,mathtools}
\usepackage[dvipsnames]{xcolor}
\RequirePackage[colorlinks,linkcolor=red!80!black,citecolor=green!50!black,urlcolor=blue,breaklinks]{hyperref}
\usepackage[normalem]{ulem}
\usepackage{epsfig,longtable,xcolor}
\usepackage{booktabs,float}
\usepackage{graphicx}
\usepackage{subcaption}
\usepackage{algorithm}
\usepackage{algpseudocode}

\usepackage{authblk}

\usepackage{fancyhdr}
\usepackage{bbm}

\usepackage[utf8]{inputenc}

\title{Fundamental limits of community detection from multi-view data: multi-layer, dynamic and partially labeled block models}
\author{Xiaodong Yang\footnote{\textit{email}: xyang@g.harvard.edu}}
\author{Buyu Lin\footnote{\textit{email}: blin@g.harvard.edu}} 
\author{Subhabrata Sen\footnote{\textit{email}: subhabratasen@fas.harvard.edu} }
\affil{Department of Statistics, Harvard University}
\date{\today}

\usepackage{fancyhdr}
\numberwithin{equation}{section}
\usepackage{graphicx}
\usepackage{pgfplots}
\pgfplotsset{width=0.6\textwidth,compat=1.13}
\usepackage{epsfig}
\usepackage{chngcntr}

\usepackage{authblk}
\usepackage{amsthm}
\usepackage{tikz}
\usetikzlibrary{fit, shapes.geometric, positioning}
\newtheorem{Theorem}{Theorem}[section]
\newtheorem{Definition}{Definition}[section]
\newtheorem{Example}{Example}[section]
\newtheorem{Corollary}{Corollary}[section]
\newtheorem{Conjecture}
{Conjecture}[section]
\newtheorem{Proposition}{Proposition}[section]
\newtheorem{Lemma}{Lemma}[section]
\newtheorem{Assumption}{Assumption}[section]
\newtheorem{Remark}{Remark}[section]

\newcommand{\bM}{\mathbf{M}}

\newcommand{\bN}{\mathbf{N}}

\newcommand{\bV}{\mathbf{V}}
\newcommand{\bW}{\mathbf{W}}

\newcommand{\bX}{\mathbf{X}}
\newcommand{\bx}{\mathbf{x}}

\newcommand{\bY}{\mathbf{Y}}

\newcommand{\by}{\mathbf{y}}

\newcommand{\bZ}{\mathbf{Z}}

\newcommand{\Tr}{\mathrm{Tr}}

\newcommand{\MMSE}{\mathsf{MMSE}}
\newcommand{\DMSE}{\mathsf{DMSE}}
\newcommand{\MSE}{\mathsf{MSE}}
\newcommand{\AMP}{\mathsf{AMP}}

\newcommand{\R}{\mathbb{R}}
\newcommand{\E}{\mathbb{E}}
\newcommand{\bb}{{\mathbf{b}}}
\newcommand{\bA}{\mathbf{A}}
\newcommand{\bG}{\mathbf{G}}

\newcommand{\diag}{\text{diag}}
\newcommand{\val}{\mathrm{val}}
\newcommand{\blam}{\boldsymbol{\lambda}}
\newcommand{\bgamma}{\boldsymbol{\gamma}}
\newcommand{\lambdal}{\lambda^{(l)}}
\newcommand{\bv}{\mathbf{v}}

\newcommand{\hbX}{\widehat{\mathbf{X}}}

\newcommand{\indsim}{\overset{\mathrm{ind}}{\sim}}
\newcommand{\nconverge}{\underset{n\rightarrow\infty}{\longrightarrow}}
\def\toW{\overset{W}{\to}}

\def\limval{\operatorname{lim-val}}

\newcommand{\ML}{\mathsf{ML}}
\newcommand{\Dyn}{\mathsf{Dyn}}
\newcommand{\Semi}{\mathsf{Semi}}

\begin{document}

\maketitle

\begin{abstract}
    Multi-view data arises frequently in modern network analysis e.g. relations of multiple types among individuals in social network analysis, longitudinal measurements of interactions among observational units, annotated networks with noisy partial labeling of vertices etc. We study community detection in these disparate settings via a unified theoretical framework, and investigate the fundamental thresholds for community recovery. We characterize the mutual information between the data and the latent parameters, provided the degrees are sufficiently large. Based on this general result, (i) we derive a sharp threshold for community detection in an inhomogeneous multilayer block model \citep{chen2022global}, (ii) characterize a sharp threshold for weak recovery in a dynamic stochastic block model \citep{matias2017statistical}, and (iii) identify the limiting mutual information in an unbalanced partially labeled block model. Our first two results are derived modulo coordinate-wise convexity assumptions on specific functions---we provide extensive numerical evidence for their correctness. Finally, we introduce iterative algorithms based on Approximate Message Passing for community detection in these problems.
 \end{abstract}


\section{Introduction}
The discovery of latent communities is a fundamental task in modern network analysis. This problem has attracted widespread attention in probability, statistics, computer science, statistical physics and social sciences over the last three decades,  leading to a detailed understanding of the basic statistical thresholds for signal recovery, and the introduction of statistically optimal algorithms for community detection. We refer the interested reader to \cite{abbe2017community} for a survey of the recent progress in this research endeavor. 

Networks are traditionally used to represent pairwise relations among interacting units. However, modern network data is increasingly sophisticated. For example, one frequently observes multiple types of relations among the observational units in modern networks. In particular, in the context of online social networks, one might observe interactions among individuals over professional and personal social networks \citep{kivela2014multilayer}. As a second prominent scenario, one might observe the evolution of a network over time. Such data is valuable in biology \citep{bakken2016comprehensive}, planning \citep{li2017fundamental}, epidemiology \citep{masuda2017introduction} etc. Finally, in addition to the network, one often has access to the true community assignments for some vertices in the network---this setting is common in recommendation systems \citep{shapira2013facebook}, the study of co-citation networks \citep{ji2022co}, protein classification \citep{weston2003semi} etc.  
These diverse examples may be conceptually unified under the setting of ``multi-view" network data. The multi-view formalism has become increasingly important for supervised learning in genomics and proteomics \citep{ding2022cooperative}; the network settings described above can be neatly merged through this framework---indeed, in the multi-type application described above, each type represents an independent view of the data. Similarly, in the temporal setting, the network observed at each time point represents an additional view. Finally, for the partially labeled setting, one can conceptually decouple the data into two correlated views: the first view comprises the unlabeled network, while the second view collects the noisy, incomplete vertex information. The community detection problem remains vital in these multi-view network problems; however, our understanding of the basic statistical limits of these problems are quite restrictive. We refer the interested reader to Section \ref{sec:lit_review} for a review of the current state-of-the-art. 

The current paper introduces a general framework to study community detection in multi-view network data. The present work has two main contributions: 
\begin{itemize}
    \item[(i)] We derive the asymptotic per-vertex mutual information between the latent signal and the observed data in the limit of large system size, under an additional diverging average degree assumption. 
    \item[(ii)] We introduce a family of iterative algorithms based on Approximate Message Passing (AMP) for signal recovery under our general framework. 
\end{itemize}
Using the asymptotic mutual information, we conjecture a precise threshold for community detection in an inhomogeneous multi-layer stochastic block model \cite{chen2022global}. We also identify the community detection threshold for temporal networks and observe an interesting dichotomy in this setting---in one case, community detection is possible with enough temporal observations; in the complementary regime, community detection is impossible, even with long term temporal observations of the networks. Finally, we characterize the limiting mutual information in an unbalanced partially labeled block model.

\subsection{A general model}\label{subsec:general model setup}
We start with a general inference problem that will encompass the specific applications of interest. To this end, we start with a latent vector of the form $(x,y,z)$, where (i) $x \in  \{\pm 1\}^L$ represents a latent signal directly associated with the observed data, (ii) $ y \in \{ \pm 1\}^{L_1}$ represents an implicit latent signal to be indirectly inferred from the data and (iii) $z \in \mathcal{Z}$ represents some observed side information for each node.

\begin{figure}
  \centering
  \begin{subfigure}[b]{0.45\textwidth}
    \centering
    \begin{tikzpicture}[scale=0.7, transform shape][level distance=1.5cm,
                        level 1/.style={sibling distance=1cm}]
      \node[circle, draw, minimum size=1.2cm] (root) {$Y$}
        child {node[circle, draw] (x1) {$X^{(1)}$}
            child {node[circle, draw] {$A^{(1)\prime}$}}}
        child {node[circle, draw] (x2) {$X^{(2)}$}
            child {node[circle, draw] {$A^{(2)\prime}$}}}
        child {node[circle, draw] (x3) {$X^{(3)}$}
            child {node[circle, draw] {$A^{(3)\prime}$}}}
        child {node[circle, draw] (x4) {$X^{(4)}$}
            child {node[circle, draw] {$A^{(4)\prime}$}}}
        child {node[circle, draw] (x5) {$X^{(5)}$}
            child {node[circle, draw] {$A^{(5)\prime}$}}};
      \node[draw, fit=(root) (x1) (x2) (x3) (x4) (x5), inner sep=0.1cm, red, dashed, rounded corners, line width=0.8pt, label=above:$\textcolor{red}{p_{\ML}}$] {};
    \end{tikzpicture}
    \caption{Illustration of $p_\ML$ with $L=5$.}\label{fig:ML illustration}
  \end{subfigure}
  \hspace{0.5cm} 
  \begin{subfigure}[b]{0.45\textwidth}
    \centering
    \begin{tikzpicture}[scale=0.7, transform shape]
      \node[circle, draw] (node1) {$X^{(1)}$};
      \node[circle, draw, right=of node1] (node2) {$X^{(2)}$};
      \node[circle, draw, right=of node2] (node3) {$X^{(3)}$};
      \node[circle, draw, right=of node3] (node4) {$X^{(4)}$};
      \node[circle, draw, right=of node4] (node5) {$X^{(5)}$};

      \node[circle, draw, below=of node1] (node1a) {$A^{(1)\prime}$};
      \node[circle, draw, below=of node2] (node2a) {$A^{(2)\prime}$};
      \node[circle, draw, below=of node3] (node3a) {$A^{(3)\prime}$};
      \node[circle, draw, below=of node4] (node4a) {$A^{(4)\prime}$};
      \node[circle, draw, below=of node5] (node5a) {$A^{(5)\prime}$};
      
      \draw (node1) -- (node2) -- (node3) -- (node4) -- (node5);
      \draw (node1) -- (node1a);
      \draw (node2) -- (node2a);
      \draw (node3) -- (node3a);
      \draw (node4) -- (node4a);
      \draw (node5) -- (node5a);
      \node[draw, fit=(node1) (node2) (node3) (node4) (node5), inner sep=0.15cm, red, dashed, rounded corners, line width=0.8pt, label=above:$\textcolor{red}{p_{\Dyn}}$] {};
    \end{tikzpicture}
    \caption{Illustration of $p_\Dyn$ with $L=5$.}\label{fig:Dyn illustration}
  \end{subfigure}
  \caption{Illustration of effective scalar channels~\eqref{eq:scalar channel} under inter-layer priors of multilayer and dynamic SBMs respectively in Example~\ref{eg:multilayer SBM} and \ref{eg:dynamic SBM}.}\label{fig:prior illustration}
\end{figure}
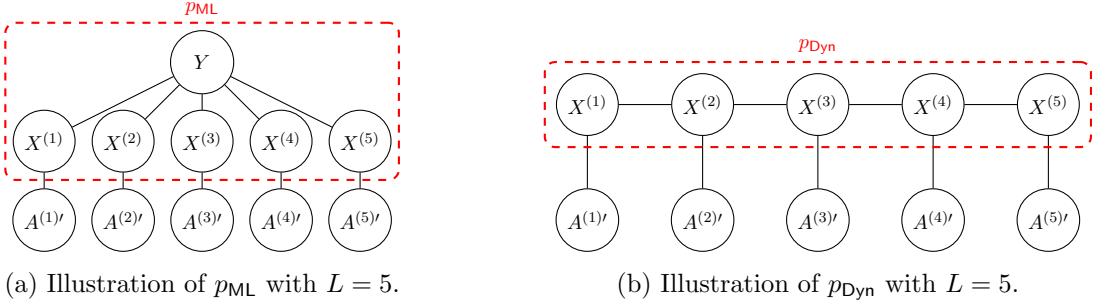

We formulate our inference problem in the context of an explicit generative model: assume a prior $p(x,y,z)$ on the variables $(x,y,z)$. For the simplicity of presentation, we restrict ourselves to a balanced setting in which marginals $p\left(x^{(l)}\right)$, $p\left(y^{(l_1)}\right)$ are all uniform on $\{\pm 1\}$ under the prior.
Our observations are specified as follows---generate $(X_i,Y_i,Z_i)\indsim p$ for $i\in[n]$; then sample $L$ random graphs $\{\bG^{(l)}:l\in[L]\}$ with adjacency matrix
\begin{equation}\label{eq:original graph model}
    G_{i,j}^{(l)}\indsim\text{Bern}\left(\frac{a^{(l)}}{n}\right)\mathbbm{1}\left\{X_i^{(l)}=X_j^{(l)}\right\}+\text{Bern}\left(\frac{b^{(l)}}{n}\right)\mathbbm{1}\left\{X_i^{(l)}\neq X_j^{(l)}\right\},\quad\forall 1\le i<j\le n.
\end{equation}
The parameters $a^{(l)}$, $b^{(l)}$ can depend on $n$---we suppress this dependence for notational convenience. 
We observe the $L$ graphs $\{\bG^{(l)}: l \in [L]\}$; in addition, we also observe the $Z_i$ variable as side information for each node $i \in [n]$. We seek to recover $\{(X_i,Y_i),i\in[n]\}$ jointly from $\bG$ and $\{Z_i,i\in[n]\}$. Three running examples of our framework are given below.
\begin{Example}[Inhomogeneous Multilayer SBM]\label{eg:multilayer SBM}
    Set $L_1=1$ and $\mathcal{Z}=\emptyset$ so that there is no per-node side information. The prior distribution $p_\ML(x,y)$ on $\{\pm 1\}^{L+1}$ is specified as follows: first, set $y$ to be uniform $p_\ML(y=\pm 1)=\frac{1}{2}$ and then let $(x^{(1)},\ldots,x^{(L)})$ be conditionally independent,
    \begin{equation*}
        p_\ML(x^{(l)}=y)=1-\rho,\quad p_\ML(x^{(l)}=-y)=\rho,\quad l\in[L].
    \end{equation*}
    In this example, $y$ is called the global membership and $x^{(l)}$ is referred to as the individualized membership of a certain layer. See Figure~\ref{fig:ML illustration} for a graphical illustration of $p_\ML$.
\end{Example}


We refer the interested reader to  \cite{chen2022global} for a survey of prior work and state-of-the-art results in this model.
Before proceeding further, we provide a quick glimpse into the high-level takeaway messages from the analysis in this paper. As a concrete consequence of our general investigations, we derive a sharp conjecture on the threshold for community detection in this model. 

\begin{Conjecture}
Set $d^{(l)} :=  \frac{a^{(l)}+ b^{(l)}}{2}$ and $d = \min_{1\leq l \leq L} d^{(l)}$. Define $\lambda_n^{(l)}= \frac{(a^{(l)} - b^{(l)})^2}{4 d^{(l)} (1- d^{(l)}/n)}$ and assume that $\lambda_n^{(l)} \to \lambda^{(l)}$ as $n \to \infty$.  Fix $0\leq \rho < 1/2$. If $d \to \infty$, weak recovery of the global and layer-wise memberships are possible if and only if 
\begin{align}
    \max\left\{ \sum_{l=1}^{L} \frac{(1-2\rho)^4 \lambda^{(l)}}{1-(1-(1-2\rho)^4)\lambda^{(l)}} , \boldsymbol{\lambda} \right\}>1. \label{eq:conj1}
\end{align}
\end{Conjecture}

Note that for $L=1$, this threshold reduces to the celebrated threshold for community detection in the Stochastic Block Model \citep{krzakala2013spectral,mossel2012stochastic,mossel2018proof,bordenave2015non}. In addition, if $\rho=0$, the threshold reduces to $\sum_{l=1}^{L}\lambda^{(l)}>1$, recovering the threshold established in \cite{ma2023community}. We provide a partial proof of this conjecture. We establish that the latent communities can be recovered weakly for almost all $\boldsymbol{\lambda}$ satisfying \eqref{eq:conj1}, provided the minimum degree $d$ is sufficiently large. Conversely, we establish that if $d \to \infty$, weak recovery is impossible if \eqref{eq:conj1} is violated under an additional convexity assumption (see Conjecture~\ref{conjecture:ML}). We provide strong numerical evidence for the veracity of this conjecture. We refer to Section~\ref{sec:application} for rigorous results and detailed discussions on this example.


\begin{Example}[Dynamic SBM]\label{eg:dynamic SBM}
    This example is popular in modeling  networks observed over time \cite{matias2017statistical}. Setting $L_1=0$ and $\mathcal{Z}=\emptyset$,  $(x^{(1)},\ldots,x^{(L)})$ forms a time-homogeneous Markov chain evolving from $p_\Dyn(x^{(1)}=\pm 1)=\frac{1}{2}$ according to
    \begin{equation*}
       p_\Dyn(x^{(l+1)}=x^{(l)})=1-\rho,\quad p_\Dyn(x^{(l)}=-x^{(l-1)})=\rho,\quad l\in[L-1].
    \end{equation*}
    We can allow $a^{(l)}=b^{(l)}=0$---this corresponds to the setting where the $l^{th}$ layer is invisible to the statistician. See Figure~\ref{fig:Dyn illustration} for a graphical illustration of $p_\Dyn$.
\end{Example}

As in the previous example, we set $d^{(l)}= (a^{(l)}+ b^{(l)})/2$ and $\lambda_n^{(l)}= (a^{(l)}- b^{(l)})^2/(4 d^{(l)}(1-d^{(l)}/n)$. We assume $\lambda^{(l)}_n \to \lambda^{(l)}$ as $n \to \infty$. 
We conjecture the following information theoretic threshold for community detection in the Dynamic SBM in the special case $\lambda^{(1)}=\cdots= \lambda^{(L)}$. 

\begin{Conjecture}
    Set $d = \min_{1\leq l \leq L} \frac{a^{(l)}+ b^{(l)}}{2}$ and assume that $d \to \infty$. In addition, assume that $\lambda^{(1)} = \cdots = \lambda^{(L)}=:\lambda$. Fix $0\leq \rho < \frac{1}{2}$. In this setting, weak recovery is possible if and only if 
    \begin{equation}
        \lambda > \frac{1-2(1-2\rho)^2\cos\theta_\ast+(1-2\rho)^4}{1-(1-2\rho)^4}, \label{eq:conj2}
    \end{equation}
    where $\theta_\ast\in(0, \pi)$ is the minimum solution of equation
    \begin{equation*}
        0 = \sin[(L+1)\theta_\ast]-2(1-2\rho)^2\sin[ L\theta_\ast]+(1-2\rho)^4\sin[(L-1)\theta_\ast].
    \end{equation*}
\end{Conjecture}

We provide a partial proof of this conjecture in this work. We refer to Section~\ref{sec:application} for an in-depth discussion of the our results. 

Our final example focuses on community detection in partially labeled Stochastic Block Models (SBM). Semi-supervised community detection has been investigated across diverse communities e.g. statistics \cite{cai2020weighted,jiang2022semi,leng2019semi}, information theory \cite{kanade2016global}, statistical physics \cite{allahverdyan2010community,saade2018fast,zhang2014phase} etc. We focus on a SBM with unbalanced partially observed labels. 

\begin{Example}[SBM with Unbalanced Partially Observed Labels]\label{eg:semi sbm}
    Fix $\varepsilon_+,\varepsilon_-\in[0,1]$. Then $p(x=1)=p(x=-1)=\frac{1}{2}$, with $z$ taking values in $\{1,\ast,-1\}$,
    \begin{equation}
        \begin{aligned}
        p(z=1|x=1)&=1-p(z=\ast|x=1)=\varepsilon_+,\\
        p(z=-1|x=-1)&=1-p(z=\ast|x=-1)=\varepsilon_-.
    \end{aligned}
    \end{equation}
\end{Example}
The variable $x$ represents the latent community assignment of a vertex. The variable $z$ captures the side information regarding the latent assignment---if $\{z=*\}$, the latent label is unobserved. Otherwise, $\{z=x\}$, and the community assignment $x$ is observed by the scientist. Parameters $\varepsilon_+,\varepsilon_-$ govern the proportion of labeled vertices within each community. Prior work in this area focuses mainly on the setting $\varepsilon_+=\varepsilon_-$ i.e., vertices of both communities are equally likely to be revealed to the observer.
In contrast, our model allows $\varepsilon_+\neq\varepsilon_-$, which control the probability of revealed labels in the two groups respectively. It is particularly interesting to study the unbalanced setting $\varepsilon_+\ll\varepsilon_-$ or $\varepsilon_+\gg\varepsilon_-$---in either case, the labeled vertices belong  predominantly to one of the communities. 
We characterize the limiting mutual information between the data and the latent signal in this model in Section \ref{sec:application}. In addition, we introduce a suitable algorithm based on Approximate Message passing, and discuss its optimality from an information theoretic perspective.


\vspace*{2mm}
\noindent\textbf{A Spiked Matrix Surrogate:} Denote $d^{(l)}:=\frac{a^{(l)}+b^{(l)}}{2}$ as the average degree for each layer. Using the generative model \eqref{eq:original graph model}, the adjacency matrices maybe decomposed as 
%
\begin{equation*}
    \frac{\mathbf{G}^{(l)}-d^{(l)}/n}{\sqrt{d^{(l)}(1-d^{(l)}/n)/n}}=\sqrt{\frac{\lambda^{(l)}_n}{n}}\bX^{(l)}\bX^{(l)\top}+\mathbf{H}^{(l)},
\end{equation*}
with $\lambda^{(l)}_n:=\frac{(a^{(l)}-b^{(l)})^2}{4d^{(l)}(1-d^{(l)}/n)}$ being the effective SNR per layer. Conditioned on each $\bX^{(l)}$, the noise matrix $\mathbf{H}^{(l)}$ has independent off-diagonal entries satisfying that for each $1\le i<j\le n$,
\begin{align}
    &\quad\E[H_{i,j}^{(l)}]=0, \qquad|H_{i,j}^{(l)}| \leq \frac{1}{\sqrt{d^{(l)}(1-d^{(l)}/n)/n}},\notag\\
    &\E[(H_{i,j}^{(l)})^2]\in\left\{\frac{a^{(l)}(1-a^{(l)}/n)}{d^{(l)}(1-d^{(l)}/n)},\frac{b^{(l)}(1-b^{(l)}/n)}{d^{(l)}(1-d^{(l)}/n)}\right\}.
\end{align}
The parameters $(d^{(l)},\lambdal_n)$ provide a more useful parametrization of $(a^{(l)},b^{(l)})$ for our analysis. 
Throughout this paper, we focus on an asymptotic setting where each $\lambdal_n\rightarrow\lambdal$ has a finite limit as $n\rightarrow\infty$. 
In our analysis, we focus initially on the following set of gaussian spike matrices
\begin{equation}\label{eq:spiked matrix each layer}
    \mathbf{A}^{(l)}=\sqrt{\frac{\lambda^{(l)}}{n}}\bX^{(l)}\bX^{(l)\top}+\mathbf{W}^{(l)},
\end{equation}
with $\mathbf{W}^{(l)}$ being sampled independently according to the Gaussian Orthogonal Ensemble, $W_{i,j}^{(l)}\sim\mathcal{N}(0,1+\mathbbm{1}\{i=j\})$ for $i<j$. We will show that this gaussian model provides a tractable surrogate for the inference problem introduced above, provided the average degrees $d^{(l)}$ are sufficiently large. Our approach will generalize similar results in the context of the traditional community detection problem \cite{deshpande2017asymptotic,lelarge2019fundamental,wang2022universality}. 


\vspace*{2mm}
\noindent\textbf{Effective Scalar Channel:} In both \eqref{eq:original graph model} and \eqref{eq:spiked matrix each layer}, each observed entry $G^{(l)}_{i,j}$ (or $A^{(l)}_{i,j}$) relates to two independent realizations $(X_i,Y_i,Z_i)$ and $(X_j,Y_j,Z_j)$. 
We will show that the pairwise observation model is asymptotically equivalent to the following scalar observation channel: 
\begin{equation}\label{eq:scalar channel}
    (X,Y,Z)\sim p(\cdot),\quad A^{(l)\prime}=\sqrt{\lambda^{(l)} q_l}X^{(l)}+W^{(l)\prime},
\end{equation}
for $l\in[L]$ and $W^{(l)\prime}\overset{\text{i.i.d.}}{\sim}\mathcal{N}(0,1)$, with a specific $\boldsymbol{q}=(q_1,\ldots,q_L)\ge0$. Within this channel, one observes $\left(A^\prime,Z\right)$ and wants to infer $(X,Y)$. As a multivariate white noise additive channel, it has been studied in-depth in information theory \cite{guo2005mutual}. See Figure~\ref{fig:prior illustration} for a visualization of the effective scalar channels under different inter-layer priors specified in Example~\ref{eg:multilayer SBM} and~\ref{eg:dynamic SBM}.

\subsection{Main Results}
For the spiked gaussian matrix model \eqref{eq:spiked matrix each layer}, the posterior distribution is given as  
\begin{equation*}
    p(\bx,\by|\bZ,\bA)\propto \prod_{i=1}^np(x_i,y_i|Z_i)\cdot \exp\left[-\frac{1}{2}\sum_{l=1}^L\sum_{i<j}\left(A^{(l)}_{i,j}-\sqrt{\frac{\lambdal}{n}}x_i^{(l)}x_j^{(l)}\right)^2\right],
\end{equation*}
for $\bx\in\{\pm1\}^{n\times L},\by\in\{\pm1\}^{n\times L_1}$. In physics parlance, the posterior distribution corresponds to a  Boltzman Gibbs distribution $p(\bx,\by|\bZ,\bA)\propto\prod_{i}p(x_i,y_i|Z_i)\cdot e^{H(\bx)}$ with Hamiltonian defined as
\begin{equation*}
    H(\bx)=\sum_{l=1}^L\sum_{i<j}\sqrt{\frac{\lambda^{(l)}}{n}}A_{i,j}^{(l)}x_i^{(l)}x_j^{(l)}=\sum_{l=1}^L\sum_{i<j}\sqrt{\frac{\lambda^{(l)}}{n}}W_{i,j}^{(l)}x_i^{(l)}x_j^{(l)}+\frac{\lambda^{(l)}}{n}X_i^{(l)}X_j^{(l)}x_i^{(l)}x_j^{(l)}.
\end{equation*}
Set $Z_n(\blam)=\sum_{\bx,\by}\prod_{i}p(x_i,y_i|Z_i)\cdot e^{H(\bx)}$ to be the normalizing constant of this Gibbs distribution.
The free energy of this matrix model is formally defined as
\begin{equation}\label{eq:free energy def}
    F_n(\blam)=\frac{1}{n}\mathbb{E}\log Z_n(\blam).
\end{equation}
For convenience, we also define the free energy of the scalar channel \eqref{eq:scalar channel} as
\begin{equation}\label{eq:scalar channel free energy}
    \mathcal{F}(\boldsymbol{\lambda},\mathbf{q})=\mathbb{E}\log \left(\sum_{x,y}p(x,y|Z)e^{\sum_{l=1}^L\lambda^{(l)}q_lX^{(l)}x^{(l)}+\sqrt{\lambda^{(l)}q_l}W^{(l)\prime}x^{(l)}}\right).
\end{equation}
We elaborate on the explicit structure of the scalar free energy $\mathcal{F}$ for the examples \ref{eg:multilayer SBM}, \ref{eg:dynamic SBM} and \ref{eg:semi sbm} in Section~\ref{sec:application}. 
Our first result derives the limit of the free energy $F_n(\boldsymbol{\lambda})$.

\begin{Proposition}\label{prop:RS prediction on free energy}
    For any $\boldsymbol{\lambda}\in[0,\infty)^L$
    \begin{equation}\label{eq:RS prediction on free energy}
        \lim_{n\rightarrow\infty}F_n(\boldsymbol{\lambda})=\frac{\sum_{l=1}^{L}\lambdal}{4}+\sup_{\mathbf{q}\ge 0}\left[\mathcal{F}(\boldsymbol{\lambda},\mathbf{q})-\sum_{l=1}^{L}\frac{\lambda^{(l)}(q_l^2+2q_l)}{4}\right].
    \end{equation}
\end{Proposition}
\noindent 
As a direct corollary of this proposition, we derive an asymptotic formula for the normalized mutual information and  translate it to the estimation error of co-membership matrices $\bX^{(l)}\bX^{(l)\top}$ or $\bY^{(l_1)}\bY^{(l_1)\top}$. The error is measured by MMSE (minimal mean squared squared error) defined formally as
\begin{equation}\label{eq:def mmse}
    \MMSE\left(\bX^{(l)}\bX^{(l)\top};\bA,\bZ\right)=\frac{2}{n(n-1)}\sum_{i<j}\mathbb{E}\left(X_i^{(l)}X_j^{(l)}-\mathbb{E}\left[X_i^{(l)}X_j^{(l)}|\bA,\bZ\right]\right)^2.
\end{equation}
We define $\MMSE\left(\bY^{(l_1)}\bY^{(l_1)\top};\bA,\bZ\right)$ analogously.


\begin{Theorem}\label{thm:asymptotic MI}
Let $\bA$ be generated from the gaussian spiked matrix model \eqref{eq:spiked matrix each layer}.
\begin{itemize}
    \item[(i)] The normalized mutual information between latent variables $(\bX, \bY)$ and observed variables $(\bA, \bZ)$ has a finite limit,
    \begin{equation}\label{eq:limiting MI}
       \lim_{n \to \infty} \frac{1}{n}I(\bX,\bY;\bA,\bZ) =  i_p(x,y;z)+\frac{\sum_l\lambdal}{4}-\sup_{\mathbf{q}\ge 0}\left[\mathcal{F}(\boldsymbol{\lambda},\mathbf{q})-\sum_{l=1}^{L}\frac{\lambda^{(l)}(q_l^2+2q_l)}{4}\right],
    \end{equation}
    where $i_p(x,y;z)$ is the mutual information between $(x,y)$ and $z$ under prior $p(\cdot)$.
    \item[(ii)] There exists a set $D \subset (0, +\infty)^L$ of Lebesgue measure zero such that for all $\boldsymbol{\lambda} \in (0, +\infty)^L\setminus D$
    there exists a unique maximizer $\mathbf{q}^\ast=\left(q^{\ast}_1,\ldots,q^\ast_L\right)$ of the RHS in \eqref{eq:limiting MI}. Further, this maximizer satisfies
    \begin{equation}\label{eq:limiting mmse}
        \lim_{n\rightarrow\infty}\MMSE\left(\bX^{(l)}\bX^{(l)\top};\bA,\bZ\right)= 1-\left(q^{\ast}_{l}\right)^2,\,\,\,\,\,\forall\,\, l\in[L].
    \end{equation}
    \item[(iii)] For any $\boldsymbol{\lambda}\in(0,\infty)^L\setminus D$ and $l\in[L_1]$, the optimal estimation error of $\bY^{(l_1)}$ converges to
    \begin{equation}\label{eq:limiting mmse implicit}
        \lim_{n\rightarrow\infty}\MMSE\left(\bY^{(l_1)}\bY^{(l_1)\top};\bA,\bZ\right)= 1-\mathbb{E}\left[Y^{(l)}\mathbb{E}\left(Y^{(l_1)}|Z,\sqrt{\boldsymbol{\lambda}\odot\mathbf{q}^\ast}\odot X+W^\prime\right)\right]^2.
    \end{equation}
\end{itemize}
\end{Theorem}


Our next result establishes universality for the limiting mutual information. Specifically, we establish that the limiting mutual information in the gaussian spike matrix model \eqref{eq:spiked matrix each layer} equals the limiting mutual information in the random graph model \eqref{eq:original graph model}, provided the average degrees diverge to infinity. 

\begin{Proposition}\label{prop:Universality MI qualitative}
Fix $\boldsymbol{\lambda}\in[0,\infty)^L$. Under the asymptotics that (i) $\lambdal_n\rightarrow\lambdal $ for any $l\in[L]$ and (ii) $\min_l d^{(l)}(1-d^{(l)}/n)\rightarrow\infty$ as $n\rightarrow\infty$,
    \begin{equation*}
        \lim_{n\rightarrow\infty}\frac{1}{n}I(\bX,\bY;\mathbf{G},\bZ)=\lim_{n\rightarrow\infty}\frac{1}{n}I(\bX,\bY;\bA,\bZ),
    \end{equation*}
    where the RHS is given in Theorem~\ref{thm:asymptotic MI}. 
\end{Proposition}

We proceed to investigate Bayes optimal estimation of community memberships $(\bX,\bY)$ from $\bZ$ and the observed networks $\bG$. A core question is whether or not $\bG$ contributes to the recovery of $\bX$ or $\bY$. If not, a dummy constant estimator $\left\{\E[X_i|Z_i]\right\}_{i\in[n]}$ only based on $\bZ$ would achieve the following matrix mean squared error,
\begin{align}
    \DMSE\left(\bX^{(l)}\bX^{(l)\top};\bZ\right)&=\frac{2}{n(n-1)}\sum_{i<j}\mathbb{E}\left(X_i^{(l)}X_j^{(l)}-\mathbb{E}\left[X_i^{(l)}X_j^{(l)}|\bZ\right]\right)^2\notag\\
    &=1-\left\{\E\left[X^{(l)}\E\left(X^{(l)}|Z\right)\right]\right\}^2:=\DMSE_l(p)\label{eq:dummy mmse in general model}
\end{align}
which is fixed over different $n$ and only depends on joint prior $p(x,y,z)$. Mathematically, we seek to determine whether inequalities
\begin{equation*}
    \MMSE\left(\bX^{(l)}\bX^{(l)\top};\bA,\bZ\right)\le \DMSE_l(p),\quad \MMSE\left(\bX^{(l)}\bX^{(l)\top};\bG,\bZ\right)\le \DMSE_l(p),
\end{equation*}
are asymptotically tight.  
The next proposition establishes that the sharp threshold in terms of $\blam$ is universal between the graph model~\eqref{eq:original graph model} and its spiked matrix counterpart~\eqref{eq:spiked matrix each layer}.
\begin{Proposition}\label{prop:universality mmse}
    Fix $\boldsymbol{\lambda}\in[0,\infty)^L$. Under the asymptotics that (i) $\blam_n\rightarrow\blam$ and (ii) $\min_l d^{(l)}(1-d^{(l)}/n)\rightarrow\infty$ as $n\rightarrow\infty$, for any $l\in[L]$, 
    \begin{equation*}
        \lim_{n\rightarrow\infty}\MMSE\left(\bX^{(l)}\bX^{(l)\top};\bA,\bZ\right)<\DMSE_l(p),
    \end{equation*}
    if and only if 
    \begin{equation*}
        \lim_{n\rightarrow\infty}\MMSE\left(\bX^{(l)}\bX^{(l)\top};\bG,\bZ\right)<\DMSE_l(p).
    \end{equation*}
    The same conclusion holds if we replace $\bX^{(l)}$ with any $\bY^{(l_1)}$.
\end{Proposition}

\begin{Remark}\label{rmk:def of weak recovery}
    The dummy mean squared error $\DMSE_l(p)$ measures the contribution of per-vertex side information $\bZ$. In the context of Example~\ref{eg:dynamic SBM} and \ref{eg:multilayer SBM} where $\bZ$ is absent and the prior on $X^{(l)}$ is uniform, it holds that $\DMSE_l(p)=1$.
\end{Remark}

\subsection{Message Passing Algorithm}

\begin{algorithm}[tb]
\caption{Coupled Approximate Message Passing}\label{alg:coupled AMP}
\begin{algorithmic}
\Require a sequence of non-linear denoisers $\left\{\mathcal{E}_t:\mathbb{R}^{L+L_2}\rightarrow\mathbb{R}^{L}\right\}_{t\ge 0}$, initialization $\mathbf{m}^0\in\mathbb{R}^{n\times L}$.
\For{$t \ge 0$}
    \State Update each coordinate $l\in[L]$ and node $i\in[n]$ by
    \begin{equation}\label{eq:coupled AMP update}
        m_{i,l}^{t+1}=\sum_{k=1}^n\sqrt{\frac{\lambda^{(l)}}{n}}A^{(l)}_{i,k}\mathcal{E}^{(l)}_t(m^t_{k},Z_k)-\lambda^{(l)}\mathsf{d}_t^{(l)}\mathcal{E}^{(l)}_{t-1}(m^{t-1}_{i},Z_i),
    \end{equation}
    \State where the Onsager term per layer $l\in[L]$ is given by
    \begin{equation}\label{eq:coupled AMP Onsager}
        \mathsf{d}_t^{(l)}=\frac{1}{n}\sum_{k=1}^n\partial_l\mathcal{E}^{(l)}_t(m^t_{k},Z_k).
    \end{equation}
\EndFor
\State Stop at $t_*\in\mathbb{N}$, and output $\widehat{\bx}^{t_*}=\mathcal{E}_{t_*}\left(\mathbf{m}^{t_*},\bZ\right)$.
\end{algorithmic}
\end{algorithm}

In this section, we turn to a study of algorithms for signal recovery in the general inference problem \eqref{eq:original graph model}. We devise signal recovery algorithms based on Approximate Message Passing (AMP). AMP algorithms are a class of efficient iterative algorithms which involve matrix-vector products, followed by element-wise non-linearities. These algorithms were originally introduced in the context of compressed sensing \cite{donoho2009message}, and have been applied widely over the past decade in diverse inference problems \cite{bayati2011dynamics,lesieur2017constrained,richard2014statistical,li2023approximate}. We refer the interested reader to the recent surveys \cite{feng2022unifying,montanari2022short} for an in-depth discussion of these algorithms and their applications. We present the relevant class of AMP algorithms in Algorithm~\ref{alg:coupled AMP}, and present a heuristic derivation of this algorithm in Section~\ref{subsec:heurestic derivation}. The main power of AMP algorithms arise from their flexibility; indeed, the scientist can choose the sequence of non-linear functions $\{ \mathcal{E}_t: \mathbb{R}^{L} \to \mathbb{R}^{L}\}$ according to the specific application, subject to minor restrictions on their regularity. In our discussion below, we first introduce the AMP algorithm for the gaussian spiked matrix model \eqref{eq:spiked matrix each layer}; subsequently, we establish that the performance of AMP algorithms on the original graph problem \eqref{eq:original graph model} match that on the gaussian spiked model \eqref{eq:spiked matrix each layer}, provided the average degree is large enough. Our algorithmic universality analysis slightly generalizes the framework introduced recently in \cite{wang2022universality}. Finally, we establish that the proposed AMP algorithm is provably optimal (from an information theoretic perspective) under very general conditions.


The main attraction of AMP algorithms stems from their analytic tractability. Specifically, their performance can be rigorously described using a low-dimensional scalar recursion, referred to as state evolution \cite{donoho2009message}. We present the state evolution for Algorithm \ref{alg:coupled AMP}. 
\begin{Definition}[State evolution iterates]\label{def:SE iterates}
    Suppose  $\mathcal{E}_t$ is Lipschitz continuous for $t \geq 0$. The state evolution iterates for Algorithm~\ref{alg:coupled AMP} corresponds to the bivariate sequence $\mu^t,\kappa^t\in\mathbb{R}^L$, indexed by $t\ge 0$. The sequence is recursively defined as follows
    \begin{align}
        \mu^{t+1}_l &= \E\left[ X^{(l)} \mathcal{E}_t^{(l)}\left(\widetilde{m}^t,Z\right)\right],\label{eq:coupled-AMP-SE-mu}\\
        \kappa^{t+1}_l &= \E\left[ \mathcal{E}_t^{(l)}\left(\widetilde{m}^t,Z\right)^2 \right],\label{eq:coupled-AMP-SE-kappa}
    \end{align}
    where $(X,Y,Z) \sim p(\cdot)$ and $\widetilde{m}^t=\left(\widetilde{m}^t_1,\ldots,\widetilde{m}^t_L\right)$ is an $L$-dim variable with $\widetilde{m}^t_l|X\indsim\mathcal{N}\left(\lambdal\mu^t_lX^{(l)},\lambdal\kappa^t_l\right)$.
\end{Definition}
\begin{Theorem}[State evolution]\label{thm:coupled-AMP-SE}
    For any pseudo Lipschitz test function $\psi:\mathbb{R}^{L}\times\{\pm 1\}^{L+L_2}\rightarrow\mathbb{R},(m,X,Z)\rightarrow\psi(m,X,Z)$ as  $n\rightarrow\infty$,
    \begin{equation*}
        \frac{1}{n}\sum_{i=1}^n\psi(m^t_i,X_i,Z_i)\overset{\text{P}}{\rightarrow}\mathbb{E}\left[\psi(\widetilde{m}^t,X,Z)\right]. 
    \end{equation*}
\end{Theorem}
\noindent As evident from the SE recursions \eqref{eq:coupled-AMP-SE-mu} and \eqref{eq:coupled-AMP-SE-kappa}, the sequence $\mu^t$ measures the correlation between $\mathbf{m}^t$ and $\bX$ while the sequence $\kappa^t$ measures the conditional variance. Therefore, the Bayes optimal choice of $\mathcal{E}_t$ is given as 
\begin{equation*}
    \mathcal{E}_t\left(\widetilde{m}^t,Z\right)=\mathbb{E}\left[X|\widetilde{m}^t,Z\right]=\frac{\sum_{x} p(x|Z)\exp\left(\sum_{l=1}^L\widetilde{m}^t_l\mu^t_l x^{(l)}/\kappa_l^t\right)x}{\sum_{x}p(x|Z)\exp\left(\sum_{l=1}^L\widetilde{m}^t_l\mu^t_l x^{(l)}/\kappa_l^t\right)},
\end{equation*}
where the joint distribution of $(X,Y,Z,\widetilde{m}^t)$ is given in Definition~\ref{def:SE iterates}. In this case, we further find $\mu^{t+1}_l=\kappa^{t+1}_l$ for any $t,l$ using Nishimori identities.
In conclusion, we can use a single state evolution iterate to describe the asymptotic distribution of this \textit{Bayes optimal} coupled AMP algorithm,
\begin{equation*}
    q^{t+1}_l=\E\left[ \mathbb{E}\left[X^{(l)}\mid\widetilde{m}^t,Z\right]^2 \right]=\E\left[X^{(l)}\frac{\sum_{x} p(x|Z)\exp\left(\sum_{l=1}^L\widetilde{m}^t_l x^{(l)}\right)x^{(l)}}{\sum_{x}p(x|Z)\exp\left(\sum_{l=1}^L\widetilde{m}^t_l x^{(l)}\right)}\right],
\end{equation*}
where $(X,Y,Z,\widetilde{m}^t)$ is jointly drawn from the effective scalar channel \eqref{eq:scalar channel} with profile $\mathbf{q}^t$.
\begin{Corollary}
    For Bayes optimal AMP and any  $t\geq 0$, $\widehat{\bx}^t=\mathbb{E}\left[X|\mathbf{m}^t,Z\right]$ has the following asymptotic mean squared error,
    \begin{equation*}
        \MSE_\AMP(t,n,\boldsymbol{\lambda}):=\frac{2}{n(n-1)}\sum_{i<j}\mathbb{E}\left(X_i^{(l)}X_j^{(l)}-\widehat{x}^{t(l)}_i\widehat{x}^{t(l)}_j\right)^2\nconverge 1-\left(q^{t}_{l}\right)^2.
    \end{equation*}
    Moreover, if $\mathbf{q}^t$ converges to the unique maximizer $\mathbf{q}^\ast$ in \eqref{eq:RS prediction on free energy} as $t\rightarrow\infty$, the coupled AMP algorithm is Bayes optimal.
\end{Corollary}
\noindent Combined with Theorem~\ref{thm:asymptotic MI}, in some cases, Algorithm~\ref{alg:coupled AMP} provably achieves minimal MSE for the spiked matrices model \eqref{eq:spiked matrix each layer} as the number of iterations $t \to \infty$. 

\paragraph{Algorithmic Universality.}
We translate the AMP algorithm and associated guarantees from the gaussian spiked model \eqref{eq:spiked matrix each layer} to the original random graph model \eqref{eq:original graph model} in this section.  
Upon appropriate rescaling of the  adjacency matrices,
\begin{equation}\label{eq:graph low rank component}
    \bar{\mathbf{G}}^{(l)}:=\frac{\mathbf{G}^{(l)}-d^{(l)}/n}{\sqrt{d^{(l)}(1-d^{(l)}/n)/n}}=\sqrt{\frac{\lambda^{(l)}}{n}}\bX^{(l)}\bX^{(l)\top}+\mathbf{H}^{(l)}.
\end{equation}
One can implement Algorithm~\ref{alg:coupled AMP} with the matrices $\bar{\bG}$ in place of the spiked gaussian matrix $\bA$. 
Denote the recursive output as $\left(\bar{\mathbf{m}}^1,\ldots,\bar{\mathbf{m}}^t\right)$. The next proposition asserts that the asymptotic distribution of $\bar{\mathbf{m}}$ is still characterized by the prescribed state evolution iterates in Definition~\ref{def:SE iterates}, as long as average degrees $d^{(l)}$ are $\Omega(\log n)$.

\begin{Proposition}\label{prop:algorithmic universality}
Suppose that $\min_l d^{(l)}(1-d^{(l)}/n)\ge C\log n$ for some constant $C>0$ independent of $n$. Further let the non-linear denoisers used in Algorithm~\ref{alg:coupled AMP} be Lipschitz with polynomial growth condition i.e., for some integer $p\in\mathbb{N}$,
\begin{equation}\label{eq:polygrowth}
    \|\mathcal{E}_t(m_1,\ldots,m_L,z)\| \leq C(1+\|(m_1,\ldots,x_L,z)\|_2^p) \text{ for a constant } C>0.
\end{equation}
Consequently, for any pseudo Lipschitz test function $\psi:\mathbb{R}^{L}\times\{\pm 1\}^{L+L_2}\rightarrow\mathbb{R},(m,X,Z)\rightarrow\psi(m,X,Z)$ as $n\rightarrow\infty$,
\begin{equation}\label{eq:algorithmic universality}
    \frac{1}{n}\sum_{i=1}^n\psi(m^t_i,X_i,Z_i)-\frac{1}{n}\sum_{i=1}^n\psi(\bar{m}^t_i,X_i,Z_i)\overset{\text{P}}{\rightarrow}0.
\end{equation}
\end{Proposition}

\paragraph{Benefits and Novelty of Aggregating Information.}
The information-theoretic threshold for recovering community memberships from a single stochastic block model has been characterized rigorously in  \cite{mossel2012stochastic,mossel2018proof,bordenave2015non}. Based on a single $\bA^{(l)}$, one requires $\lambdal>1$ for weak recovery. Thus based on the marginal information, weak recovery is possible if and only if $\max_l \lambdal >1$. 

Our Algorithm~\ref{alg:coupled AMP} is novel in that it \textit{aggregates all networks (spiked matrices) available from multiple sources} in a non trivial way tailored to the \textit{prior knowledge about relationships among observed layers.} To be precise, information is jointly extracted via recursive application of prior-specific denoisers $\mathcal{E}$ onto $\left\{\mathbf{A}^{(l)},l\in[L]\right\}$. As shown in Section~\ref{sec:application} for specific examples, Algorithm~\ref{alg:coupled AMP} can recover node memberships even when $\max_{l}\lambdal\le 1$. This algorithm is also Bayes optimal in diverse settings---we verify this for our applications on a case by case basis. 
In the following, we make comparisons to existing related AMP algorithms.
\begin{itemize}
    \item[(i)] \citet{ma2023community} proposed an AMP algorithm designed for a homogeneous multilayer stochastic block model (a special case of Example~\ref{eg:multilayer SBM} with $\rho=0$), in which $\bX^{(l)}=\bY$ always holds. Their algorithm would sum up all observed networks and then run a single-layer AMP. Algorithm~\ref{alg:coupled AMP} is more flexible in the sense that it allows across-layer node memberships to differ. If one observes additional contextual side information in the form of gaussian mixtures as in \cite{deshpande2018contextual,lu2023contextual,ma2023community}, one might also incorporate an orchestrated version of Algorithm~\ref{alg:coupled AMP} to exploit this side information, in the same spirit as \cite{ma2023community}.
    
    \item[(ii)] In \cite{gerbelot2021graph}, graph-based approximate message passing iterations are proposed to unify recent advances \cite{berthier2020state,manoel2017multi,aubin2018committee} in applying AMP to complicated hierarchical problems. In this class of iterations, multiple observations are connected via rectangular sensing matrices. It is noteworthy that Algorithm~\ref{alg:coupled AMP} is not included in this class, since we only get to observe symmetric sensing matrices and the latent variables are connected via the prior. As detailed in Remark~\ref{remark:comparison to gerbelot-berthier}, we are able to borrow an intermediate result of \cite{gerbelot2021graph} to simplify the proof of Theorem~\ref{thm:coupled-AMP-SE}.
    
    \item[(iii)] We note that a very similar AMP algorithm has been proposed in the context of the Matrix Tensor product model in \cite{rossetti2023approximate}. However, the authors focus explicitly on the gaussian spike matrix setting in their work. In contrast, we are motivated by inference problems on sparse networks, and  first establish universality for these algorithms. Finally, our work yields the sharp information theoretic thresholds for these problems, which has not been explored in prior work.  
    

\end{itemize}

\subsection{Related Literature}
\label{sec:lit_review}

We take this opportunity to survey the current state-of-the-art on our applications of interest. 

\paragraph{Multilayer SBM:} Community detection for multilayer networks has attracted considerable attention recently. Most existing works model the observed networks using the multilayer Stochastic Block Model (SBM). Most early works in the area assume that the community assignments are constant across layers (see e.g. \cite{bhattacharyya2020general,lei2020consistent,ma2023community} and references therein). As argued in \cite{chen2022global}, this assumption is unrealistic in practice. In many applications, one assumes that the layerwise community assignments are unequal, but correlated. \cite{chen2022global} introduces a tractable model for this setting, and analyzes the performance of a two-step algorithm which employs spectral clustering followed by a MAP refinement. We study the same model as \cite{chen2022global}, but as remarked earlier, their theoretical analysis focuses on the high SNR setting where exact recovery is possible. In contrast, we focus on weak recovery, and study this problem in the low SNR regime. Our analysis is closest in spirit to the one in \cite{ma2023community}. However, we go beyond the homogeneous community assignment critical in this work. Finally, the recent work \cite{lei2023computational} studies testing for  latent community structures in multilayer networks.


\paragraph{Dynamical SBM:} The dynamical Stochastic Block Model (dynamical SBM) studied here was introduced in \cite{yang2011detecting}. 
Several algorithms for community detection have been developed in this context---\cite{matias2017statistical, longepierre2019consistency} introduce a  variational EM algorithm for community detection in this setting, while  \cite{keriven2022sparse, pensky2019spectral} study spectral clustering for community recovery in this model. Finally  \cite{bhattacharjee2020change} studies change point detection in the setting of dynamic SBMs. In contrast to these existing works, we focus on the sharp information theoretic threshold for weak recovery in the dynamic SBM under a low SNR setting. 

\paragraph{Semi-supervised Community Detection:} In certain practical settings, the community assignment of some vertices might be available a  priori. Having access to the true community assignment of a few vertices can be valuable in the context of community detection. For example, in the unlabeled balanced community recovery problem, one can only hope to recover the latent communities up to a global permutation. However, once some vertices are labeled, one can actually recover the true labeling. In addition, the presence of the partial labeling shifts the threshold for weak recovery, and promises to improve the performance of local algorithms. This has prompted a thorough investigation of semi-supervised community detection (see e.g. \cite{cai2020weighted,jiang2022semi,leng2019semi,kanade2016global,allahverdyan2010community,saade2018fast,zhang2014phase} and references therein). Most prior work in this direction assumes that the true labels of $\epsilon$-fraction of the vertices are observed in addition to the network. \cite{cai2020weighted,saade2018fast,zhang2014phase} analyze the performance of local algorithms such as belief propagation in this setting. \cite{kanade2016global} formalizes several folklore conjectures in this area, and rigorously establishes that one can improve over the Kesten-Stigum threshold with additional partial labels in a setting with a large number of communities. Finally, we refer to \cite{jiang2022semi,leng2019semi} for recent progress in analyzing these questions on more realistic models of real networks.

\section{Implications for Representative Models}\label{sec:application}
We discuss the applications of our general results to specific applications of interest in this Section. To this end, we introduce the notion of weak recovery. Recall the inference problem \eqref{eq:original graph model}. We focus on a setting with no side information $\mathbf{Z}$, and assume that the prior on $\mathbf{X}^{(l)}$ is uniform under the prior $p$.

\begin{Definition}\label{def:weak recovery}
For $1 \leq l \leq L$, we say that $\bX^{(l)} \in \{ \pm 1\}^n$ is weakly recoverable if there exists an estimator  $\hat{\bX}^{(l)} := \hat{\bX}^{(l)}(\bG) \in \{\pm 1\}^n$ such that  
    \begin{equation*}
        \liminf_{n\to \infty} \mathbb{E}\left[\left|\frac{1}{n}\left\langle \widehat{\mathbf{X}}^{(l)},\mathbf{X}^{(l)}\right\rangle\right|\right]>0.
    \end{equation*}
    Weak recovery of $\bY^{(l)}$, $1\leq l \leq L_1$ is defined analogously. 
\end{Definition}
The latent signal is weakly recoverable if it can be reconstructed better than random guessing. Our definition of weak recovery is consistent with common notions of weak recovery (Definition 7 in \cite{abbe2017community}).
The following lemma relates weak recovery to the estimation of the $n\times n$ co-membership matrix $\mathbf{X}^{(l)}\mathbf{X}^{(l)\top}$.
\begin{Lemma}[Lemma 3.5, 3.6, \cite{deshpande2017asymptotic}]\label{lemma:MSE to weak recovery}
    Weak recovery of $\mathbf{X}^{(l)}$ is possible if and only if there exists an estimator $\widehat{\mathbf{C}}(\mathbf{G})\in\mathbb{R}^{n\times n}$ such that
    \begin{equation*}
        \MSE(\mathbf{X}^{(l)}\mathbf{X}^{(l)\top}, \widehat{\mathbf{C}})=\frac{2}{n(n-1)}\sum_{i<j}\mathbb{E}\left(X_i^{(l)}X_j^{(l)}-\widehat{C}_{i,j}\right)^2<1,
    \end{equation*}
    where RHS corresponds a random guess of the latent variables. 
\end{Lemma}
We are comparing $\MSE(\mathbf{X}^{(l)}\mathbf{X}^{(l)\top}, \widehat{\mathbf{C}})$ against $1$ as a random guess would yield a dummy mean squared error of $1$ i.e., the  prior variance of every $X_i^{(l)}X_j^{(l)}$. Since we adopt mean squared error, the best possible estimator would be posterior mean $\mathbb{E}\left[\bX^{(l)}\bX^{(l)\top}|\mathbf{G}\right]$, whose error is studied in \eqref{eq:limiting mmse}. So we turn to study the minimal mean squared error of $\bX^{(l)}\bX^{(l)\top}$ defined in \eqref{eq:def mmse}.

\begin{Remark}
    At the end of Section~\ref{subsec:general model setup}, we investigate whether 
    $\MMSE\left(\bX^{(l)}\bX^{(l)\top};\bG,\bZ\right)<\DMSE_l(p)$; informally, one seeks to understand whether a combination of the networks $\bG$ and side information $\bZ$ can strictly improve over the statistical performance of estimators derived solely based on the side information $\bZ$. In the context of multilayer (Example~\ref{eg:multilayer SBM}) and dynamic SBMs (Example~\ref{eg:dynamic SBM}) where no side information $\bZ$ is available and the prior is marginally uniform, it holds that $\DMSE_l(p)=1$. Thus this threshold matches the one for weak recovery introduced in Definition~\ref{def:weak recovery}. 
\end{Remark}


Theorem~\ref{thm:asymptotic MI} part (ii) implies that for  $\boldsymbol{\lambda} \notin D$ the question of weak recovery boils down to characterizing the behavior of the unique maximizer  $q^\ast_l$ \eqref{eq:limiting MI}. For notational compactness, we define  \eqref{eq:limiting MI} by
\begin{equation}\label{eq:opt task}
     G(\mathbf{q}) = \mathcal{F}(\boldsymbol{\lambda},\mathbf{q})-\sum_{l=1}^{L}\frac{\lambda^{(l)}(q_l^2+2q_l)}{4},
\end{equation}
where $\mathcal{F}(\boldsymbol{\lambda},\mathbf{q})$ is the free energy function defined in  \eqref{eq:scalar channel free energy}. For this section only, for any function $f(x,y)$ on $\{\pm 1\}^{L+L_1}$, define
\begin{equation*}
\left\langle f(x,y)\right\rangle_{\mathrm{sc}}=\mathbb{E}\left[f(X,Y)\mid Z,\sqrt{\boldsymbol{\lambda}\odot\mathbf{q}}\odot X+W^\prime\right]
\end{equation*}
as posterior mean after observing $\left(A^\prime, Z\right)$ from the scalar channel \eqref{eq:scalar channel}.
The next proposition collects some analytic results on the objective function $G$; the proof of these identities are deferred to Section~\ref{subsec:opt derivative}. Finally, we derive first order conditions for any stationary point of $G$, and establish that these conditions are satisfied by any maximizer, even at the boundary.  

    
\begin{Proposition}\label{prop:opt derivative}
\begin{itemize}
    \item[(i)] For any  $\boldsymbol{\lambda} \in [0, \infty)^L$, 
    \begin{align*}
        \partial_l G(\mathbf{q})&=\frac{\lambdal}{2}\mathbb{E}\left[\left\langle x^{(l)}\right\rangle_{\mathrm{sc}}^2\right]-\frac{\lambdal q_l}{2},\\
        \partial_{l_1}\partial_{l_2} G(\mathbf{q})&=\frac{\lambda^{(l_1)}\lambda^{(l_2)}}{2}\mathbb{E}\left\{\left(\left\langle x^{(l_1)}x^{(l_2)}\right\rangle_{\mathrm{sc}}-\left\langle x^{(l_1)}\right\rangle_{\mathrm{sc}}\left\langle x^{(l_2)}\right\rangle_{\mathrm{sc}}\right)^2\right\}-\frac{\lambda^{(l_1)}}{2}\mathbbm{1}\{l_1=l_2\}.
    \end{align*}
    Note that the RHS depends on $\mathbf{q}$ through the bracket $\langle\cdot\rangle_{\mathrm{sc}}$. 
    \item[(ii)] Define a mapping $T:[0,\infty)^L\rightarrow[0,\infty)^L$ with $l$-th coordinate given by
    \begin{equation}\label{eq:general se mapping}
        T_l(\boldsymbol{\gamma})=\mathbb{E}\left[\left\langle x^{(l)}\right\rangle_{\mathrm{sc}}^2\right]=\mathbb{E}\left[X^{(l)}\frac{\sum_{x,y}p(x,y|Z)x^{(l)}e^{\sum_{l=1}^L\gamma_lX^{(l)}x^{(l)}+\sqrt{\gamma_l}W^{(l)\prime}x^{(l)}}}{\sum_{x,y}p(x,y|Z)e^{\sum_{l=1}^L\gamma_lX^{(l)}x^{(l)}+\sqrt{\gamma_l}W^{(l)\prime}x^{(l)}}}\right].
    \end{equation}
    Any local maximizer $\mathbf{q}\in[0,\infty)^L$ of \eqref{eq:opt task} satisfies the first order stationary point condition $T(\boldsymbol{\lambda}\odot\mathbf{q})=\mathbf{q}$. 
    \end{itemize}
\end{Proposition}
\begin{Remark}
    We collect two immediate consequences of the above lemma.
    \begin{itemize}
        \item[1.] For $1\leq l 
        \leq L$, $T_l$ is non-decreasing in each argument as for any $l_1,l_2\in[L]$, there holds
        \begin{equation*}
            \partial_{l_2} T_{l_1}(\bgamma)=\mathbb{E}\left[\mathrm{Cov}(X^{(l_1)},X^{(l_2)}\big|\sqrt{\bgamma}\odot X+W^\prime)^2\right]\ge 0.
        \end{equation*}
        \item[2.] For any $\blam\in(0,+\infty)^L$, global maximizers $\mathbf{q}^\ast$ exist and are in a compact subset $[0,1]^L$. This follows as  $\partial_l G(\mathbf{q})\le \lambdal(1-q_l)/2<0$ as long as $q_l>1$. Consequently $G(\mathbf{q})< \max_{\mathbf{q}\in[0,1]^L} G(\mathbf{q})$ whenever $\max_l q_l>1$.
    \end{itemize}
\end{Remark}

If the prior $p(x,y,z)$ factorizes in special ways, it is possible that a certain layer $\bX^{(l_1)}$ is weakly recoverable while another layer $\bX^{(l_2)}$ is not. However, this situation does not happen under the applications of interest, as demonstrated in the following proposition. 


\begin{Proposition}\label{prop:all-in-all-out}
    Suppose there is no side information $Z$ and every pair of coordinates is correlated in the prior $p(x,y)$ and $\boldsymbol{\lambda}\notin D$. Then either all the layers can be weakly recovered, or none can be recovered weakly.
\end{Proposition}
Specifically, when $\rho \neq 1/2$, inhomogeneous multilayer SBM (Example~\ref{eg:multilayer SBM}) and dynamic SBM (Example~\ref{eg:dynamic SBM}) both satisfy the sufficient condition in Proposition~\ref{prop:all-in-all-out}.

\subsection{Applications to Inhomogeneous Multilayer SBM}

In this section, we specialize our general results to the inhomogenous multilayer SBM. The per-node prior $p$ for this example is specified in Example~\ref{eg:multilayer SBM} and admits a hierarchical structure  $p(x,y)=p(y)\prod_l p\left(x^{(l)}|y\right)$. As far as the authors know, most existing works on multilayer networks \cite{paul2020spectral, bhattacharyya2020general, kivela2014multilayer} assume the underlying community structure to be homogeneous across layers. Relatively few existing works \cite{chen2022global, valles2016multilayer} address possible label corruptions for each layer. The inhomogeneous model here is the same as that in \cite{chen2022global}, but they focus on high SNR settings where exact recovery is possible. They study community detection using spectral clustering, and analyze their method when the SNR parameters $\boldsymbol{\lambda}$ to diverge to infinity.  
Under this input prior, the effective scalar channel in \eqref{eq:scalar channel} reduces to
\begin{equation}\label{eq:scalar channel ML}
    Y\sim \text{Unif}\left\{\pm 1\right\}, \quad X^{(l)}\indsim \mathsf{BSC}_{\rho}(Y), \quad A^{(l)\prime}\indsim\sqrt{\lambdal q_l}X^{(l)}+\mathcal{N}\left(0,1\right),
\end{equation}
where observations can be interpreted as $L$ independent measurements of the composition of a binary symmetric channel and a Gaussian channel. We recall the general free energy functional \eqref{eq:scalar channel free energy} and note that in this special case, it reduces to
\begin{align*}
    \mathcal{F}_{\ML,\rho}(\boldsymbol{\lambda},\mathbf{q})&=\mathbb{E}\log\left\{\sum_{y=\pm 1}\frac{1}{2}\prod_{l=1}^L \left[\cosh\left(\lambdal q_l+\sqrt{\lambdal q_l}W^{(l)\prime}\right)\right.\right.\\
    &\qquad\qquad\qquad \left.\left.+\left(1-2\rho\right)yX^{(l)}\sinh\left(\lambdal q_l+\sqrt{\lambdal q_l}W^{(l)\prime}\right)\right]\right\},
\end{align*}
where the expectation is evaluated with respect to $W^{(l)\prime}\indsim\mathcal{N}(0,1)$.
Consequently, Theorem~\ref{thm:asymptotic MI} and Proposition~\ref{prop:Universality MI qualitative} immediately identify the asymptotic normalized mutual information between observed multilayer graph $\mathbf{G}$ and the underlying latent variables $(\mathbf{X},\mathbf{Y})$.
\begin{Theorem}\label{thm:multilayer joint MI}
    Under the asymptotics that (i) $\blam_n\rightarrow\blam$ and (ii) $\min_l d^{(l)}(1-d^{(l)}/n)\rightarrow\infty$ as $n\rightarrow\infty$, we have,
    \begin{align}
        &\quad\lim_{n\rightarrow\infty}\frac{1}{n}I(\bX,\bY;\mathbf{G})=\lim_{n\rightarrow\infty}\frac{1}{n}I(\bX,\bY;\bA)\notag\\
        &=\frac{\sum_l\lambdal}{4}-\sup_{\mathbf{q}\ge 0}\left[\mathcal{F}_{\ML,\rho}(\boldsymbol{\lambda},\mathbf{q})-\sum_{l=1}^{L}\frac{\lambda^{(l)}(q_l^2+2q_l)}{4}\right],\label{eq:multilayer joint MI}
    \end{align}
    where we recall $\mathbf{A}$ to be the  spiked gaussian matrices in \eqref{eq:spiked matrix each layer}. 
\end{Theorem}

\noindent In this specific setting, one might be interested in estimating the global membership $\bY$ or the local membership $\mathbf{X}^{(l)}$ of a specific layer. The next result answers this question. Our proof exploits the hierarchical structure in the prior $p$, depicted in Figure~\ref{fig:ML illustration}; we establish this result in 
Section~\ref{subsec:pf multilayer}.


\begin{Corollary}\label{cor:global and local MI multilayer}
    (i) For global estimation, define a free energy functional corresponding to the estimation of $\mathbf{X}^{(l)}$ from $\bG^{(l)}$ conditioned on $\bY$,
    \begin{align}
        \mathcal{F}^{(l)}_{\ML,\rho}(\lambdal,q_l)&=\mathbb{E}\log\left\{ \cosh\left(\lambdal q_l+\sqrt{\lambdal q_l}W^{(l)\prime}\right)\right.\notag\\
        &\qquad\qquad \left.+\left(1-2\rho\right)B\sinh\left(\lambdal q_l+\sqrt{\lambdal q_l}W^{(l)\prime}\right)\right\},\label{eq:free energy multilayer conditional global}
    \end{align}
    where the expectation is taken with respect to $B\sim2\mathrm{Bern}(\rho)-1$ and independent standard normal $W^{(l)\prime}$. Under the same asymptotics as Theorem~\ref{thm:multilayer joint MI}, we have, 
    \begin{align}
        \lim_{n\rightarrow\infty}\frac{1}{n}I(\bY;\mathbf{G})&=\sum_{l=1}^L\sup_{q_l\ge 0}\left[\mathcal{F}^{(l)}_{\ML,\rho}(\lambdal,q_l)-\frac{\lambda^{(l)}(q_l^2+2q_l)}{4}\right]\notag\\
        &\qquad -\sup_{\mathbf{q}\ge 0}\left[\mathcal{F}_{\ML,\rho}(\boldsymbol{\lambda},\mathbf{q})-\sum_{l=1}^{L}\frac{\lambda^{(l)}(q_l^2+2q_l)}{4}\right].\label{eq:global MI multilayer}
    \end{align}
    (ii) As for a specific individualized layer $l\in[L]$, consider another free energy functional corresponding to inferring $(\bX^{(-l)},\bY)$ from $\bG^{(-l)}$ conditioned on $\bX^{(l)}$,
    \begin{align}
    \bar{\mathcal{F}}^{(l)}_{\ML,\rho}(\boldsymbol{\lambda}^{(-l)},\mathbf{q}_{-l})&=\mathbb{E}\log\left\{\sum_{b=\pm 1}p(b)\prod_{l_1\neq l}^L \left[\cosh\left(\lambda^{(l_1)}q_{l_1}+\sqrt{\lambda^{(l_1)}q_{l_1}}W^{(l_1)\prime}\right)\right.\right.\notag\\
    &\qquad\qquad\qquad \left.\left.+\left(1-2\rho\right)bX^{(l)}X^{(l_1)}\sinh\left(\lambda^{(l_1)}q_{l_1}+\sqrt{\lambda^{(l_1)}q_{l_1}}W^{(l_1)\prime}\right)\right]\right\},\label{eq:free energy multilayer conditional local}
    \end{align}
    where we take $1-p(b=1)=p(b=-1)=\rho$. The expectation is taken with respect to $(X,Y)\sim p_{\ML}$ and independent standard normals $W^{(l)}$. Under the same asymptotics as Theorem~\ref{thm:multilayer joint MI}, it follows that
    \begin{align}
        \lim_{n\rightarrow\infty}\frac{1}{n}I(\bX^{(l)};\mathbf{G})&=\frac{\lambdal}{4}+\sup_{\mathbf{q}_{-l}\ge 0}\left[\bar{\mathcal{F}}^{(l)}_{\ML,\rho}(\boldsymbol{\lambda}^{(-l)},\mathbf{q}_{-l})-\sum_{l_1\neq l}\frac{\lambda^{(l_1)}(q_{l_1}^2+2q_{l_1})}{4}\right]\notag\\
        &\qquad -\sup_{\mathbf{q}\ge 0}\left[\mathcal{F}_{\ML,\rho}(\boldsymbol{\lambda},\mathbf{q})-\sum_{l=1}^{L}\frac{\lambda^{(l)}(q_l^2+2q_l)}{4}\right].\label{eq:local MI multilayer}
    \end{align}
\end{Corollary}

\paragraph{AMP Denoiser and State Evolution.}
After deducing the asymptotics of mutual information, we turn to an investigation of the algorithmic questions related to community detection. We describe the specializations of  Algorithm~\ref{alg:coupled AMP} to this multilayer model. In the \textit{Bayes optimal} setting, given the state evolution iterates $\{\mathbf{q}^t\}$, the non-linear denoisers $\mathcal{E}_t:\R^{L}\rightarrow\R^L$ are chosen as 
\begin{align*}
    &\mathcal{E}^{(l)}_t(m_{1:L})=\mathbb{E}\left[X^{(l)}\bigg| \sqrt{\lambda^{(l_1)} q_{l_1}^t}A^{(l_1)\prime}=m_l,\forall l_1\in[L]\right]
\end{align*}
for each $l\in[L]$ in the spirit of optimizing the reconstruction accuracy.  The following lemma provides an explicit expression for this denoiser.  
\begin{Lemma}\label{lemma:ML denoiser}
    For this multilayer model $p_{\ML}$, denoisers admit such a closed-form expression
    \begin{equation}\label{eq:ML denoiser}
        \mathcal{E}_{\ML,\rho}^{(l)}(m_{1:L})=\frac{\tanh\left(m_l+\bar{\rho}\right)\prod_{l^\prime}\cosh\left(m_{l^\prime}+\bar{\rho}\right)+\tanh\left(m_l-\bar{\rho}\right)\prod_{l^\prime}\cosh\left(m_{l^\prime}-\bar{\rho}\right)}{\prod_{l^\prime}\cosh\left(m_{l^\prime}+\bar{\rho}\right)+\prod_{l^\prime}\cosh\left(m_{l^\prime}-\bar{\rho}\right)},
    \end{equation}
    with $\bar{\rho}=\frac{1}{2}\log\left(\frac{1-\rho}{\rho}\right)$. Consequently, Bayes optimal AMP does not have to make denoisers vary over iterations.
\end{Lemma}

Armed with this denoiser, we obtain an explicit representation SE parameters. 
Specifically, $\{\mathbf{q}^t\in\R^L:t\in\mathbb{N}\}$ is recursively updated as 
\begin{equation*}
    q_l^{t+1}=\mathbb{E}\left[X^{(l)}\mathcal{E}_{\ML,\rho}^{(l)}(\boldsymbol{\lambda}\odot\mathbf{q}^t\odot X+\sqrt{\boldsymbol{\lambda}\odot\mathbf{q}^t}\odot W^\prime)\right]=T_l^{\ML}(\boldsymbol{\lambda}\odot\mathbf{q}^t),
\end{equation*}
where $T^\ML$ is given in \eqref{eq:general se mapping} by specifying the multilayer prior $p_\ML$. As a result, by recursively updating $\mathbf{q}^{t+1}\leftarrow T^{\ML}(\boldsymbol{\lambda}\odot\mathbf{q}^t)$, state evolution iterates $\mathbf{q}^t$ would converge to a fixed point of mapping $\mathbf{q}\mapsto T^{\ML}(\boldsymbol{\lambda}\odot\mathbf{q})$. If in addition this fixed point is also the global unique global maximizer of the free energy functional, our coupled AMP algorithm is indeed Bayes optimal.

\paragraph{Weak Recovery Thresholds.} Denote $D$ as the set of $\boldsymbol{\lambda}$ such that \eqref{eq:opt task} does not have a unique maximizer; Theorem~\ref{thm:asymptotic MI} establishes that $D$ has measure zero. For any $\boldsymbol{\lambda}\notin D$, write $\mathbf{q}^\ast(\boldsymbol{\lambda})$ as the unique global maximizer. As a result, understanding whether $q_l^\ast(\boldsymbol{\lambda})=0$ or not determines whether or not weak recovery is possible for the $l$-th layer for $\boldsymbol{\lambda}$. 

Proposition~\ref{prop:all-in-all-out} identifies an all-or-nothing phenomenon across all individualized memberships $\bX^{(l)}$ and the global membership $\bY$. Thus we can focus on ``the" threshold for weak recovery in this example. 
We proceed to derive explicit conditions for weak recovery by checking whether or not $\mathbf{0}$ is a \textit{global} maximizer of \eqref{eq:multilayer joint MI}. Since $\mathbf{0}$ is always a saddle point, we have to examine the second-order properties of the functional at this point. The following theorem establishes a sufficient condition for weak recovery by identifying the condition under which $\mathbf{0}$ is not even a \textit{local} maximizer.

\begin{Proposition}\label{prop:multilayer-threshold}
    Suppose all $\lambdal\le 1$ so that weak recovery is impossible based on any individual layer. With $0\le \rho<1/2$ and $\boldsymbol{\lambda}\notin D$, the origin $\mathbf{0}$ is not a local maximizer of \eqref{eq:multilayer joint MI} if and only if
    \begin{equation}\label{eq:threshold-inhomo-SBM}
        \sum_{l = 1}^L \frac{(1-2\rho)^4\lambda^{(l)}}{1-(1-(1-2\rho)^4)\lambda^{(l)}}>1.
    \end{equation}
    Therefore, when \eqref{eq:threshold-inhomo-SBM} holds, weak recovery is (information-theoretically) possible for any individualized layer $\bX^{(l)}$ and global membership $\bY$.
\end{Proposition}

\begin{figure}
    \begin{subfigure}{0.58\textwidth}
        \centering
        \includegraphics[width=\linewidth]{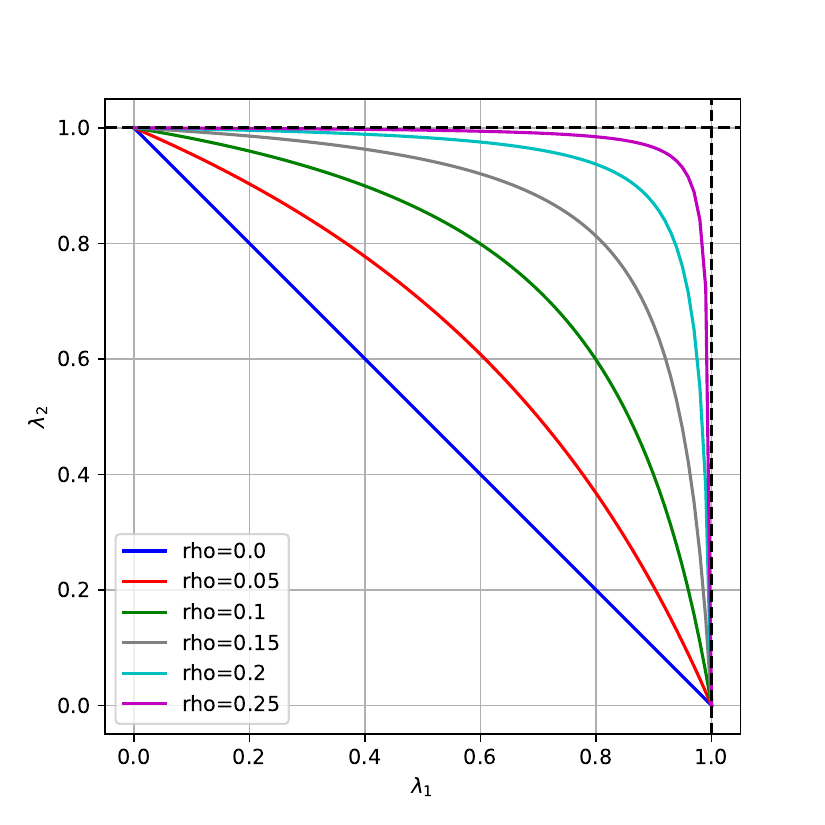}
        \caption{Feasible region of SNRs in 2-Layer SBM}
        \label{fig:feasible_region_2layer}
    \end{subfigure}
    \begin{subfigure}{0.42\textwidth}
        \centering
        \includegraphics[width=\linewidth]{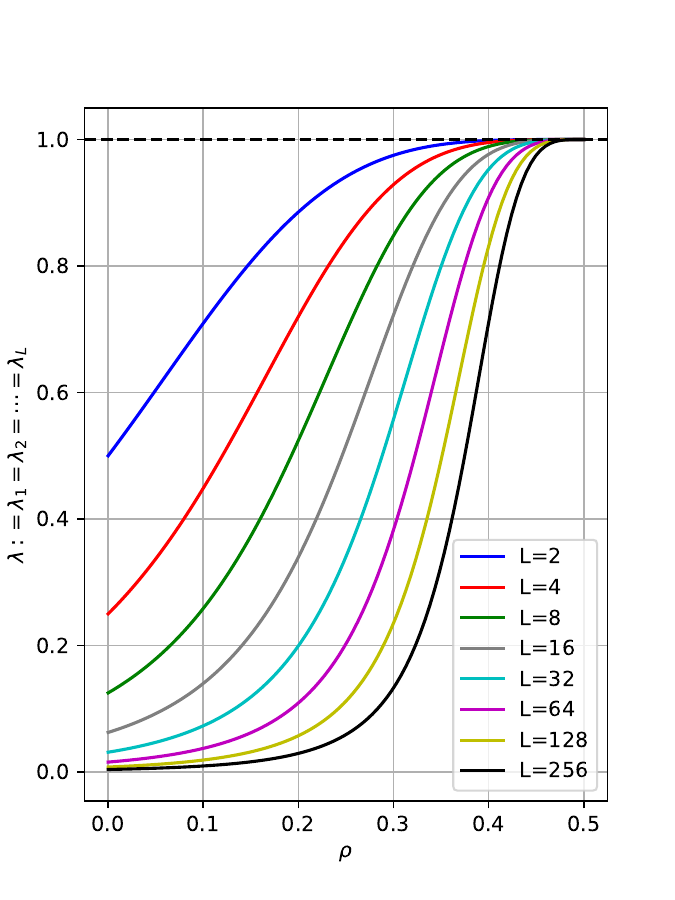}
        \caption{Weak recovery thresholds}
        \label{fig:lam_curve_multilayer}
    \end{subfigure}
    \caption{Illustration of Weak Recovery Thresholds \eqref{eq:threshold-inhomo-SBM} for Multilayer Inhomogeneous SBM}
\end{figure}

Subsequently, one may ask if \eqref{eq:threshold-inhomo-SBM} is a necessary condition for weak recovery. In terms of the free energy functional, we wonder whether or not $\mathbf{0}$ is the \textit{global} maximizer if it is already a \textit{local} one. After extensive numerical experiments, we make the following conjecture on the structure of the map $T^{\ML}$. 
\begin{Conjecture}\label{conjecture:ML}
    For the balanced binary multilayer prior $p_{\ML}$ with $0\le \rho\le 1/2$, every coordinate $\boldsymbol{\gamma}\mapsto T_l^\ML(\boldsymbol{\gamma})$ of its state evolution mapping is strictly concave on every straight line passing through $\mathbf{0}$, i.e. $t\in[0,+\infty)\mapsto T_l^\ML(t\boldsymbol{\gamma})$ is strictly concave for every $l\in[L]$ and non-zero $\gamma\in[0,\infty)^L$.
\end{Conjecture}

\begin{figure}
    \centering
    \makebox[\linewidth][c]{
        \includegraphics[width=1.1\linewidth]{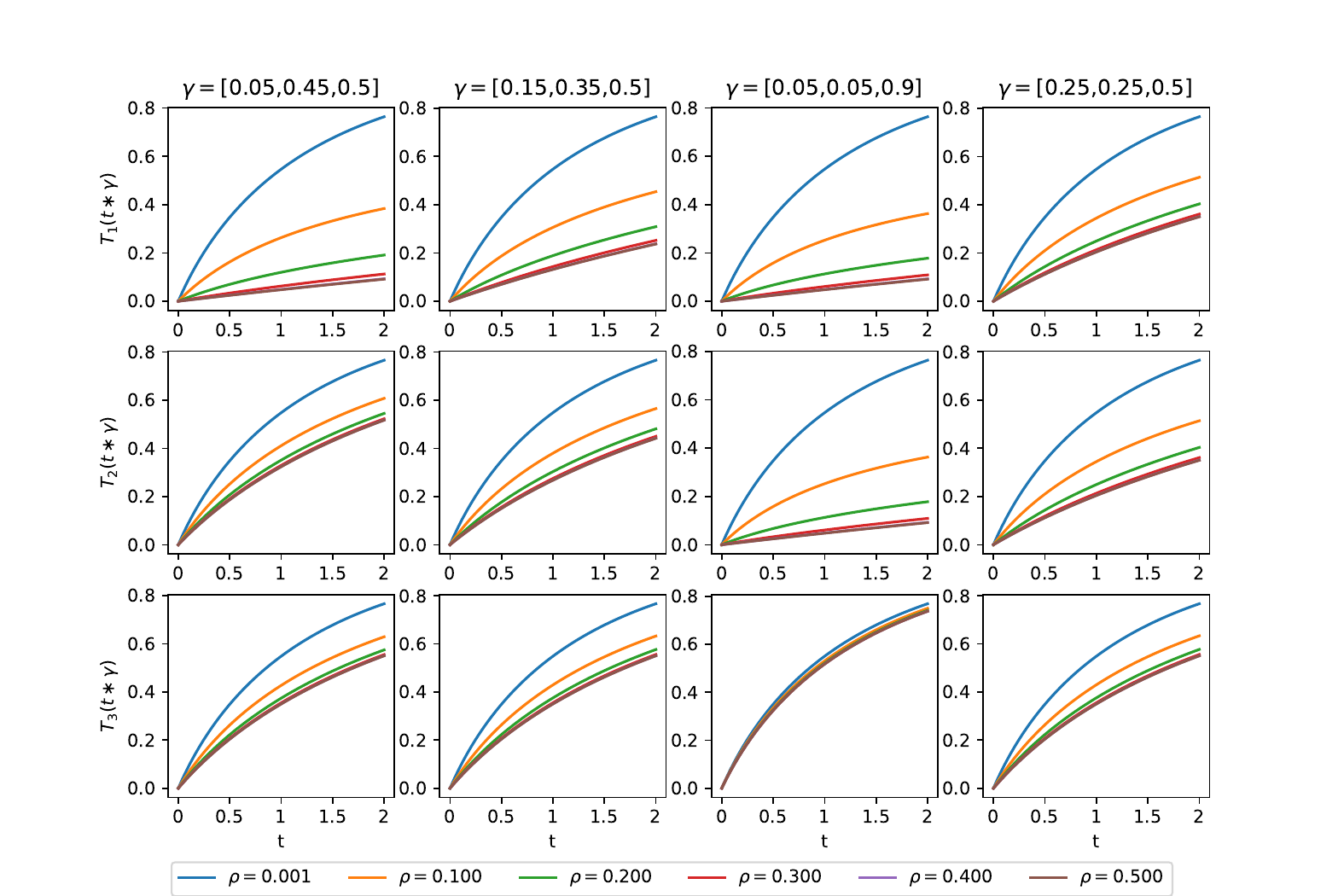}
    }
    \caption{This figure plots the mapping $t\in[0,2]\mapsto T_l^{\ML,\rho}(t\boldsymbol{\gamma})$ for four different choices of $\gamma\in\mathbb{R}_+^3$ and six different choices of $\rho\in[0,1/2]$, in the context of a $3$-layer inhomogeneous SBM as Example~\ref{eg:multilayer SBM}. Conjecture~\ref{conjecture:ML} claims these mappings to be concave all along $t\in[0,\infty)$, and is a key condition to show the optimality of Algorithm~\ref{alg:coupled AMP} as in Proposition~\ref{prop:ML conjectured}.}
    \label{fig:conjecture ML}
\end{figure}
See Figure~\ref{fig:conjecture ML} for a numerical evaluation of the state evolution mappings we conjecture about. In particular for the extremal cases $\rho=0$ or $\rho=1/2$, this multilayer model can both be reduced to a single layer setting, thus Conjecture~\ref{conjecture:ML} is strictly proved in \cite{deshpande2017asymptotic}. As shown by the numerical evaluations, we believe this phenomenon to be ubiquitous for every $\rho\in[0,1/2]$.

Our last proposition for this multilayer model is to prove necessity of our weak recovery threshold and Bayes optimality of the coupled AMP algorithm.
\begin{Proposition}\label{prop:ML conjectured}
    Suppose Conjecture~\ref{conjecture:ML} is true. Then the following holds.
    \begin{itemize}
        \item[(i)] For any\footnote{Compared to Theorem~\ref{thm:asymptotic MI}, this Conjecture~\ref{conjecture:ML} helps us rule out the bad set $D$ of measure zero.} $\boldsymbol{\lambda}\in(0,\infty)^L$, there exists a unique $\mathbf{q}^\ast$ to the optimization problem \eqref{eq:multilayer joint MI} and \eqref{eq:limiting mmse} holds.
        \item[(ii)] If condition \eqref{eq:threshold-inhomo-SBM} fails, $\mathbf{0}$ would be a \textit{global} maximizer to \eqref{eq:multilayer joint MI}. Thus, weak recovery is impossible.
        \item[(iii)] If condition \eqref{eq:threshold-inhomo-SBM} holds, when running Algorithm~\ref{alg:coupled AMP}, whenever initialization $\mathbf{m}^0$ is better than random guessing, i.e.
        \begin{equation*}
            \lim_{n\rightarrow\infty}\frac{1}{n}\mathbb{E}\left[\bX^{(l)\top}\mathbf{m}^0_{\cdot,l}\right]>0
        \end{equation*}
        for some $l\in[L]$, state evolution iterates $\mathbf{q}^t$ converges to $\mathbf{q}^\ast$ exponentially. So the algorithm is Bayes optimal and runs in polynomial time.
    \end{itemize}
\end{Proposition}

\noindent Figure~\ref{fig:feasible_region_2layer} depicts the feasible region of $(\lambda^{(1)},\lambda^{(2)})$ in the context of a $2$-layer model with varying $\rho$. For each $\rho$, weak recovery is possible in the region above the curve. Figure~\ref{fig:lam_curve_multilayer} then displays the threshold when $\lambda^{(1)}=\lambda^{(2)}=\cdots=\lambda^{(L)}$. Given $L$ of layers, weak recovery is (information-theoretically) possible only when $\lambda$ is above the curve. 
Before moving further, we collect some other observations by considering special cases. 
\begin{itemize}
    \item[1.] If $\rho=0$, the model degenerates to a homogeneous multilayer setting studied in \cite{ma2023community}. Our threshold reduces to $\sum_{l}\lambdal>1$, recovering the threshold derived in their work. 
    \item[2.] In \eqref{eq:threshold-inhomo-SBM}, RHS decreases with respect to increasing $\rho$, and approaches $0$ as $\rho\rightarrow1/2$. In this limit, different layers are uncorrelated, and no recovery should be possible since we already set all $\lambdal\le 1$.
    \item[3.] If we restrict to the case all $\lambdal$ equal $\lambda$, then the threshold translates to $\lambda>\frac{1}{1+(L-1)(1-2\rho)^4}$. We  interpret different layers as i.i.d. noisy observations on the global membership $\bY$. No matter how small effective SNR $\lambda>0$, conducting enough experiments (i.e. letting $L\rightarrow\infty$) eventually leads to weak recovery. 
\end{itemize}

\begin{figure}
    \centering
    \makebox[\linewidth][c]{
        \includegraphics[width=1.2\linewidth]{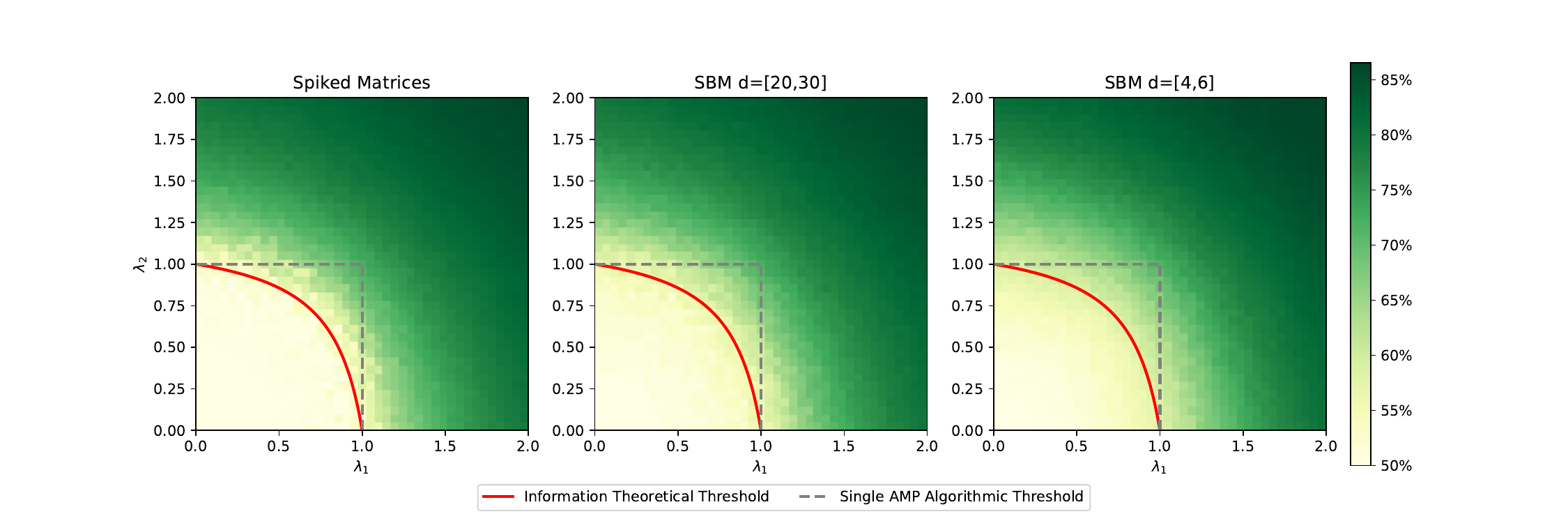}
    }
    \caption{Simulation results on an inhomogenous $2$-layer model, formally defined in Example~\ref{eg:multilayer SBM}. Algorithm~\ref{alg:coupled AMP} is tested on spiked matrices, sparse graphs with $d=(20,30)$ or $d=(4,6)$ respectively. Color brightness reflect global membership recovery accuracy averaged over $10$ repetitions. The empirical feasibility boundary matches well with the theoretically predicted red curve. Although our theory applies only to denser graphs with $d\rightarrow\infty$ as $n$ increases, the simulation suggests the boundary to be critical even for extremely sparse graphs: $d$ is set to $\sim5$ compared to $n=10000$.}
    \label{fig:simu_ML}
\end{figure}

\paragraph{Simulation Studies.} We close our discussion on this example with a simulation study on an inhomogenous $2$-layer model with fixed $\rho=0.1$, as shown in Figure~\ref{fig:simu_ML}. Tested models include the spiked matrix counterpart~\eqref{eq:spiked matrix each layer}, a sparse graph setting with average degrees $d=(20,30)$ and an extremely sprse one with $d=(4,6)$.
We run Algorithm~\ref{alg:coupled AMP} with Bayes optimal denoisers specified in Lemma~\ref{lemma:ML denoiser} for $1600$ different choices of $(\lambda^{(1)},\lambda^{(2)})\in[0,2]^2$ and $n=10000$ nodes. For each configuration, the algorithm starts with a warm initialization\footnote{About $10\%$ of $\{m_{i,l}^0:i\in[n],l\in[L]\}$ equal the underlying true community membership, with rest set to $0$.} $\mathbf{m}^0$ and proceeds by $100$ iterations to output an estimate $\hat{\bx}\in\R^{n\times L}$ for individualized memberships of all layers. 
For each node, we then pass $\hat{x}_{i}\in\R^L$ to another global denoiser to derive an estimate $\hat{y}_i\in\{\pm 1\}$ on the global community membership $Y_i$. Lastly the performance is measured by the proportion $\frac{1}{n}|\{i\in[n]:\hat{y}_i=Y_i\}|$ of accurately recovered nodes. A random guess should yield $50\%$ accuracy, while our warm initialization yields $10\%/2+50\%=55\%$ accuracy in signal recovery. Below the detection threshold, the algorithm forgets leaked initialization information with increasing iterations. Repeating the simulation process $10$ times for each configuration of $(\lambda^{(1)},\lambda^{(2)})$, the averaged recovery accuracy is reflected by color brightness in Figure~\ref{fig:simu_ML}.

\subsection{Applications to Dynamic SBM}
We now proceed to apply our results to dynamical stochastic block models, introduced in Example~\ref{eg:dynamic SBM}.
Under this input prior $p_\Dyn$, the effective scalar channel \eqref{eq:scalar channel} turns into a hidden Markov model indexed by $l\in[L]$ as demonstrated in Figure~\ref{fig:Dyn illustration}, with binary discrete states $X^{(1:L)}$ and Gaussian observations $A^{(1:L)\prime}$: 
\begin{align}
    X^{(1)}\sim \text{Unif}\left\{\pm 1\right\},\quad &A^{(1)\prime}\sim\sqrt{\lambda^{(1)} q_1}X^{(1)}+\mathcal{N}\left(0,1\right),\notag\\
    X^{(l)}\sim \mathsf{BSC}_{\rho}\left(X^{(l-1)}\right), \quad &A^{(l)\prime}\sim\sqrt{\lambdal q_l}X^{(l)}+\mathcal{N}\left(0,1\right),\quad 2\le l\le L.\label{eq:scalar channel dyn}
\end{align}
As defined in \eqref{eq:scalar channel free energy}, the free energy functional of this scalar channel is given by
\begin{align*}
    \mathcal{F}_{\Dyn,\rho}(\boldsymbol{\lambda},\mathbf{q})&=\mathbb{E}\log\left\{\sum_{x \in \{ \pm 1\}^L}\frac{1}{2} \exp\left(\lambda^{(1)} q_1X^{(1)}x^{(1)}+\sqrt{\lambda^{(1)} q_1}W^{(1)\prime}x^{(1)}\right)\right.\\
    &\qquad\qquad\qquad \prod_{l=2}^L\left. p\left(x^{(l)}|x^{(l-1)}\right)\exp\left(\lambdal q_lX^{(l)}x^{(l)}+\sqrt{\lambdal q_l}W^{(l)\prime}x^{(l)}\right)\right\}.
\end{align*}
This chain structure prevents us from deriving more explicit expression for the free energy functional. We emphasize that this recursive summation can be evaluated in $O(L)$ time, instead of the worst case $O(2^L)$ time. Consequently, Theorem~\ref{thm:asymptotic MI} and Proposition~\ref{prop:Universality MI qualitative} immediately identify the asymptotic normalized mutual information between observed graphs $\mathbf{G}$ and underlying latent variables $\mathbf{X}$.
\begin{Theorem}\label{thm:dyn joint MI}
    Under the asymptotics that (i) $\blam_n\rightarrow\blam$ and (ii) $\min_l d^{(l)}(1-d^{(l)}/n)\rightarrow\infty$ as $n\rightarrow\infty$, we have,
    \begin{align}
        &\quad\lim_{n\rightarrow\infty}\frac{1}{n}I(\bX;\mathbf{G})=\lim_{n\rightarrow\infty}\frac{1}{n}I(\bX;\bA)\notag\\
        &=\frac{\sum_l\lambdal}{4}-\sup_{\mathbf{q}\ge 0}\left[\mathcal{F}_{\Dyn,\rho}(\boldsymbol{\lambda},\mathbf{q})-\sum_{l=1}^{L}\frac{\lambda^{(l)}(q_l^2+2q_l)}{4}\right],\label{eq:dyn joint MI}
    \end{align}
    where we recall $\mathbf{A}$ to be the approximating spiked matrices in \eqref{eq:spiked matrix each layer}. 
\end{Theorem}

\noindent One might also be interested in inferring community labels $\bX^{(l)}$ at a certain moment $l\in[L]$ only. A direct corollary provides an exact formula for the limiting asymptotic mutual information between $\bX^{(l)}$ and $\bG$.

\begin{Corollary}\label{cor:local dyn}
Due to the sequential structure depicted in Figure~\ref{fig:Dyn illustration}, we denote $\bX^{(>l)}=\left(\bX^{(l+1)},\ldots,\bX^{(L)}\right)$ and $\bX^{(<l)}=\left(\bX^{(1)},\ldots,\bX^{(l-1)}\right)$ respectively. Then $\bG^{(>l)}$ and $\bG^{(<l)}$ follow naturally. Define free energy functionals as
\begin{align*}
    \mathcal{F}_{\Dyn,\rho}(\boldsymbol{\lambda}^{(>l)},\mathbf{q}_{>l})&=\mathbb{E}\log\Bigg\{\sum_{x^{(>l)} \in \{ \pm 1\}^{L-l}} p\left(x^{(l+1)}|X^{(l)}\right)\prod_{l^\prime=l+2}^L p\left(x^{(l^\prime)}|x^{(l^\prime-1)}\right)\\
    &\qquad\qquad\qquad \prod_{l^\prime=l+1}^L\exp\left(\lambda^{(l^\prime)} q_{l^\prime}X^{(l^\prime)}x^{(l^\prime)}+\sqrt{\lambda^{l^\prime} q_{l^\prime}}W^{(l^\prime)\prime}x^{(l^\prime)}\right)\Bigg\}.
\end{align*}
Moreover, $\mathcal{F}_{\Dyn,\rho}(\boldsymbol{\lambda}^{(<l)},\mathbf{q}_{<l})$ is defined analogously. Then under the same asymptotics as Theorem~\ref{thm:dyn joint MI}, we have,
\begin{align}
    \lim_{n\rightarrow\infty}\frac{1}{n}I(\bX^{(l)};&\bG)=\frac{\lambdal}{4}-\sup_{\mathbf{q}\ge 0}\left[\mathcal{F}_{\Dyn,\rho}(\boldsymbol{\lambda},\mathbf{q})-\sum_{l=1}^{L}\frac{\lambda^{(l)}(q_l^2+2q_l)}{4}\right]\notag\\
        &+\sup_{\mathbf{q}_{>l}\ge 0}\left[\mathcal{F}_{\Dyn,\rho}(\boldsymbol{\lambda}^{(>l)},\mathbf{q}_{>l})-\sum_{l^\prime>l}\frac{\lambda^{(l^\prime)}(q_{l^\prime}^2+2q_{l^\prime})}{4}\right]\notag\\
        &+\sup_{\mathbf{q}_{<l}\ge 0}\left[\mathcal{F}_{\Dyn,\rho}(\boldsymbol{\lambda}^{(<l)},\mathbf{q}_{<l})-\sum_{l^\prime<l}\frac{\lambda^{(l^\prime)}(q_{l^\prime}^2+2q_{l^\prime})}{4}\right].\label{eq:local dyn}
\end{align}
\end{Corollary}

\paragraph{AMP Denoiser and State Evolution.}
To implement our Algorithm~\ref{alg:coupled AMP}, we have to specify the denoisers. 
Similar to the previous example, given SE parameters $\mathbf{q}^t$, the
Bayes optimal denoisers should be chosen as
\begin{align*}
    &\mathcal{E}^{(l)}_{\Dyn,\rho}(m_{1:L})=\mathbb{E}\left[X^{(l)}\bigg| \sqrt{\lambdal q_l^t}A^{(l)\prime}=m_l,\forall l\in[L]\right]
\end{align*}
for each $l\in[L]$.
Under the chain structure of prior $p_\Dyn$ (Figure~\ref{fig:Dyn illustration}), an explicit formula for denoisers $\mathcal{E}_\Dyn$ seems out of reach. However, the following lemma provides a  dynamical programming approach to efficiently evaluate the denoiser $O(L)$ time. Our approach borrows the idea of a  Kalman filter \cite{welch1995introduction}, developed originally for a continuous state space model with Gaussian transition and observation.
\begin{Lemma}\label{lemma:dyn denoiser}
    For dynamical model $p_\Dyn$, denoisers $\mathcal{E}_{\Dyn,\rho}$ can be explicitly evaluated by Algorithm~\ref{alg:denoiser dyn}. 
\end{Lemma}
Upon specifying the appropriate class of denoisers for the problem, we obtain an explicit update scheme for the SE parameters. 
In particular, $\{\mathbf{q}^t\in\R^L:t\in\mathbb{N}\}$ is recursively given by
\begin{equation*}
    q_l^{t+1}=\mathbb{E}\left[X^{(l)}\mathcal{E}_{\Dyn,\rho}^{(l)}(\boldsymbol{\lambda}\odot\mathbf{q}^t\odot X+\sqrt{\boldsymbol{\lambda}\odot\mathbf{q}^t}\odot W^\prime)\right]=T_l^{\Dyn}(\boldsymbol{\lambda}\odot\mathbf{q}^t),
\end{equation*}
where $T^\Dyn$ is given in \eqref{eq:general se mapping} by specifying dynamical prior $p_\Dyn$. As a result, in an idealized setting, by recursively updating $\mathbf{q}^{t+1}\leftarrow T^{\Dyn}(\boldsymbol{\lambda}\odot\mathbf{q}^t)$, state evolution iterates $\mathbf{q}^t$ would converge to a fixed point of mapping $\mathbf{q}\mapsto T^{\Dyn}(\boldsymbol{\lambda}\odot\mathbf{q})$. If in addition this fixed point is also the global unique global maximizer of the free energy functional, our coupled AMP algorithm is indeed Bayes optimal.

\begin{algorithm}[tb]
\caption{Node-wise Denoiser for Dynamical SBM}\label{alg:denoiser dyn}
\begin{algorithmic}
\Require Sequence of estimates $\{m_{1:L}\}$ to be denoised, flipping rate $\rho\in[0,1/2]$.
\State Initialize table $g_{1:L}(\pm 1)$ by $g_1(1)=\frac{1}{2}\exp(m_1) $ and $g_1(-1)=\frac{1}{2}\exp(-m_1)$.
\For{$l=2,\ldots,L$}
    \State Recursively update table $g$ by
    \begin{align*}
        g_{l}(1)&=(1-\rho)\exp(m_l)g_{l-1}(1)+\rho\exp(m_l)g_{l-1}(-1),\\
        g_{l}(-1)&=\rho\exp(-m_l)g_{l-1}(1)+(1-\rho)\exp(-m_l)g_{l-1}(-1).
    \end{align*}    
\EndFor
\State Start constructing $\mathcal{E}\in\R^L$ by firstly setting $\mathcal{E}^{(L)}=\frac{g_L(1)-g_L(-1)}{g_L(1)+g_L(-1)}$.
\For{$l=L-1,\ldots,1$}
    \State Recursively update $\mathcal{E}$ by
    \begin{align*}
        \mathcal{E}^{(l)}&=\frac{1+\mathcal{E}^{(l+1)}}{2}\cdot\frac{(1-\rho)\exp(m_{l+1})g_{l}(1)-\rho\exp(m_{l+1})g_{l}(-1)}{g_{l+1}(1)}\\
        &\quad+\frac{1-\mathcal{E}^{(l+1)}}{2}\cdot\frac{\rho\exp(-m_{l+1})g_{l}(1)-(1-\rho)\exp(-m_{l+1})g_{l}(-1)}{g_{l+1}(-1)}.
    \end{align*}
\EndFor
\State Return $\mathcal{E}_{\Dyn,\rho}(m)=\left(\mathcal{E}^{(1)},\ldots,\mathcal{E}^{(L)}\right)$.
\end{algorithmic}
\end{algorithm}

\paragraph{Weak Recovery Thresholds.}
We turn to a study of weak recovery thresholds in this example. As before, we use $D$ to denote the set of $\boldsymbol{\lambda}$ such that \eqref{eq:opt task} does not have a unique maximizer. Theorem~\ref{thm:asymptotic MI} ensures that $D$ has measure zero. For any $\boldsymbol{\lambda}\notin D$, define $\mathbf{q}^\ast(\boldsymbol{\lambda})$ to be  the unique global maximizer. As a result, weak recovery thresholds reduce to determining whether $q_l^\ast(\boldsymbol{\lambda})=0$ or not. 

Recall that an all-or-nothing phenomenon is also identified for the dynamical model when $\rho<1/2$ in Proposition~\ref{prop:all-in-all-out}. We can thus focus on determing ``the" weak recovery threshold. To this end, we proceed to derive explicit conditions for weak recovery by checking the global optimality of $\mathbf{0}$ in \eqref{eq:dyn joint MI}. Since $\mathbf{0}$ is always a saddle point, we have to examine the second-order properties of the functional at this point. 
By Proposition~\ref{prop:opt derivative}(i), the optimization objective in \eqref{eq:dyn joint MI} has an explicit form
\begin{equation*}
    \nabla^2 G(\mathbf{0})=\frac{1}{2}\begin{pmatrix}
        \lambda^{(1)2}-\lambda^{(1)} & \lambda^{(1)}\lambda^{(2)}(1-2\rho)^2 & & \lambda^{(1)}\lambda^{(L)}(1-2\rho)^{2(L-1)}\\
        \lambda^{(1)}\lambda^{(2)}(1-2\rho)^2 & \lambda^{(2)2}-\lambda^{(2)} & & \lambda^{(2)}\lambda^{(L)}(1-2\rho)^{2(L-2)}\\
        & & \ddots & \\
        \lambda^{(1)}\lambda^{(L)}(1-2\rho)^{2(L-1)} & & & \lambda^{(L)2}-\lambda^{(L)}
    \end{pmatrix}.
\end{equation*}
for its Hessian matrix at the origin $\mathbf{0}$. 
The following remark collects a condition under which $\mathbf{0}$ is not even a \textit{local} maximizer, as a sufficient condition for weak recovery.
\begin{Remark}\label{remark:threshold-dynamic-SBM}
Weak recovery for dynamic SBM is feasible if $\boldsymbol{\lambda}$ admits a vector $\mathbf{v}\in(0,\infty)^L$ such that $ \mathbf{v}^\top \nabla^2 G(\mathbf{0})\mathbf{v}>0$.
\end{Remark}

For general $\boldsymbol{\lambda}$, this condition can be efficiently checked by solving the above quadratic programming problem using existing numerical packages. In the special case that all $\lambda^{(l)}=\lambda$, the Hessian matrix $\nabla^2 G(\mathbf{0})$ reduces to a Toeplitz matrix of a special form, referred to as the Kac-Murdock-Szego matrices in the literature of quadratic forms \cite{trench2010spectral,grenander1958toeplitz}. In this case, one can derive an explicit formula for $\mathbf{0}$ to be a local maximizer.
\begin{Proposition}\label{prop:dynamic-threshold}
    Assume $d = \min_{1\leq l \leq L} \frac{a^{(l)}+ b^{(l)}}{2} \to \infty$. If $\lambda^{(1)} = \lambda^{(2)} = \cdots = \lambda^{(L)} := \lambda$, then for $0<\rho<1/2$, weak recovery is possible for the dynamic SBM if
    \begin{equation}
        \lambda > \frac{1-2(1-2\rho)^2\cos\theta_\ast+(1-2\rho)^4}{1-(1-2\rho)^4}, \label{eq:threshold_dynamic}
    \end{equation}
    where $\theta_\ast\in(0, \pi)$ is the minimum solution of equation
    \begin{equation*}
        0 = \sin[(L+1)\theta_\ast]-2(1-2\rho)^2\sin[ L\theta_\ast]+(1-2\rho)^4\sin[(L-1)\theta_\ast].
    \end{equation*}
\end{Proposition}


\begin{Corollary}
Assume $\lambda^{(1)} = \lambda^{(2)} = \cdots = \lambda^{(L)} := \lambda$. Then for $0<\rho<1/2$, as $L\rightarrow\infty$, the threshold converges to 
\begin{equation*}
    \lambda>\frac{1-(1-2\rho)^2}{1+(1-2\rho)^2}.
\end{equation*}   
\end{Corollary}
\begin{proof}
It is proved in \cite{trench2010spectral,grenander1958toeplitz} that $0<\theta_\ast<\pi/(L+1)$ for Kac-Murdock-Szego matrices. Therefore, as $L\rightarrow\infty$, the minimal solution $\theta_\ast\rightarrow0$. 
In turn, the threshold reduces to $\frac{1-(1-2\rho)^2}{1+(1-2\rho)^2}$.
\end{proof}

\begin{figure}
    \centering
    \includegraphics[width=0.4\linewidth]{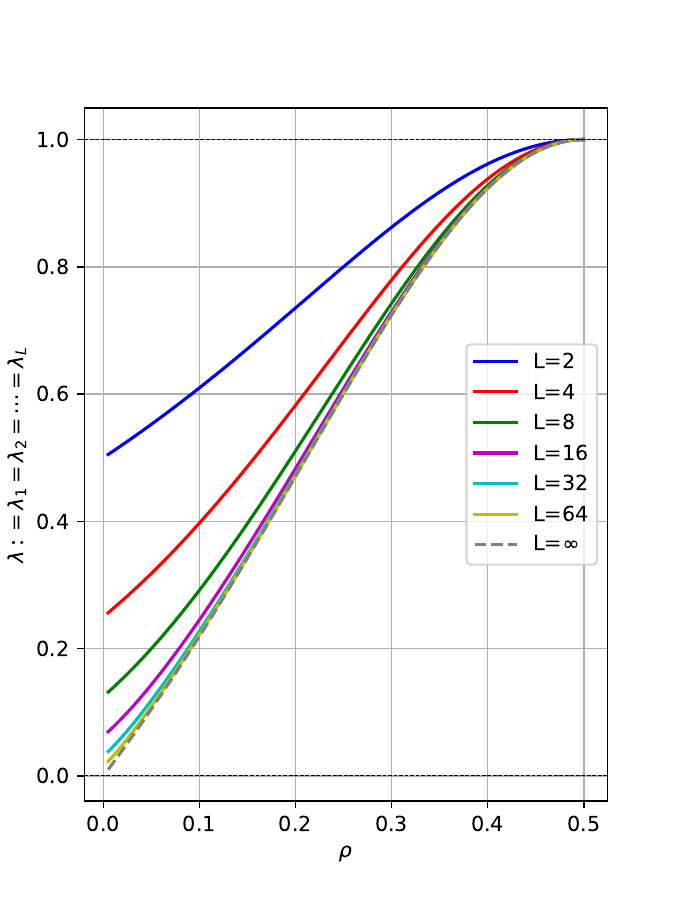}
    \caption{Weak recovery thresholds of dynamic SBM when $\lambda^{(1)}=\cdots=\lambda^{(L)}$.}
    \label{fig:lam_curve_dyn}
\end{figure}

In sharp contrast to inhomogeneous multilayer SBM, the requirement on $\lambda$ does not vanish as $L\rightarrow\infty$. Even if we get to observe infinite many networks at every moment $l\in\mathbb{N}$, we still need the signal strength to be larger than a threshold in order to weakly recover the community structure, as shown in Figure~\ref{fig:lam_curve_dyn}.

Subsequently, one may wonder if Remark~\ref{remark:threshold-dynamic-SBM} provides a necessary condition for weak recovery. In terms of free energy functional, this corresponds to the global optimality of $\mathbf{0}$. Based on  extensive numerical experiments, we come up with a relevant conjecture.
\begin{Conjecture}\label{conjecture:dyn}
    For this balanced binary multilayer prior $p_{\Dyn}$ with $0\le \rho\le 1/2$, every coordinate $\boldsymbol{\gamma}\mapsto T_l^\Dyn(\boldsymbol{\gamma})$ of its state evolution mapping is strictly concave in every straight line passing $\mathbf{0}$, i.e. $t\in[0,+\infty)\mapsto T_l^\Dyn(t\boldsymbol{\gamma})$ is strictly concave for every $l\in[L]$ and non-zero $\gamma\in[0,\infty)^L$.
\end{Conjecture}

\begin{figure}
    \centering
    \makebox[\linewidth][c]{
        \includegraphics[width=1.1\linewidth]{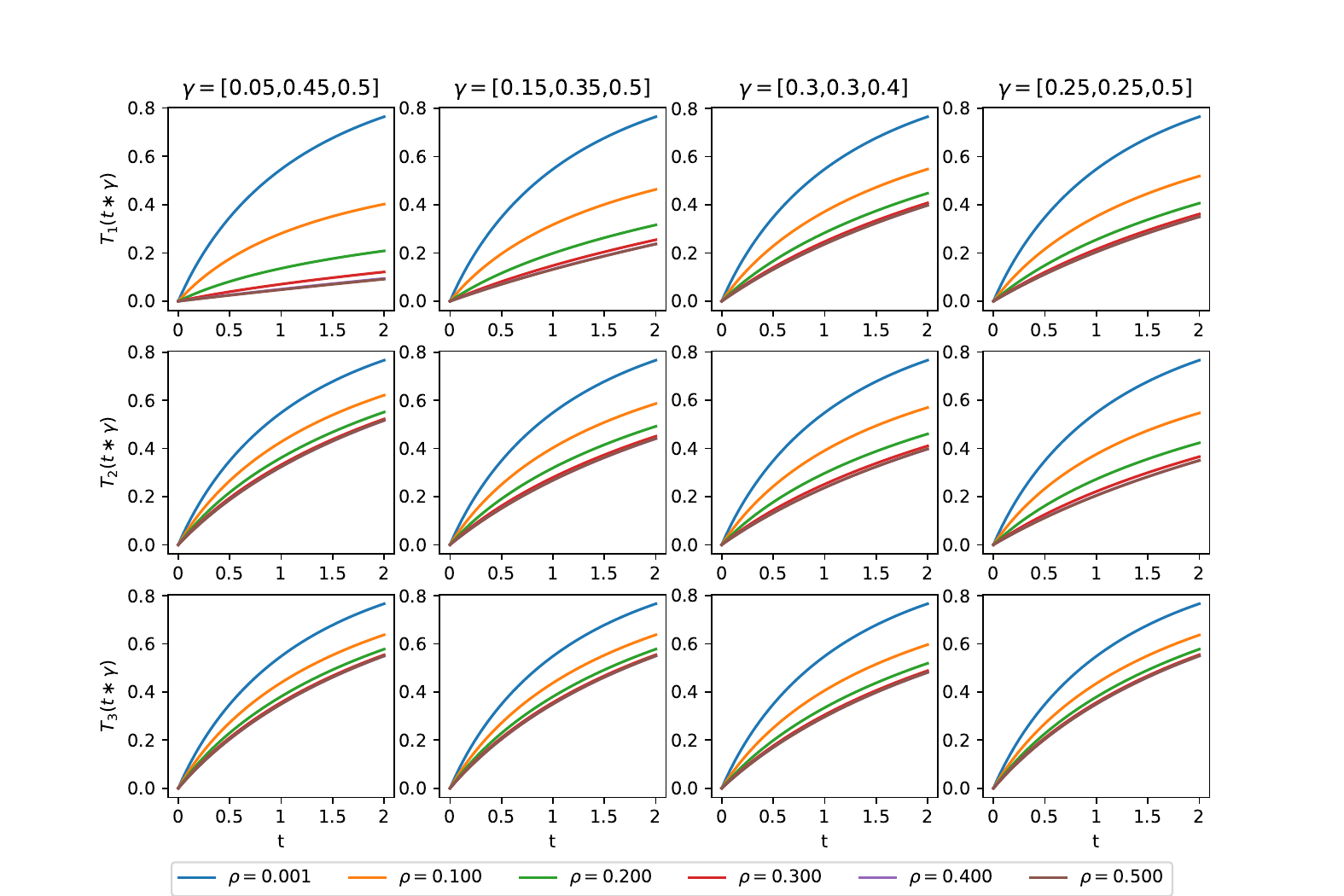}
    }
    \caption{This figure plots the mapping $t\in[0,2]\mapsto T_l^{\Dyn,\rho}(t\boldsymbol{\gamma})$ for four different choices of $\gamma\in\mathbb{R}_+^3$ and six different choices of $\rho\in[0,1/2]$, in the context of a $3$-epoch dynamical SBM (Example~\ref{eg:dynamic SBM}). Conjecture~\ref{conjecture:dyn} claims these mappings to be concave all along $t\in[0,\infty)$, and is a key condition to show the optimality of Algorithm~\ref{alg:coupled AMP} in Proposition~\ref{prop:dyn conjectured}.}
    \label{fig:conjecture Dyn}
\end{figure}
Our last proposition for this multilayer model is to prove necessity of our weak recovery threshold and Bayes optimality of the coupled AMP algorithm.
\begin{Proposition}\label{prop:dyn conjectured}
    Suppose Conjecture~\ref{conjecture:dyn} is true. Then we have,
    \begin{itemize}
        \item[(i)] For any\footnote{Compared to Theorem~\ref{thm:asymptotic MI}, this Conjecture~\ref{conjecture:dyn} helps us rule out the set $D$ of measure zero.} $\boldsymbol{\lambda}\in(0,\infty)^L$, there exists a unique $\mathbf{q}^\ast$ to the optimization problem \eqref{eq:dyn joint MI} and \eqref{eq:limiting mmse} holds.
        \item[(ii)] If \eqref{eq:threshold_dynamic} fails, $\mathbf{0}$ would be a \textit{global} maximizer to \eqref{eq:dyn joint MI}. Thus, weak recovery is impossible. 
        \item[(iii)] Assume that \eqref{eq:threshold_dynamic} holds and that the initialization $\mathbf{m}^0$ is better than random guessing, i.e.
        \begin{equation*}
            \lim_{n\rightarrow\infty}\frac{1}{n}\mathbb{E}\left[\bX^{(l)\top}\mathbf{m}^0_{\cdot,l}\right]>0
        \end{equation*}
        for some $l\in[L]$. Assume that one implements  Algorithm~\ref{alg:coupled AMP} with denoisers computed from Algorithm~\ref{alg:denoiser dyn}. In this case, the state evolution iterates $\mathbf{q}^t$ converge to $\mathbf{q}^\ast$ exponentially. Thus the algorithm is Bayes optimal and runs in polynomial time.
    \end{itemize}
\end{Proposition}

\begin{figure}
    \centering
    \makebox[\linewidth][c]{
        \includegraphics[width=1.2\linewidth]{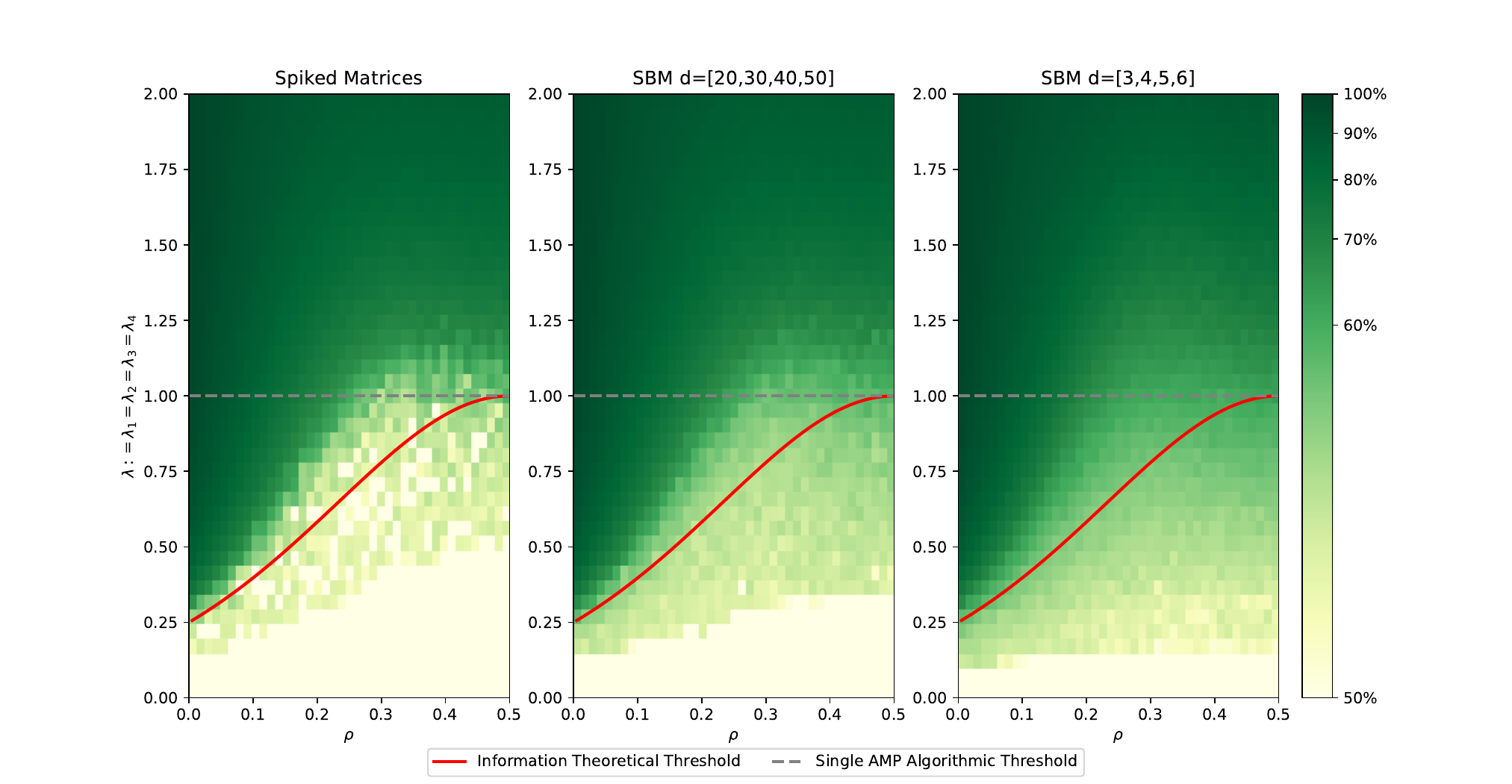}
    }
    \caption{Simulation results on a $4$-epoch dynamical model, formally defined in Example~\ref{eg:dynamic SBM}. Algorithm~\ref{alg:coupled AMP} is tested on spiked matrices, sparse graphs with $d=(20,30,40,50)$ or $d=(3,4,5,6)$ respectively. Color brightness reflects first-moment membership recovery accuracy averaged over $10$ repetitions. The empirical feasibility boundary matches well with the theoretically predicted red curve. Although our theory applies only to denser graphs with growing $d\rightarrow\infty$ as $n$ increases, the simulation suggests the boundary to be critical even for extremely sparse graphs: $d$ is set to $\sim5$ compared to $n=10000$.}
    \label{fig:simu_Dyn}
\end{figure}

\paragraph{Simulation Studies.} We close our discussion on this example with a simulation study on a $4$-epoch dynamical model with varying $\lambda:=\lambda^{(1)}=\lambda^{(2)}=\lambda^{(3)}=\lambda^{(4)}$, as shown in Figure~\ref{fig:simu_Dyn}. Tested models include the spiked matrix counterpart~\eqref{eq:spiked matrix each layer}, a sparse graph setting with average degrees $d=(20,30,40,50)$ and a very sparse setting with $d=(3,4,5,6)$.
We run Algorithm~\ref{alg:coupled AMP} with Bayes optimal denoisers Algorithm~\ref{alg:denoiser dyn} specified in Lemma~\ref{lemma:dyn denoiser} for over $1600$ different choices of $(\rho,\lambda)\in[0,0.5]\times[0,2]$ and $n=10000$ nodes. For each configuration, the algorithm starts with a warm initialization\footnote{About $10\%$ of $\{m_{i,l}^0:i\in[n],l\in[L]\}$ equal the underlying true community membership, with rest set to $0$.} $\mathbf{m}^0$ and proceeds by $100$ iterations to output an estimate $\hat{\mathbf{m}}\in\R^{n\times L}$ for community memberships at all $L=4$ moments. 
For each node, we then pass $\hat{m}_{i}\in\R^L$ to denoiser $\mathcal{E}^{(1)}_{\Dyn,\rho}$ to derive an estimate $\hat{x}_i\in\{\pm 1\}$ on its community membership $X_i^{(1)}$ at the first moment. Lastly the performance is measured by the proportion $\frac{1}{n}|\{i\in[n]:\hat{x}_i=X_i^{(1)}\}|$ of accurately recovered nodes. A random guess should yield $50\%$ accuracy, while our warm initialization yields $10\%/2+50\%=55\%$ accuracy in signal recovery. When below the detection threshold, the algorithm forgets leaked initial information with increasing iterations. Repeating the simulation process $10$ times for each configuration of $(\rho,\lambda)$, the averaged recovery accuracy is reflected by color brightness in Figure~\ref{fig:simu_Dyn}.

\subsection{Semi-supervised Community Detection}
We investigate community recovery in the context of the partially labeled block model \eqref{eg:semi sbm}. Note that in contrast to most prior work, we allow unbalanced partial labelings i.e. the observed labels might be predominantly from a specific community.  


To make resulting formula more compact, we use $\varepsilon=\frac{\varepsilon_++\varepsilon_-}{2}$ and $\delta=\frac{\varepsilon_--\varepsilon_+}{2-2\varepsilon}$ to reparametrize prior $p_\Semi$ in Example~\ref{eg:semi sbm}. Then $\varepsilon$ governs the proportion of labeled vertices, and $\delta$ governs the unbalance of revealed labels in two groups. As a direct application of  Theorem~\ref{thm:asymptotic MI} and Proposition~\ref{prop:Universality MI qualitative} we obtain 
\begin{align}
    \mathcal{F}_{\Semi}(\lambda,q)&=\mathbb{E}\log \left(\sum_{x}p(x|Z)e^{\lambda qXx+\sqrt{\lambda q}W^{\prime}x}\right)\notag\\
    &=\varepsilon\lambda q+\left(1-\varepsilon\right)\mathbb{E}\left\{\log\left[\cosh\left(\lambda qX+\sqrt{\lambda q}W^\prime\right)+\delta \sinh\left(\lambda qX+\sqrt{\lambda q}W^\prime\right)\right]\Big| Z=\ast\right\}\notag\\
    &=\varepsilon\lambda q+\frac{1-\varepsilon}{2}\mathbb{E}\log\left[1+(1-\delta^2)\sinh^2\left(\lambda q+\sqrt{\lambda q}W^\prime\right)\right]\notag\\
    &\qquad\qquad\qquad+\frac{(1-\varepsilon)\delta}{2}\E\log\left[\frac{1+\delta\tanh\left(\lambda q+\sqrt{\lambda q}W^\prime\right)}{1-\delta\tanh\left(\lambda q+\sqrt{\lambda q}W^\prime\right)}\right].
\end{align}
Moreover, the mutual information between $X$ and $Z$ under $p_\Semi$ can be computed by
\begin{equation*}
    i_p(x;z)=\varepsilon\log2 +\frac{1-\varepsilon}{2}\left[(1+\delta)\log(1+\delta)+(1-\delta)\log(1-\delta)\right].
\end{equation*}

\begin{Theorem}\label{thm:semi joint MI}
    When the average degree $d_n = \frac{a +b}{2} $ satisfies $d_n(1-d_n/n) \to \infty$, 
    \begin{align}
        &\lim_{n\rightarrow\infty}\frac{1}{n}I(\bX;\mathbf{G},\mathbf{Z})=\lim_{n\rightarrow\infty}\frac{1}{n}I(\bX;\bA,\mathbf{Z})
        =i_p(x;z)+\frac{\lambda}{4}-\sup_{q\ge 0}\left[\mathcal{F}_{\Semi}(\lambda,q)-\frac{\lambda(q^2+2q)}{4}\right],\label{eq:semi joint MI}
    \end{align}
    where we recall $\mathbf{A}$ to be the approximating spiked matrices in \eqref{eq:spiked matrix each layer}. 
\end{Theorem}
Figure~\ref{fig:mi semi} below displays the curve $\lambda\mapsto \lim I(\bX;\bG,\bZ)/n$ for $4$ different choices of $(\varepsilon_+,\varepsilon_-)$.
We next turn our attention to the algorithmic question of combining the network information with the partial vertex labels. We note that this algorithmic question has motivated significant research in the existing literature (see e.g. \cite{kanade2016global,cai2020weighted,zhang2014phase,saade2018fast} and references therein). For sparse graphs, the main proposal has been to use belief propagation style algorithms (or associated spectral algorithms based on the nonbacktracking operator) to combine the graph information with the revealed vertex labels. Here we instead employ an algorithm based on the AMP formalism. The main advantage of AMP based algorithms over belief propagation stems from their tractability---in this case, we can rigorously analyze the signal recovery performance of this algorithm using the state evolution framework. 

\paragraph{AMP Denoiser and State Evolution.}
To implement our Algorithm~\ref{alg:coupled AMP} in this setting, we have to specify the non-linear denoisers. 
We choose the Bayes optimal denoisers 
\begin{align}
    \mathcal{E}_{\Semi}(m,z)&=\mathbb{E}\left[X\bigg| \sqrt{\lambda q^t}A^{\prime}=m,Z=z\right]\notag\\
    &=z\mathbbm{1}\{z\neq\ast\}+\frac{\tanh(m)+\delta}{1+\delta\tanh(m)}\mathbbm{1}\{z=\ast\}.\label{eq:semi denoiser}
\end{align}
In turn, the state evolution parameters $\{q^t:t\in\mathbb{N}\}$ evolve as 
\begin{align*}
    q^{t+1}&=\mathbb{E}\left[X\mathcal{E}_{\Semi}(\lambda q^t X+\sqrt{\lambda q^t}W^\prime,Z)\right]\\
    &=\varepsilon+(1-\varepsilon)\mathbb{E}\left[X\frac{\tanh(\lambda q^tX+\sqrt{\lambda q^t}W^\prime)+\delta}{1+\delta\tanh(\lambda q^tX+\sqrt{\lambda q^t}W^\prime)}\Bigg|Z=\ast\right]\\
    &=\varepsilon+(1-\varepsilon)\mathbb{E}\left[\frac{(1-\delta^2)\tanh(\lambda q^t+\sqrt{\lambda q^t}W^\prime)+\delta^2\left(1-\tanh^2(\lambda q^t+\sqrt{\lambda q^t}W^\prime)\right)}{1-\delta^2\tanh^2(\lambda q^t+\sqrt{\lambda q^t}W^\prime)}\right]=T^{\Semi}(\lambda q^t),
\end{align*}
where $T^\Semi$ is given in \eqref{eq:general se mapping} by specifying prior $p_\Semi$. By recursively updating $q^{t+1}\leftarrow T^{\Semi}(\lambda q^t)$, the state evolution iterates $q^t$ usually converge to a fixed point of mapping $q\mapsto T^{\Semi}(\lambda q)$. In addition, if this fixed point is also the unique global maximizer of the free energy functional, our coupled AMP algorithm is indeed Bayes optimal.

\paragraph{Optimal Recovery}
As shown in Figure~\ref{fig:se mapping semi}, the curve $\gamma\mapsto T^\Semi(\gamma)$ is increasing and seemingly concave. Therefore, for each $\lambda>0$, the straight line $\gamma/\lambda$ has a unique intersection with this curve, thus justifying Bayes optimality of our AMP algorithm. 

Although the monotonicity of $T^\Semi(\gamma)$ should hold generally, its concavity might not hold in general. When $\varepsilon_+=\varepsilon_-$ (i.e. $\delta=0$), \cite{deshpande2017asymptotic} rigorously established the concavity of this mapping. As demonstrated in \cite{montanari2021estimation} and Figure 1 of \cite{guo2005mutual}, we anticipate $T^\Semi(\gamma)$ to be convex near $\gamma\approx 0$ if $|\delta|\to 1$. 
On the information-theoretic side, the curvature of $T^\Semi(\gamma)$ determines the optimal recovery performance achievable by an algorithm.  
On an algorithmic side, this issue also determines the convergence of state evolution iterates, thus deciding practical feasibility of optimal statistical  recovery. We leave a careful study of this phenomenon for future work. 

\begin{figure}
    \begin{subfigure}{0.5\textwidth}
        \centering
        \includegraphics[width=\linewidth]{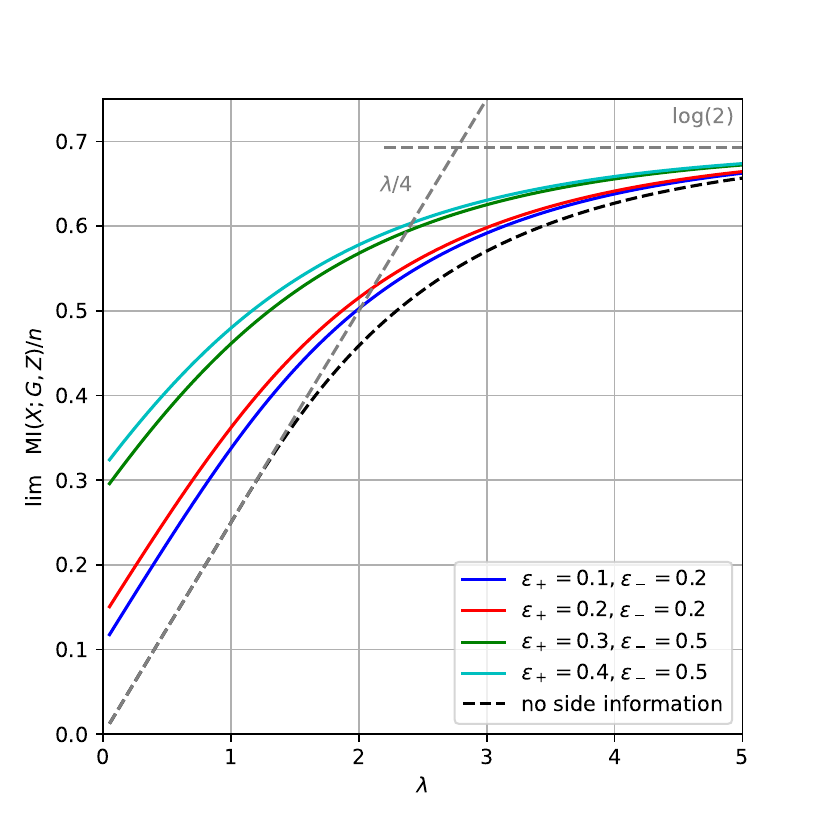}
        \caption{Limiting Mutual Information}
        \label{fig:mi semi}
    \end{subfigure}
    \begin{subfigure}{0.5\textwidth}
        \centering
        \includegraphics[width=\linewidth]{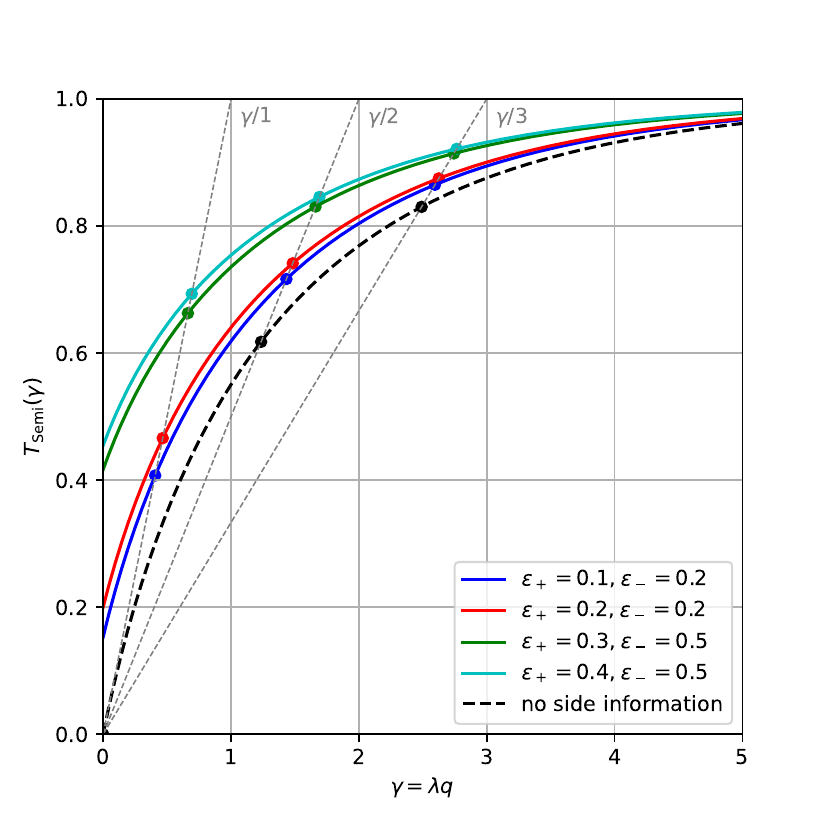}
        \caption{State Evolution Mapping}
        \label{fig:se mapping semi}
    \end{subfigure}
    \caption{Illustration of $\lim I(\bX;\bG,\bZ)/n$ and $T_\Semi$ for SBM with Partially Observed Labels Example~\ref{eg:semi sbm}. For each $\lambda$, the intersection of $\gamma\mapsto T_\Semi(\gamma)$ and straight line $\gamma/\lambda$ measures the limiting overlap between $\bX$ and $\E(\bX|\bG,\bZ)$, as illustrated by colored dots in (b).}
\end{figure}

\vspace{2mm}
\noindent
\textbf{Organization:} The rest of the paper is structured as follows. We establish Proposition~\ref{prop:RS prediction on free energy} in Section~\ref{sec:RS Prediction of Free Energy}. Armed with this result, we turn to a proof of Theorem~\ref{thm:asymptotic MI}, Proposition~\ref{prop:Universality MI qualitative} and Proposition~\ref{prop:universality mmse} in Section~\ref{sec:mi_mmse}. Section~\ref{sec:AMP_derivation} focuses on the derivation and analysis of the general Approximate Message Passing Algorithm (Algorithm \ref{alg:coupled AMP}). Finally, we establish Proposition~\ref{prop:opt derivative}, Proposition~\ref{prop:all-in-all-out} and the results specific to the examples in Section~\ref{sec:misc_proofs}. 

\vspace{2mm}
\noindent
\textbf{Notation:} Usual Bachmann-Landau notation $O(\cdot),o(\cdot),\Theta(\cdot)$ is used throughout the manuscript. We use $C$ to denote positive constants, and its dependence may differ across positions of appearance. For any matrix $A$, $\|A\|$ denotes its operator norm and $A_{i,j}$ denotes a certain entry of it. Throughout the paper, $L$ always represents the number of networks we get to observe, and $n$ is the number of vertices.
Capital English letters $X,Y,Z,A,W,\ldots$ denote random variables whose dimension is finite (it could dependent on $L$, but remain fixed as $n$ increases); while their boldface version $\bX,\bY,\bZ,\bA,\bW,\ldots$ denotes $n$ or $n^2$ repetitions depending on the context. Lowercase letters $x,y,z,\ldots$ denote corresponding realizations of these random variables. We use superscripts $X^{(l)}$ to refer to different layers and subscripts $X_i$ to count though $n$ units. Finally, we write $\overset{\text{P}}{\rightarrow}$ for convergence in probability and $a\odot b$ for coordinate-wise product of two vectors $a,b\in\R^L$.


\section{Replica Symmetry of the limiting Free Energy}\label{sec:RS Prediction of Free Energy}

We establish Proposition~\ref{prop:RS prediction on free energy} in this section. We establish the lower bound using Guerra's interpolation technique, while the upper bound is established using the Aizenman-Simms-Starr scheme. This technique originated in the study of mean-field spin glasses (see e.g. \cite{panchenko2013sherrington}). In the special setting of inference problems, our approach mirrors the one introduced in the seminal work of \cite{lelarge2019fundamental}. We note that alternative techniques e.g. the adaptive interpolation techniques have been developed to derive limiting mutual information in Bayesian inference problems \cite{barbier2019adaptive,barbier2019adaptivelater}. One could, in principle, approach our problem using these alternative techniques. 

\subsection{Lower Bound: Guerra's Interpolation Technique}\label{subsec:Guerra's interpolation}

Fix $q_l\geq 0$, and introduce $t\in[0,1]$ to interpolate between the spiked matrix model and white noise scalar channel,
\begin{equation}
    \left\{\begin{aligned}
        &(X_i,Y_i,Z_i)\indsim p(\cdot), &1\le i\le n;\\
        &A_{i,j}^{(l)}=\sqrt{\frac{\lambda^{(l)} t}{n}}X_i^{(l)}X_j^{(l)}+W_{i,j}^{(l)}, &1\leq i<j\leq n,1\le l\le L;\\
        &A^{(l)\prime}_i=\sqrt{(1-t)\lambda^{(l)} q_l}X^{(l)}_i+W^{(l)\prime}_i,&1\leq i\leq n,1\le l\le L.
    \end{aligned}\right.
\end{equation}
We define an intermediate Hamiltonian for any $\bx\in\{\pm 1\}^{n\times L}$,
\begin{equation}
    \begin{aligned}
        H_n(\mathbf{x},t)=&\sum_{l=1}^L\sum_{i<j}\sqrt{\frac{\lambda^{(l)} t}{n}}W_{i,j}^{(l)}x_i^{(l)}x_j^{(l)}+\frac{\lambda^{(l)} t}{n}X_i^{(l)}X_j^{(l)}x_i^{(l)}x_j^{(l)}\\
        &+\sum_{l=1}^L\sum_{i=1}^n\sqrt{(1-t)\lambda^{(l)} q_l}W^{(l)\prime}_i x_i^{(l)}+(1-t)\lambda^{(l)} q_lX_i^{(l)}x_i^{(l)}.
    \end{aligned}
\end{equation}
We introduce the corresponding  intermediate free energy at time $t\in[0,1]$ as
\begin{equation}
    \psi(t)=\frac{1}{n}\mathbb{E}\left\{\log\left[\sum_{\bx,\by}\prod_{i=1}^np(x_i,y_i|Z_i)\cdot \exp\left(H_n(\mathbf{x},t)\right)\right]\right\}.
\end{equation}
With this definition, one can directly verify that $\psi(1)= F_n(\boldsymbol{\lambda})$ and $\psi(0)$ corresponds exactly to the scalar channel introduced in \eqref{eq:scalar channel}. Differentiating in $t$, we obtain 
\begin{align}
    \psi^\prime(t)=&\frac{1}{n}\mathbb{E}\left\langle\frac{\partial}{\partial t}H_n(\mathbf{x},t)\right\rangle\\
    =&\frac{1}{n}\mathbb{E}\left\langle\sum_{l=1}^L\sum_{i<j}\frac{1}{2}\sqrt{\frac{\lambda^{(l)}}{nt}}W_{i,j}^{(l)}x_i^{(l)}x_j^{(l)}+\frac{\lambda^{(l)}}{n}X_i^{(l)}X_j^{(l)}x_i^{(l)}x_j^{(l)}\right\rangle\notag\\
    &-\frac{1}{n}\mathbb{E}\left\langle\sum_{l=1}^L\sum_{i=1}^n\frac{1}{2}\sqrt{\frac{\lambda^{(l)} q_l}{1-t}}W^{(l)\prime}_i x_i^{(l)}+\lambda^{(l)} q_lX_i^{(l)}x_i^{(l)}\right\rangle.\label{eq:derivative of intermediate free energy}
\end{align}
Gaussian integration by parts  yields
\begin{align*}
    \mathbb{E}W_{i,j}^{(l)}\left\langle x_i^{(l)}x_j^{(l)}\right\rangle=\mathbb{E}\frac{\partial}{\partial W_{i,j}^{(l)}}\left\langle x_i^{(l)}x_j^{(l)}\right\rangle=\sqrt{\frac{\lambda^{(l)}t}{n}}\mathbb{E}\left[\left\langle x_i^{(l)2}x_j^{(l)2}\right\rangle-\left\langle x_i^{(l)}x_j^{(l)}\right\rangle^2\right].
\end{align*}
Let $(\bx,\by)$ and $(\bar{\bx},\bar{\by})$ be two independent configurations drawn from the Gibbs distribution. The first term is trivial since $x_i^{(l)2}=1$.
The second term above can be re-expressed as 
\begin{equation*}
    \mathbb{E}\left[\left\langle x_i^{(l)}x_j^{(l)}\right\rangle^2\right]=\mathbb{E}\left[\left\langle x_i^{(l)}x_j^{(l)}\bar{x}_i^{(l)}\bar{x}_j^{(l)}\right\rangle\right]=\mathbb{E}\left[\left\langle x_i^{(l)}x_j^{(l)}X_i^{(l)}X_j^{(l)}\right\rangle\right],
\end{equation*}
where the last equality follows by Nishimori identity. Therefore,
\begin{align*}
    \mathbb{E}W_{i,j}^{(l)}\left\langle x_i^{(l)}x_j^{(l)}\right\rangle=&\sqrt{\frac{\lambda^{(l)}t}{n}}\left[1-\mathbb{E}\left\langle X_i^{(l)}X_j^{(l)}x_i^{(l)}x_j^{(l)}\right\rangle\right].
\end{align*}
Similar computations yield
\begin{align*}
    \mathbb{E}W_{i}^{(l)\prime}\left\langle y_i^{(l)}\right\rangle=&\sqrt{\lambda^{(l)}q_l(1-t)}\left[1-\mathbb{E}\left\langle X_i^{(l)}x_i^{(l)}\right\rangle\right].
\end{align*}
Plugging back into \eqref{eq:derivative of intermediate free energy}, we can provide a lower bound
\begin{align}
    \psi^\prime(t)=&\frac{1}{2n}\mathbb{E}\left[\left\langle\sum_{l=1}^L\sum_{i<j}\frac{\lambda^{(l)}}{n}X_i^{(l)}X_j^{(l)}x_i^{(l)}x_j^{(l)}-\sum_{l=1}^L\sum_{i=1}^{n}\lambda^{(l)} q_lX_i^{(l)}x_i^{(l)}\right\rangle\right]+\sum_{l=1}^L \lambdal \left(\frac{n-1}{4n}-\frac{q_l}{2}\right)\notag\\
    =&\frac{1}{4}\sum_{l=1}^L\mathbb{E}\left[\lambda^{(l)}\left\langle\left(\bx^{(l)}.\bX^{(l)}\right)^2-2q_l\bx^{(l)}.\bX^{(l)}\right\rangle\right]+\sum_{l=1}^L \frac{\lambdal(1-2q_l)}{4}+o(1)\notag\\
    =&\frac{1}{4}\sum_{l=1}^L\mathbb{E}\left[\lambda^{(l)}\left\langle\left(\bx^{(l)}.\bX^{(l)}-q_l\right)^2\right\rangle\right]+\sum_{l=1}^L \frac{\lambdal(1-2q_l-q_l^2)}{4}+o(1)\notag\\
    \ge&\sum_{l=1}^L \frac{\lambdal(1-2q_l-q_l^2)}{4}+o(1).\label{eq:derivative of intermediate free energy 2}
\end{align}
where we define $\bx^{(l)}.\bX^{(l)}=\frac{1}{n}\sum_{i=1}^n x_i^{(l)}X_i^{(l)}$ as the empirical overlap.
Because $t=0$ corresponds to repeating channel \eqref{eq:scalar channel} independently for $n$ times,
\begin{equation*}
    \psi(0)=\mathbb{E}\log\left[\sum_{(x,y)\in\{\pm 1\}^{L+L_1}}p(x,y|Z)\exp\left(\sum_{l=1}^L\sqrt{\lambda^{(l)}q_l}A^{(l)\prime}x^{(l)}\right)\right].
\end{equation*} 
Since $F_n(\boldsymbol{\lambda})=\psi(0)+\int_{0}^1\psi^\prime(t)\mathrm{d}t$ and $\psi^\prime(t)$ is lower bounded in \eqref{eq:derivative of intermediate free energy 2}, we find
\begin{equation}\label{eq:replica symmetric lower bound}
    \lim_{n\rightarrow\infty}F_n(\boldsymbol{\lambda})\ge\frac{\sum_{l=1}^{L}\lambdal}{4}+\sup_{\mathbf{q}\ge 0}\left[\mathcal{F}(\boldsymbol{\lambda},\mathbf{q})-\sum_{l=1}^{L}\frac{\lambda^{(l)}(q_l^2+2q_l)}{4}\right].
\end{equation}

\subsection{Overlap Concentration via Perturbation by BEC}\label{subsec:overlap concentration with perturb}

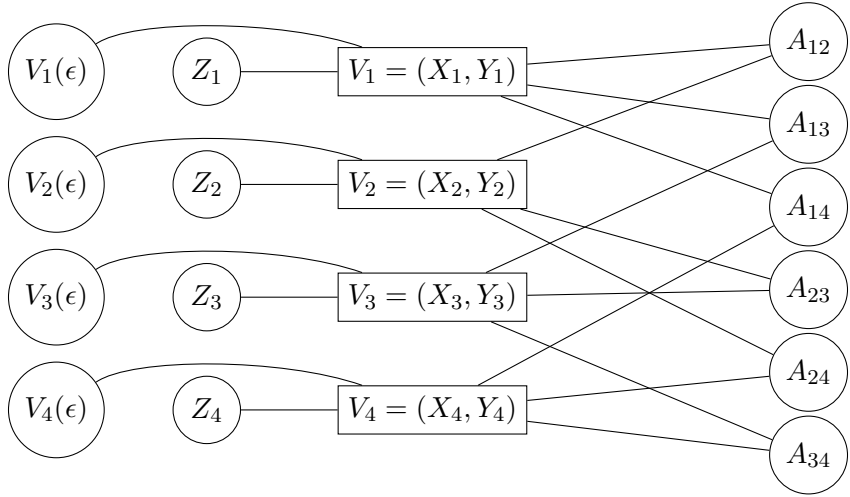
\begin{figure}
\centering
\begin{tikzpicture}
\centering
  \pgfmathsetmacro{\n}{4}
  \pgfmathsetmacro{\m}{6}
  \foreach \i in {1,...,\n} {
    \node[circle, draw] (Xcircle\i) at (-5,-\i*3/2) {$V_{\i}(\epsilon)$};
    \node[circle, draw] (Zcircle\i) at (-3,-\i*3/2) {$Z_{\i}$};
  }
  \foreach \i in {1,...,\n} {
    \node[rectangle, draw] (square\i) at (0,-\i*3/2) {$V_{\i}=(X_{\i},Y_{\i})$};
  }
  
  \pgfmathsetmacro{\np}{\n-1}
  \foreach \x in {1,...,\np} {
    \pgfmathsetmacro{\xp}{int(\x + 1)}
    \foreach \y in {\xp,...,\n} {
      \pgfmathsetmacro{\j}{\n*(\x-1) -int(\x*(\x+1)/2) + \y}
      \node[circle, draw, minimum size=1cm] (circle\x\y) at (5,-\j*1.1) {$A_{\x \y}$};
    }
  }

  \foreach \i in {1,...,\n} {
    \draw (square\i) -- (Zcircle\i);
    \draw (Xcircle\i) .. controls (-4,-\i*3/2+0.7) and (-2,-\i*3/2+0.7) .. (square\i);
  }
  \foreach \x in {1,...,\np} {
    \pgfmathsetmacro{\xp}{int(\x+1)}
    \foreach \y in {\xp,...,\n} {
      \draw (square\x) -- (circle\x\y);
      \draw (square\y) -- (circle\x\y);
    }
  }
\end{tikzpicture}
\caption{Factor graph of the observation model when $n=4$.}\label{fig:factor graph illustration}
\end{figure}

For the converse upper bound, it will be useful to establish concentration of the overlap in this inference problem. To this end, inspired by the techniques introduced in \cite{andrea2008estimating,lelarge2019fundamental}, we add a small perturbation to the model via a binary erasure channel. This will force correlation decay among the different coordinates in a sample. To set the stage of the small perturbation, we introduce a factor graph $(V,F,Z,E)$ of our observation model, which is defined to be a bipartite graph.
To be precise, $|V|=n$ are vertices corresponding to variable nodes, $V_i=(X_i,Y_i)$ for each $i\in[n]$; $|F|=n(n-1)/2$ corresponds to function nodes standing for memoryless observations, $F_a=(A_{i,j}^{(1)},\ldots,A_{i,j}^{(L)})$ for each pair $a=(i,j)$ with $1\leq i<j\leq n$. Edge set $E$ characterizes how each $F_a$ depends on $V$, namely $F_{(i,j)}$ only connects to $V_i$ and $V_j$.
Moreover, we also add nodes $Z_i$ corresponding to additionally observed side information. Each $Z_i$ connects to $V_i$ only. See Figure~\ref{fig:factor graph illustration} for an illustration of this factor graph, where circles indicate observed variables and squares indicate variables to be recovered.

\vspace{3mm}
\noindent\textit{Remark:} Compared to the sparse graph in \cite{andrea2008estimating}, our bipartite graph has asignificantly larger number of function nodes  ($n(n-1)/2$ compared to the number $n$ of variable nodes). However, each function node's signal strength $\frac{\boldsymbol{\lambda}}{n}$ decays with respect to $n$, thus facilitating applications of lemmas from \cite{andrea2008estimating}.
\vspace{3mm}

 We are now ready to formally define the perturbation in the form of erasure channels. Suppose we additionally observe,
\begin{equation}\label{eq:perturbation def}
    V_i(\epsilon)=\mathsf{BEC}_\epsilon(V_i)=\left\{\begin{aligned}
        & V_i&\text{with probability } \epsilon,\\
        & \ast&\text{otherwise}.
    \end{aligned}\right.
\end{equation}
Upon observing this additional information, we first justify that the perturbed free energy does not deviate too much from the original free energy. For any configuration $\bv=(v_1,\ldots,v_n)\in\{\pm 1\}^{n\times (L+L_1)}$, define $\widetilde{\bv}=(\widetilde{v}_1,\ldots,\widetilde{v}_n)$ by
\begin{equation*}
    \widetilde{v}_i=v_i\mathbbm{1}\{V_i(\epsilon)=\ast\}+V_i\mathbbm{1}\{V_i(\epsilon)\neq\ast\}.
\end{equation*}
This notion allows a very convenient expression for the normalizing constant $Z_{n,\epsilon}$ of $p(\cdot|\mathbf{V}(\epsilon),\bZ,\mathbf{A})$,
\begin{align*}
    Z_{n,\epsilon}:=&\sum_{\bv}p(\bv|\mathbf{V}(\epsilon),\bZ)e^{H_n(\bv)}=\sum_{\bv}\left(\prod_{i:V_i(\epsilon)\neq\ast}\mathbbm{1}\{v_i=V_i\}\right)\left(\prod_{i:V_i(\epsilon)=\ast}p(v_i|Z_i)\right)e^{H_n(\bv)}.
\end{align*}
Denote $\mathcal{I}=\{i\in[n]:V_i(\epsilon)=\ast\}$, then we find
\begin{align*}
    Z_{n,\epsilon}=&\sum_{\bv_{\mathcal{I}}}\left(\prod_{i\in\mathcal{I}}p(v_i|Z_i)\right)e^{H_n(\bv_{\mathcal{I}}\oplus \bV_{\mathcal{I}^c})}\\
    =&\sum_{\bv_{\mathcal{I}}}\left(\sum_{\bv_{\mathcal{I}^c}}\prod_{i\notin\mathcal{I}}p(v_i|Z_i)\right)\left(\prod_{i\in\mathcal{I}}p(v_i|Z_i)\right)e^{H_n(\bv_{\mathcal{I}}, \bV_{\mathcal{I}^c})}\\
    =&\sum_{\bv}\prod_{i}p(v_i|Z_i)e^{H_n(\widetilde{\bv})}.
\end{align*}
\begin{Proposition}\label{prop:get rid of perturbation}
    Define pertubed normalized free energy as $F_{n,\epsilon}=\frac{1}{n}\mathbb{E}\log Z_{n,\epsilon}$. This quantity is Lipschitz in $\epsilon$, i.e. there exists $C>0$ (independent of $n$) such that for any $\epsilon,\epsilon^\prime\in[0,1]$, we have
    \begin{equation*}
        |F_{n,\epsilon}-F_{n,\epsilon^\prime}|\leq C|\epsilon-\epsilon^\prime|.
    \end{equation*}
\end{Proposition}
\begin{proof}
    To establish this result, we consider a more general setting in which revelation probabilities for each $V_i$ can be different: $\mathbb{P}(V_i(\boldsymbol{\epsilon})=V_i)=\epsilon_i$. We can still define $F_{n,\boldsymbol{\epsilon}}$ and $Z_{n,\boldsymbol{\epsilon}}$. Using the law of total probability, 
    \begin{align*}
        F_{n,\boldsymbol{\epsilon}}=&\frac{\epsilon_1}{n}\mathbb{E}\left[\log\left(\sum_{\bv}p(\bv|\bZ)e^{H_n(\widetilde{\bv})}\right)\bigg|V_1(\boldsymbol{\epsilon})\neq\ast\right]+\frac{1-\epsilon_1}{n}\mathbb{E}\left[\log\left(\sum_{\bv}p(\bv|\bZ)e^{H_n(\widetilde{\bv})}\right)\bigg|V_1(\boldsymbol{\epsilon})=\ast\right]\\
        =&\frac{\epsilon_1}{n}\mathbb{E}\left[\log\left(\sum_{\bv}p(\bv|\bZ)e^{H_n(V_1,\widetilde{\bv}_{-1})}\right)\right]+\frac{1-\epsilon_1}{n}\mathbb{E}\left[\log\left(\sum_{\bv}p(\bv|\bZ)e^{H_n(v_1,\widetilde{\bv}_{-1})}\right)\right].
    \end{align*}
    Consequently,
    \begin{equation*}
        \frac{\partial F_{n,\boldsymbol{\epsilon}}}{\partial\epsilon_1}=\frac{1}{n}\mathbb{E}\left[\log\left(\sum_{\bv}p(\bv|\bZ)e^{H_n(V_1,\widetilde{\bv}_{-1})}\right)\right]-\frac{1}{n}\mathbb{E}\left[\log\left(\sum_{\bv}p(\bv|\bZ)e^{H_n(v_1,\widetilde{\bv}_{-1})}\right)\right].
    \end{equation*}
    Define a Hamiltonian $H^\prime_n(\bv)=\sum_{2\le i<j}\sum_{l=1}^L \sqrt{\frac{\lambdal}{n}}A^{(l)}_{i,j}x_i^{(l)}x_j^{(l)}$ and let $\langle\cdot\rangle$ denote the expectation with respect to the Gibbs measure corresponding to $H_n^\prime$. Then we can write
    \begin{align*}
        \frac{\partial F_{n,\boldsymbol{\epsilon}}}{\partial\epsilon_1}=&\frac{1}{n}\mathbb{E}\left[\log\left\langle\exp\left(\sum_{j\ge 2}\sum_{l=1}^L \sqrt{\frac{\lambdal}{n}}A^{(l)}_{1,j}X_1^{(l)}\widetilde{x}_j^{(l)}\right)\right\rangle\right]\\
        &-\frac{1}{n}\mathbb{E}\left[\log\left\langle\sum_{v_1}p(v_1|Z_1)\exp\left(\sum_{j\ge 2}\sum_{l=1}^L \sqrt{\frac{\lambdal}{n}}A^{(l)}_{1,j}x_1^{(l)}\widetilde{x}_j^{(l)}\right)\right\rangle\right]=:T_1-T_2.
    \end{align*}
    Let $\mathbb{E}_1$ represent taking expectation over $(W_{1,j})_{j\ge 2}$ only. Since $(W_{1,j})_{j\ge 2}$ is independent of $H^\prime_n$ defined above, $\mathbb{E}_1$ can be alternated with $\langle\cdot\rangle$. In addition, applying Jenson's inequality on $\log(\cdot)$, we have
    \begin{align*}
        T_1\le &\frac{1}{n}\mathbb{E}\left[\log\left\langle\mathbb{E}_1\exp\left(\sum_{j\ge 2}\sum_{l=1}^L \frac{\lambdal}{n}X_1^{(l)2}X_j^{(l)}\widetilde{x}_j^{(l)}+\sqrt{\frac{\lambdal}{n}}W^{(l)}_{1,j}X_1^{(l)}\widetilde{x}_j^{(l)}\right)\right\rangle\right]\\
        =&\frac{1}{n}\mathbb{E}\left[\log\left\langle\exp\left(\sum_{j\ge 2}\sum_{l=1}^L \frac{\lambdal}{n}X_j^{(l)}\widetilde{x}_j^{(l)}+\frac{\lambdal}{2n}\right)\right\rangle\right]\le \frac{3\sum_{l}\lambdal}{2n}.
    \end{align*}
    On the other hand, applying  Jenson's inequality to the $\exp(\cdot)$ function, we have
    \begin{align*}
        T_1\ge &\frac{1}{n}\mathbb{E}\left[\sum_{j\ge 2}\sum_{l=1}^L \frac{\lambdal}{n}X_1^{(l)2}X_j^{(l)}\left\langle\widetilde{x}_j^{(l)}\right\rangle+\sqrt{\frac{\lambdal}{n}}W^{(l)}_{1,j}X_1^{(l)}\left\langle\widetilde{x}_j^{(l)}\right\rangle\right]\\
        =&\frac{1}{n}\mathbb{E}\left[\sum_{j\ge 2}\sum_{l=1}^L \frac{\lambdal}{n}X_j^{(l)}\left\langle\widetilde{x}_j^{(l)}\right\rangle\right]\ge \frac{\sum_{l}\lambdal}{n}.
    \end{align*}
    In the same way, we can also show $|T_2|$ is of order $O(1/n)$. As a result, we would have $\left|\frac{\partial F_{n,\boldsymbol{\epsilon}}}{\partial\epsilon_1}\right|\le \frac{C}{n}$. We set $\epsilon_1=\cdots=\epsilon_n$ to conclude that $\left|\frac{\partial F_{n,\epsilon}}{\partial\epsilon}\right|\le C$, thus finishing the proof.
\end{proof}


This perturbation ensures decorrelation in the posterior distribution, as established in 
\cite[Lemma 3.1]{andrea2008estimating}. We recall the result for the sake of completeness. 
\begin{Lemma}[\cite{andrea2008estimating}]\label{lem:vanishing mutual information between nodes}
    Let $\epsilon\sim U[0,n^{-1/2}]$ independent of everything else, the expected conditional mutual information between $V_i$ and $V_j$ decays in n. Formally, we have, 
    \begin{equation}
        \frac{1}{\sqrt{n}} \mathbb{E}_\epsilon\left[\frac{1}{n^2}\sum_{i,j}I(V_i;V_j|\mathbf{V}(\epsilon),\bZ,\mathbf{A})\right]\leq\frac{C}{n}.
    \end{equation} 
\end{Lemma}

We use this correlation decay to establish overlap concentration for the perturbed posterior. To this end, we follow the strategy introduced in \cite[Proposition 24]{lelarge2019fundamental}. 
\begin{Lemma}\label{lemma:overlap concentration after perturbation}
    Let $\mathbf{v}=(v_1,\ldots,v_n)$ and $\bar{\mathbf{v}}=(\bar{v}_1,\ldots,\bar{v}_n)$ denote  two independent configurations sampled from the posterior conditioned on $(\mathbf{V}(\epsilon),\bZ,\bA)$, with each $v_i=\left(x_i,y_i\right)$ and $\bar{v}_i=\left(\bar{x}_i,\bar{y}_i\right)$. Denote the empirical overlap and its expectation (under the posterior) respectively as 
    \begin{equation*}
        \mathbf{v}.\bar{\mathbf{v}}=\frac{1}{n}\sum_{i=1}^n v_i\bar{v}_i^\top, \quad Q=\left\langle\frac{1}{n}\sum_{i=1}^n v_i\bar{v}_i^\top\right\rangle.
    \end{equation*}
    Under the above mentioned perturbation,
    \begin{equation}
        \mathbb{E}_\epsilon\mathbb{E}\left[\left\langle\left\|\mathbf{v}.\bar{\mathbf{v}}-Q\right\|^2\right\rangle\right]\underset{n\rightarrow \infty }{\longrightarrow}0.
    \end{equation}
\end{Lemma}
\begin{proof}
    First, we expand the squared residual,
    \begin{align*}
        &\left\langle\left\|\mathbf{v}.\bar{\mathbf{v}}-Q\right\|^2\right\rangle=\left\langle\left\|\mathbf{v}.\bar{\mathbf{v}}\right\|^2\right\rangle-\left\|\left\langle\mathbf{v}.\bar{\mathbf{v}}\right\rangle\right\|^2\\
        =&\left\langle\frac{1}{n^2}\sum_{i,j}\Tr\left(v_i\bar{v}_i^\top v_j\bar{v}_j^\top\right)\right\rangle-\Tr\left(\left\langle \frac{1}{n}\sum_{i}v_i\bar{v}_i^\top\right\rangle\left\langle \frac{1}{n}\sum_{j}v_j\bar{v}_j^\top\right\rangle\right)\\
        =&\frac{1}{n^2}\sum_{i,j}\Tr\left[\left\langle v_iv_j^{\top}\right\rangle^2-\left(\left\langle v_i\right\rangle\left\langle v_j\right\rangle^\top\right)^2\right].
    \end{align*}
    Since every $v_i$ is bounded and $\Tr(A^2-B^2)=\Tr((A-B)(A+B))$ for square matrices,
    \begin{align*}
        &\left\langle\left\|\mathbf{v}.\bar{\mathbf{v}}-Q\right\|^2\right\rangle\\
        \leq&\frac{C}{n^2}\sum_{i,j}\left\|\left\langle v_iv_j^{\top}\right\rangle-\left\langle v_i\right\rangle\left\langle v_j\right\rangle^\top\right\|\\
        =&\frac{C}{n^2}\sum_{i,j}\left\|\sum_{v_i,v_j}v_iv_j^\top \left[p(V_i=v_i,V_j=v_j\big|\mathbf{V}(\epsilon),\bZ,\bA)-p(V_i=v_i\big|\mathbf{V}(\epsilon),\bZ,\bA)p(V_j=v_j\big|\mathbf{V}(\epsilon),\bZ,\bA)\right]\right\|\\
        \leq&\frac{C}{n^2}\sum_{i,j}\sum_{v_i,v_j}\left|p(V_i=v_i,V_j=v_j\big|\mathbf{V}(\epsilon),\bZ,\bA)-p(V_i=v_i\big|\mathbf{V}(\epsilon),\bZ,\bA)p(V_j=v_j\big|\mathbf{V}(\epsilon),\bZ,\bA)\right|\\
        \leq&\frac{C}{n^2}\sum_{i,j}\text{D}_{\text{TV}}(p(V_i,V_j\big|\mathbf{V}(\epsilon),\bZ,\bA),p(V_i\big|\mathbf{V}(\epsilon),\bZ,\bA)\otimes p(V_j\big|\mathbf{V}(\epsilon),\bZ,\bA)).
    \end{align*}
    Applying Pinsker's inequality,
    \begin{align*}
        &\left\langle\left\|\mathbf{v}.\bar{\mathbf{v}}-Q\right\|^2\right\rangle\\
        \leq&\frac{C}{n^2}\sum_{i,j}\text{D}_{\text{TV}}(p(V_i,V_j\big|\mathbf{V}(\epsilon),\bZ,\bA),p(V_i\big|\mathbf{V}(\epsilon),\bZ,\bA)\otimes p(V_j\big|\mathbf{V}(\epsilon),\bZ,\bA))\\
        \leq&\frac{C}{n^2}\sum_{i,j}I\left(V_i;V_j|\mathbf{V}(\epsilon),\bZ,\bA\right).
    \end{align*}
    To finish this lemma, it suffices to take expectation and apply Lemma~\ref{lem:vanishing mutual information between nodes}.
\end{proof}
Finally, using Nishimori identity we have
\begin{equation}
    \mathbb{E}_\epsilon\mathbb{E}\left[\left\langle\left\|\mathbf{v}.\mathbf{V}-Q\right\|^2\right\rangle\right]\underset{n\rightarrow \infty }{\longrightarrow}0.
\end{equation}

\subsection{Upper Bound: Aizenman-Sims-Starr scheme}\label{subsec:Aizenman-Sims-Starr}
With the convention $F_{0,\epsilon_0}=0$, we have $F_{n,\epsilon_n}=\frac{1}{n}\sum_{k=0}^{n-1}(k+1)F_{k+1,\epsilon_{k+1}}-kF_{k,\epsilon_k}$. Define
\begin{align}
    B_n=&(n+1)F_{n+1,\epsilon_{n+1}}-nF_{n,\epsilon_n}=\mathbb{E}\log Z_{n+1,\epsilon_{n+1}}-\mathbb{E}\log Z_{n,\epsilon_{n}}\notag\\
    =&\mathbb{E}\log\left(\sum_{\bv,v_{n+1}}p(\bv|\bZ)p(v_{n+1}|Z_{n+1})e^{H_{n+1}(\widetilde{\bv},\widetilde{v}_{n+1})}\right)-\mathbb{E}\log\left(\sum_{\bv}p(\bv|\bZ)e^{H_{n}(\widetilde{\bv})}\right),\label{eq:Cesaro term}
\end{align}
in which $\bv$ is the collection of first $n$ variables. Now we turn to compare
\begin{equation*}
    H_{n+1}(\bx,x_{n+1})=\sum_{l=1}^L\sum_{i<j\le n+1}\sqrt{\frac{\lambda^{(l)}}{n+1}}W_{i,j}^{(l)}x_i^{(l)}x_j^{(l)}+\frac{\lambda^{(l)}}{n+1}X_i^{(l)}X_j^{(l)}x_i^{(l)}x_j^{(l)}
\end{equation*}
with $H_n(\bx)=\sum_{l=1}^L\sum_{i<j}\sqrt{\frac{\lambda^{(l)}}{n}}W_{i,j}^{(l)}x_i^{(l)}x_j^{(l)}+\frac{\lambda^{(l)}}{n}X_i^{(l)}X_j^{(l)}x_i^{(l)}x_j^{(l)}$. 
In order to adjust the change of denominator from $n$ to $n+1$ in the coefficients before $\{W_{i,j}^{(l)},i<j\le n\}$, let's introduce $\widetilde{W}_{i,j}^{(l)}\indsim\mathcal{N}(0,1)$ independent of everything else. Noting $\sqrt{\frac{1}{n}}W_{i,j}^{(l)}\overset{d}{=}\sqrt{\frac{1}{n+1}}W_{i,j}^{(l)}+\sqrt{\frac{1}{n(n+1)}}\widetilde{W}_{i,j}^{(l)}$ in the definition of $H_n$, we can decompose $H_{n+1}$ into
\begin{equation}\label{eq:comparing Hamiltonian from n to n+1}
    H_{n+1}(\bx,x_{n+1})=H_n(\bx)-\delta(\bx)+\sum_{l=1}^L s_l(\bx)x_{n+1}^{(l)},\quad s_l(\bx):=\sum_{i=1}^n\sqrt{\frac{\lambdal}{n+1}}A_{i,n+1}^{(l)}x_i^{(l)},
\end{equation}
where $s_l(\bx)$ is used to characterize how $\bx$ influences $x_{n+1}$, while
\begin{equation*}
    \delta(\bx):=\sum_{l=1}^L\sum_{1\le i<j\le n}\frac{\lambdal}{n(n+1)}X_i^{(l)}X_j^{(l)}x_i^{(l)}x_j^{(l)}+\sqrt{\frac{\lambdal}{n(n+1)}}\widetilde{W}_{i,j}^{(l)}x_i^{(l)}x_j^{(l)}
\end{equation*}
is a technical term to adjust changes in denominator from $n$ to $n+1$. From now on, we use $\langle\cdot\rangle$ to denote the expectation over Gibbs distribution defined by Hamiltonian $H_n$. Therefore, \eqref{eq:Cesaro term} becomes,
\begin{equation}\label{eq:Cesaro term 2}
    B_n=\mathbb{E}\log\left\langle\sum_{v_{n+1}}p(v_{n+1}|Z_{n+1})e^{-\delta(\widetilde{\bx})+\sum_{l}s_l(\widetilde{\bx})\widetilde{x}_{n+1}^{(l)}}\right\rangle.
\end{equation}
Now we proceed by incorporating some intuitive mean field approximations, which will be rigorously proved in the next subsection. Expanding their respective definitions,
\begin{align*}
    s_l(\bx)x_{n+1}^{(l)}=&\sum_{i=1}^n\frac{\lambdal}{n+1}x_i^{(l)}X_i^{(l)}X_{n+1}^{(l)}x_{n+1}^{(l)}+\sqrt{\frac{\lambdal}{n+1}}W_{i,n+1}^{(l)}x_i^{(l)}x_{n+1}^{(l)}
\end{align*}
is found to be approximated by $s_l^0x_{n+1}^{(l)}+\frac{\lambdal-\lambdal Q_l}{2}$ where
\begin{equation*}
    s_l^0:=\lambdal Q_lX_{n+1}^{(l)}+\sum_{i=1}^n\sqrt{\frac{\lambdal}{n}}W_{i,n+1}^{(l)}\left\langle x_i^{(l)}\right\rangle.
\end{equation*}
Besides,
\begin{align*}
    \delta(\bx)=\sum_{l=1}^L\frac{\lambdal n}{2(n+1)}\left(\bx^{(l)}.\bX^{(l)}\right)^2-\frac{\lambdal}{2(n+1)}+\sum_{i<j}\sqrt{\frac{\lambdal}{n(n+1)}}\widetilde{W}_{i,j}^{(l)} x_i^{(l)}x_j^{(l)}
\end{align*}
will be approximated by
\begin{equation*}
    \delta^0:=\sum_{l=1}^L \left[\frac{\lambdal Q_l^2+\lambdal}{4}+\sum_{1\le i<j\le n}\frac{\sqrt{\lambdal}}{n}\widetilde{W}_{i,j}^{(l)}\left\langle x_i^{(l)}x_j^{(l)}\right\rangle\right].
\end{equation*}
These intuitions build upon Lemma~\ref{lemma:overlap concentration after perturbation}, as overlaps $\bx^{(l)}.\bX^{(l)}$ concentrate near $Q_l$ under $\langle\cdot\rangle$. Proof of the following lemma is postponed to next subsection.
\begin{Lemma}\label{lemma:mean field approx}
    Similar to the notion in \eqref{eq:Cesaro term}, we write
    \begin{equation*}
        B_n^0:=\mathbb{E}\log\sum_{v_{n+1}}p(v_{n+1}|Z_{n+1})e^{-\delta^0+\sum_{l}s_l^0\widetilde{x}_{n+1}^{(l)}+\frac{\lambdal-\lambdal Q_l}{2}}.
    \end{equation*}
    Then we have $\mathbb{E}_{\epsilon}|B_n-B_n^0|\nconverge0$.
\end{Lemma}
\begin{proof}[Proof of Proposition~\ref{prop:RS prediction on free energy}]
    Using Proposition~\ref{prop:get rid of perturbation} and Cesaro sum,
\begin{equation*}
    \limsup_{n\rightarrow\infty} F_{n}\le\limsup_{n\rightarrow\infty} F_{n,\epsilon_n}\le \limsup_{n\rightarrow\infty} B_n\le \limsup_{n\rightarrow\infty} B_n^0.
\end{equation*}
While for every $n$, since $\sum_{i=1}^n\sqrt{\frac{\lambdal}{n}}W_{i,n+1}^{(l)}\left\langle x_i^{(l)}\right\rangle$ equals in law to $\sqrt{\lambdal Q_l}\mathcal{N}(0,1)$, it holds
\begin{equation*}
    B_n^0\le \sup_{\boldsymbol{q}\ge 0}\mathbb{E}\log\left[\sum_{(x,y)\in\{\pm 1\}^{L+L_1}}p(x,y|Z)\exp\left(\sum_{l=1}^L\sqrt{\lambda^{(l)}q_l}A^{(l)\prime}x^{(l)}\right)\right]+\sum_{l=1}^{L}\frac{\lambda^{(l)}(1-2q_l-q_l^{2})}{4}.
\end{equation*}
Combined with the other bound~\eqref{eq:replica symmetric lower bound} shown in Section~\ref{subsec:Guerra's interpolation}, the proof of Proposition~\ref{prop:RS prediction on free energy} is complete. 
\end{proof}

\subsection{Proof of Lemma~\ref{lemma:mean field approx}}\label{subsec:pf mean field approx}
Define respectively that
\begin{align*}
    U_1=&\left\langle\sum_{v_{n+1}}p(v_{n+1}|Z_{n+1})e^{-\delta(\widetilde{\bx})+\sum_{l}s_l(\widetilde{\bx})\widetilde{x}_{n+1}^{(l)}}\right\rangle,\\
    \quad U_2=&\sum_{v_{n+1}}p(v_{n+1}|Z_{n+1})e^{-\delta^0+\sum_{l}s_l^0\widetilde{x}_{n+1}^{(l)}+\frac{\lambdal-\lambdal Q_l}{2}}.
\end{align*}
As a result, Lemma~\ref{lemma:mean field approx} can be restated as $\mathbb{E}_{\epsilon}|\mathbb{E}\log U_1-\mathbb{E}\log U_2|\nconverge0$.
\begin{Lemma}\label{lemma:convergence in L2}
    It holds that $\lim_{n\rightarrow}\mathbb{E}_{\epsilon}\mathbb{E}|U_1-U_2|^2=0$.
\end{Lemma}
\begin{proof}
Write $\mathbb{E}_{\widetilde{\bW}}$ as taking expectation over $(\widetilde{W}_{i,j}^{(l)})_{1\le i<j\le n}$ only, and $\mathbb{E}_{\bW_{n+1}}$ as taking taking expectation over $(W_{i,n+1}^{(l)})_{1\le i\le n}$ only. We can exchange $\mathbb{E}_{\widetilde{\bW}}$ (or $\mathbb{E}_{\bW_{n+1}}$) with $\langle\cdot\rangle$ because these standard normals are independent to $\langle\cdot\rangle$.
Plug MGF $\mathbb{E}e^{tW}=e^{t^2/2}$ for standard normal $W$ and $Q_l=\frac{1}{n}\left\langle x_i^{(l)}\right\rangle^2$, we first compute
\begin{align*}
    \mathbb{E}U_2^2=&\sum_{v_{n+1},\bar{v}_{n+1}}p(v_{n+1},\bar{v}_{n+1}|Z_{n+1})\mathbb{E}_{\widetilde{\bW},\bW_{n+1}}e^{-2\delta^0+\sum_{l}s_l^0\left(\widetilde{x}_{n+1}^{(l)}+\widetilde{\bar{x}}_{n+1}^{(l)}\right)+\lambdal-\lambdal Q_l}\\
    =&\sum_{v_{n+1},\bar{v}_{n+1}}p(v_{n+1},\bar{v}_{n+1}|Z_{n+1})\prod_l\exp\Bigg\{\frac{\lambdal-\lambdal Q_l^2}{2}-\sum_{i<j\le n}\frac{2\lambdal}{n^2}\left\langle x_i^{(l)}x_j^{(l)}\right\rangle^2\\
    &\quad +\lambdal Q_l\left(X_{n+1}^{(l)}\widetilde{x}_{n+1}^{(l)}+X_{n+1}^{(l)}\widetilde{\bar{x}}_{n+1}^{(l)}+\widetilde{x}_{n+1}^{(l)}\widetilde{\bar{x}}_{n+1}^{(l)}\right)\Bigg\}.
\end{align*}
Moreover, since
\begin{equation*}
    \sum_{i<j\le n}\frac{2\lambdal}{n^2}\left\langle x_i^{(l)}x_j^{(l)}\right\rangle^2=\lambdal\left(\bx^{(l)}.\bX^{(l)}\right)^2+o(1)=\lambdal Q_l^2+o(1),
\end{equation*}
we can further simplify
\begin{equation}
    \mathbb{E}U_2^2=\sum_{v_{n+1},\bar{v}_{n+1}}p(v_{n+1},\bar{v}_{n+1}|Z_{n+1})e^{\sum_l\frac{\lambdal-3\lambdal Q_l^2}{2}+\lambdal Q_l\left(X_{n+1}^{(l)}\widetilde{x}_{n+1}^{(l)}+X_{n+1}^{(l)}\widetilde{\bar{x}}_{n+1}^{(l)}+\widetilde{x}_{n+1}^{(l)}\widetilde{\bar{x}}_{n+1}^{(l)}\right)}+o(1).
\end{equation}
Now we use MGF again to find,
\begin{align*}
    \mathbb{E}_{\widetilde{\bW}}\exp\left[\sqrt{\frac{\lambdal}{n(n+1)}}\widetilde{W}_{i,j}^{(l)}\left(x_i^{(l)}x_j^{(l)}+\bar{x}_i^{(l)}\bar{x}_j^{(l)}\right)\right]=&\exp\left[\frac{\lambdal}{n(n+1)}\left(1+x_i^{(l)}x_j^{(l)}\bar{x}_i^{(l)}\bar{x}_j^{(l)}\right)\right],\\
    \mathbb{E}_{\bW_{n+1}}\exp\left(\sqrt{\frac{\lambdal}{n+1}}W_{i,n+1}^{(l)}\left(x_i^{(l)}x_{n+1}^{(l)}+\bar{x}_i^{(l)}\bar{x}_{n+1}^{(l)}\right)\right)=&\exp\left[\frac{\lambdal}{n+1}\left(1+x_i^{(l)}x_{n+1}^{(l)}\bar{x}_i^{(l)}\bar{x}_{n+1}^{(l)}\right)\right].
\end{align*}
Next, we compute
\begin{align*}
    \mathbb{E}U_1^2=&\left\langle\sum_{v_{n+1},\bar{v}_{n+1}}p(v_{n+1},\bar{v}_{n+1}|Z_{n+1})\mathbb{E}_{\widetilde{\bW},\bW_{n+1}}e^{-\delta(\bx)-\delta(\bar{\bx})+\sum_{l}s_l(\bx)\widetilde{x}_{n+1}^{(l)}+s_l(\bar{\bx})\widetilde{\bar{x}}_{n+1}^{(l)}}\right\rangle\\
    =&\Bigg\langle\sum_{v_{n+1},\bar{v}_{n+1}}p(v_{n+1},\bar{v}_{n+1}|Z_{n+1})\prod_l\exp\Bigg\{\frac{\lambdal}{2}-\frac{\lambdal}{2}\left(\bx^{(l)}.\bX^{(l)}\right)^2-\frac{\lambdal}{2}\left(\bar{\bx}^{(l)}.\bX^{(l)}\right)^2\\
    &\quad -\frac{\lambdal}{2}\left(\bx^{(l)}.\bar{\bx}^{(l)}\right)^2+\lambdal\left(\bx^{(l)}.\bX^{(l)}\right)X_{n+1}^{(l)}\widetilde{x}_{n+1}^{(l)}+\lambdal\left(\bar{\bx}^{(l)}.\bX^{(l)}\right)X_{n+1}^{(l)}\widetilde{\bar{x}}_{n+1}^{(l)}\\
    &\quad +\lambdal\left(\bx^{(l)}.\bar{\bx}^{(l)}\right)x_{n+1}^{(l)}\widetilde{\bar{x}}_{n+1}^{(l)}\Bigg\}\Bigg\rangle+o(1).
\end{align*}
Define for any triple $(r_1,r_2,r_3)$ where $r_j=(r_j^{(1)},\ldots,r_j^{(L)})$ the following mapping,
\begin{align*}
    f(r_1,r_2,r_3)=&\sum_{v_{n+1},\bar{v}_{n+1}}p(v_{n+1},\bar{v}_{n+1}|Z_{n+1})\exp\Bigg(\sum_{l=1}^L \frac{\lambdal}{2}+\lambdal r_1^{(l)}X_{n+1}^{(l)}\widetilde{x}_{n+1}^{(l)}\\
    &\quad +r_2^{(l)}X_{n+1}^{(l)}\widetilde{\bar{x}}_{n+1}^{(l)}+r_3^{(l)}x_{n+1}^{(l)}\widetilde{\bar{x}}_{n+1}^{(l)}-\frac{\lambdal}{2}(r_1^{(l)2}+r_2^{(l)2}+r_3^{(l)2})\Bigg).
\end{align*}
Equipped with this notion, we have 
\begin{align*}
    \mathbb{E}U_2^2=&f(Q,Q,Q)+o(1),\\
    \mathbb{E}U_1^2=&\left\langle f(\bx.\bX, \bar{\bx}.\bX, \bx.\bar{\bx})\right\rangle+o(1).
\end{align*}
For each $j\in[3]$ and $l\in[L]$, it holds that $\left|\frac{\partial f}{\partial r_j^{(l)}}\right|$ is bounded by some positive constant only depending on $\blam$. Therefore, $f$ is Lipschitz in each of its arguments. We conclude $\mathbb{E}_\epsilon \left|\mathbb{E}U_1^2-\mathbb{E}U_2^2\right|\nconverge0$ from Lemma~\ref{lemma:overlap concentration after perturbation}. We can also compute
\begin{align*}
    \mathbb{E}U_1U_2=&\left\langle\sum_{v_{n+1},\bar{v}_{n+1}}p(v_{n+1},\bar{v}_{n+1}|Z_{n+1})\mathbb{E}_{\widetilde{\bW},\bW_{n+1}}e^{-\delta(\bx)-\delta^0+\sum_{l}s_l(\bx)\widetilde{x}_{n+1}^{(l)}+s_l^0\widetilde{\bar{x}}_{n+1}^{(l)}+\frac{\lambdal-\lambdal Q_l}{2}}\right\rangle\\
    =& \left\langle f(\bx.\bX, \bb.\bX, \bx.\bb)\right\rangle+o(1),
\end{align*}
where $b^{(l)}_i=\left\langle x_i^{(l)}\right\rangle$, $\bb^{(l)}.\bx^{(l)}=\frac{1}{n}\sum_{i}b_i^{(l)}x_i^{(l)}$ and $\bb^{(l)}.\bX^{(l)}=\frac{1}{n}\sum_{i}b_i^{(l)}X_i^{(l)}$. Once again, we would conclude $\mathbb{E}_\epsilon \left|\mathbb{E}U_1U_2-\mathbb{E}U_2^2\right|\nconverge0$ from Lemma~\ref{lemma:overlap concentration after perturbation} and Nishimori identity. The desired $L^2$ convergence now follows.
\end{proof}
\begin{Lemma}\label{lemma:bounded inverse square}
    There exists a constant $C$ independent of $n$ such that $\mathbb{E}U_1^{-2}+\mathbb{E}U_2^{-2}\le C$.
\end{Lemma}
\begin{proof}
    Using Jensen inequality,
    \begin{align*}
        U_1^{-2}\le &\left(\sum_{v_{n+1}}p(v_{n+1}|Z_{n+1})e^{-\left\langle\delta(\bx)\right\rangle+\sum_{l}\left\langle s_l(\bx)\right\rangle\widetilde{x}_{n+1}^{(l)}}\right)^{-2}\\
        \le &\sum_{v_{n+1}}p(v_{n+1}|Z_{n+1})e^{2\left\langle\delta(\bx)\right\rangle-2\sum_{l}\left\langle s_l(\bx)\right\rangle\widetilde{x}_{n+1}^{(l)}}.
    \end{align*}
    The random variable $\left\langle\delta(\bx)\right\rangle-\sum_{l}\left\langle s_l(\bx)\right\rangle\widetilde{x}_{n+1}^{(l)}$ has the same law as 
    \begin{equation*}
        \sum_{l}\frac{\lambdal Q_l^2}{2}-\lambdal Q_l X_{n+1}^{(l)}x_{n+1}^{(l)}+\sqrt{\frac{3}{2}\sum_{l}\lambdal Q_l}\mathcal{N}(0,1).
    \end{equation*}
    Thus we have
    \begin{equation*}
        \mathbb{E}U_1^{-2}\le\mathbb{E} \sum_{v_{n+1}}p(v_{n+1}|Z_{n+1})\mathbb{E}e^{2\left\langle\delta(\bx)\right\rangle-2\sum_{l}\left\langle s_l(\bx)\right\rangle\widetilde{x}_{n+1}^{(l)}}\le C.
    \end{equation*}
    The bound on   $\mathbb{E}U_2^{-2}$ follows similarly. 
\end{proof}
Lastly, we use $|\log U_1-\log U_2|\le \max(U_1^{-1},U_2^{-2})|U_1-U_2|$ and Cauchy Schwarz inequality so that 
\begin{equation*}
    \mathbb{E}_{\epsilon}|\mathbb{E}\log U_1-\mathbb{E}\log U_2|\le \sqrt{\mathbb{E}_{\epsilon}\mathbb{E}U_1^{-2}+\mathbb{E}_{\epsilon}\mathbb{E}U_2^{-2}}\sqrt{\mathbb{E}_{\epsilon}\mathbb{E}|U_1-U_2|^2}\nconverge0,
\end{equation*}
where the final convergence is evident from Lemmas~\ref{lemma:convergence in L2} and~\ref{lemma:bounded inverse square}. The proof of  Lemma~\ref{lemma:mean field approx} is now complete.

\section{Proof of Theorem~\ref{thm:asymptotic MI}, Proposition~\ref{prop:Universality MI qualitative} and Proposition~\ref{prop:universality mmse}}
\label{sec:mi_mmse}

\subsection{Asymptotic Normalized MI and MMSE}
\begin{proof}[Proof of Theorem~\ref{thm:asymptotic MI} (i)]
    By definition,
    \begin{align*}
        \frac{1}{n}I(\bX,\bY;\bA|\bZ)=&-\frac{1}{n}\mathbb{E}\log\left[\frac{p(\bA|\bZ)}{p(\bA|\bX,\bY,\bZ)}\right]=-\frac{1}{n}\mathbb{E}\log\left[\frac{\sum_{\bx,\by}p(\bx,\by|\bZ)p(\bA|\bx)}{p(\bA|\bX)}\right]\\
        =&-\frac{1}{n}\mathbb{E}\log\left[\frac{\sum_{\bx,\by}p(\bx,\by|\bZ)\exp\left(\sum_{l}\sum_{i<j}\sqrt{\lambdal/n}A_{i,j}^{(l)}x_i^{(l)}x_j^{(l)}\right)}{\exp\left(\sum_{l}\sum_{i<j}\sqrt{\lambdal/n}A_{i,j}^{(l)}X_i^{(l)}X_j^{(l)}\right)}\right]\\
        =&\frac{1}{n}\mathbb{E}\left[\sum_{l}\sum_{i<j}\sqrt{\lambdal/n}A_{i,j}^{(l)}X_i^{(l)}X_j^{(l)}\right]-F_n(\boldsymbol{\lambda})\\
        =&\sum_{l}\frac{(n-1)\lambdal}{2n}-F_n(\boldsymbol{\lambda}).
    \end{align*}
    Combined with Proposition~\ref{prop:RS prediction on free energy}, we can establish
    \begin{equation}\label{eq:limiting conditional MI}
       \lim_{n \to \infty} \frac{1}{n}I(\bX,\bY;\bA|\bZ) =  \frac{\sum_l\lambdal}{4}-\sup_{\mathbf{q}\ge 0}\left[\mathcal{F}(\boldsymbol{\lambda},\mathbf{q})-\sum_{l=1}^{L}\frac{\lambda^{(l)}(q_l^2+2q_l)}{4}\right].
    \end{equation}
    Subsequently, since $(\bX,\bY,\bZ)$ comprise $n$ independent repetitions of the prior $p(x,y,z)$,
    \begin{align*}
        \frac{1}{n}I(\bX,\bY;\bA,\bZ)&=\frac{1}{n}I(\bX,\bY;\bZ)+\frac{1}{n}I(\bX,\bY;\bA|\bZ)=i_p(x,y;z)+\frac{1}{n}I(\bX,\bY;\bA|\bZ)\\
        &\nconverge i_p(x,y;z)+\frac{\sum_l\lambdal}{4}-\sup_{\mathbf{q}\ge 0}\left[\mathcal{F}(\boldsymbol{\lambda},\mathbf{q})-\sum_{l=1}^{L}\frac{\lambda^{(l)}(q_l^2+2q_l)}{4}\right],
    \end{align*}
    confirming the first assertion \eqref{eq:limiting MI} of this theorem, where $i_p(x,y;z)$ is the mutual information between $(x,y)$ and $z$ under prior $p(\cdot)$.
\end{proof}
\begin{proof}[Proof of Theorem~\ref{thm:asymptotic MI} (ii)]
    Recall the definition of the free energy from \eqref{eq:free energy def}.
    By direct differentiation, we obtain
    \begin{align*}
        \frac{\partial}{\partial \lambda^{(l)}}F_n(\boldsymbol{\lambda})=&\frac{1}{n}\mathbb{E}\frac{\sum_{\bx,\by}\prod_i p(x_i,y_i|Z_i)\cdot\exp\left(H(\bx)\right)\frac{\partial}{\partial \lambdal}H(\bx)}{\sum_{\bx,\by}\prod_i p(x_i,y_i|Z_i)\cdot\exp\left(H(\bx)\right)}\\
        =&\frac{1}{n}\mathbb{E}\frac{\sum_{\bx,\by}\prod_i p(x_i,y_i|Z_i)\cdot\exp\left(H(\bx)\right)\left[\sum_{i<j}\frac{1}{n}X_i^{(l)}X_j^{(l)}x_i^{(l)}x_j^{(l)}+\frac{1}{2\sqrt{n\lambdal}}W_{i,j}^{(l)}x_i^{(l)}x_j^{(l)}\right]}{\sum_{\bx,\by}\prod_i p(x_i,y_i|Z_i)\cdot\exp\left(H(\bx)\right)}\\
        =&\frac{1}{n}\mathbb{E}\frac{\sum_{\bx,\by}\prod_i p(x_i,y_i|Z_i)\cdot\exp\left(H(\bx)\right)\left[\sum_{i<j}\frac{1}{2n}X_i^{(l)}X_j^{(l)}x_i^{(l)}x_j^{(l)}+\frac{1}{2n}\right]}{\sum_{\bx,\by}\prod_i p(x_i,y_i|Z_i)\cdot\exp\left(H(\bx)\right)}\\
        =&\frac{1}{2n^2}\mathbb{E}\left[\sum_{i<j}\left\langle X_i^{(l)}X_j^{(l)}x_i^{(l)}x_j^{(l)}\right\rangle\right]+\frac{n-1}{4n},
    \end{align*}
    where $\langle\cdot\rangle$ refers to averages under the posterior. Note that we use Gaussian integration by parts in the third step above. Second-order derivatives can be calculated as 
    \begin{align*}
        &\frac{\partial^2}{\partial \lambda^{(l_1)}\partial\lambda^{(l_2)}}F_n(\boldsymbol{\lambda})
        =\frac{1}{2n^2}\mathbb{E}\left[\frac{\partial}{\partial \lambda^{(l_2)}}\sum_{i<j}\left\langle X_i^{(l_1)}X_j^{(l_1)}x_i^{(l_1)}x_j^{(l_1)}\right\rangle\right]\\
        =&\frac{1}{2n^2}\mathbb{E}\left\{\sum_{i_2<j_2}\left[\left\langle\left(\sum_{i<j} X_i^{(l_1)}X_j^{(l_1)}x_i^{(l_1)}x_j^{(l_1)}\right)x^{(l_2)}_{i_2}x^{(l_2)}_{j_2}\right\rangle\right.\right.\\
        &\qquad\left.\left.-\left\langle\sum_{i<j} X_i^{(l_1)}X_j^{(l_1)}x_i^{(l_1)}x_j^{(l_1)}\right\rangle\left\langle x^{(l_2)}_{i_2}x^{(l_2)}_{j_2}\right\rangle\right]\cdot\left(\frac{1}{n}X^{(l_2)}_{i_2}X^{(l_2)}_{j_2}+\frac{1}{2\sqrt{n\lambda^{(l_2)}}}W_{i_2,j_2}^{(l_2)}\right)\right\}\\
        =&\frac{1}{2n^3}\mathbb{E}\left\{\sum_{i_2<j_2}\left[\left\langle\left(\sum_{i<j} X_i^{(l_1)}X_j^{(l_1)}x_i^{(l_1)}x_j^{(l_1)}\right)x^{(l_2)}_{i_2}x^{(l_2)}_{j_2}X^{(l_2)}_{i_2}X^{(l_2)}_{j_2}\right\rangle\right.\right.\\
        &\qquad\qquad+\left\langle\sum_{i<j} X_i^{(l_1)}X_j^{(l_1)}x_i^{(l_1)}x_j^{(l_1)}\right\rangle\left\langle x^{(l_2)}_{i_2}x^{(l_2)}_{j_2}\right\rangle^2\\
        &\qquad\qquad\left.\left.-2\left\langle\sum_{i<j} X_i^{(l_1)}X_j^{(l_1)}x_i^{(l_1)}x_j^{(l_1)}\right\rangle\left\langle x^{(l_2)}_{i_2}x^{(l_2)}_{j_2}X^{(l_2)}_{i_2}X^{(l_2)}_{j_2}\right\rangle\right]\right\}\\
        =&\frac{1}{2n^3}\mathbb{E}\left[\sum_{i_1<j_1}\sum_{i_2<j_2}\left(\left\langle x_{i_1}^{(l_1)}x_{j_1}^{(l_1)}x_{i_2}^{(l_2)}x_{j_2}^{(l_2)}\right\rangle-\left\langle x_{i_1}^{(l_1)}x_{j_1}^{(l_1)}\right\rangle\left\langle x_{i_2}^{(l_2)}x_{j_2}^{(l_2)}\right\rangle\right)^2\right],
    \end{align*}
    where the third display follows from Gaussian integration by parts and the last equation follows from the Nishimori identity. To simplify the notations, we write
    \begin{equation*}
        \xi_{i,j}^{(l)}=x_{i}^{(l)}x_{j}^{(l)}-\left\langle x_{i}^{(l)}x_{j}^{(l)}\right\rangle,\quad i<j,\quad l\in[L],
    \end{equation*}
    as a centered functional of any replica $\mathbf{x}$ drawn from $\langle\cdot\rangle$. Moreover, denote
    \begin{equation*}
        \bar{\xi}_{i,j}^{(l)}=\bar{x}_{i}^{(l)}\bar{x}_{j}^{(l)}-\left\langle \bar{x}_{i}^{(l)}\bar{x}_{j}^{(l)}\right\rangle,
    \end{equation*}
    when this functional is applied to another independent replica $\bar{\mathbf{x}}$. With this notion, we can simplify $\frac{\partial^2}{\partial \lambda^{(l_1)}\partial\lambda^{(l_2)}}F_n(\boldsymbol{\lambda})$ to
    \begin{equation*}
        \frac{\partial^2}{\partial \lambda^{(l_1)}\partial\lambda^{(l_2)}}F_n(\boldsymbol{\lambda})=\frac{1}{2n^3}\mathbb{E}\left[\sum_{i_1<j_1}\sum_{i_2<j_2}\left\langle\xi_{i_1,j_1}^{(l_1)}\xi_{i_2,j_2}^{(l_2)}\bar{\xi}_{i_1,j_1}^{(l_1)}\bar{\xi}_{i_2,j_2}^{(l_2)}\right\rangle\right].
    \end{equation*}
    For any $(m_1,\ldots,m_L)\in\mathbb{R}^L$,
    \begin{align*}
        \sum_{l_1,l_2}m_{l_1}m_{l_2}\frac{\partial^2}{\partial \lambda^{(l_1)}\partial\lambda^{(l_2)}}F_n(\boldsymbol{\lambda})=\frac{1}{2n^3}\mathbb{E}\left[\left\langle\left( \sum_{l=1}^L\sum_{i<j}m_{l}\xi_{i,j}^{(l)}\bar{\xi}_{i,j}^{(l)}\right)^2\right\rangle\right]\ge 0.
    \end{align*}
    In conclusion, $F_n(\boldsymbol{\lambda})$ is jointly convex in $\boldsymbol{\lambda}$.
    Using Proposition~\ref{prop:RS prediction on free energy}, this series of multivariate convex functions converge pointwise to function
    \begin{equation}\label{eq:limiting functional of free energy}
        \boldsymbol{\lambda}\mapsto \frac{\sum_l\lambdal}{4}+\sup_{\mathbf{q}\ge 0}\left[\mathcal{F}(\boldsymbol{\lambda},\mathbf{q})-\sum_{l=1}^{L}\frac{\lambda^{(l)}(q_l^2+2q_l)}{4}\right],
    \end{equation}
    as $n\rightarrow\infty$. Thus the limiting function is also convex. 
    Set 
    \begin{equation*}
        D=\left\{\boldsymbol{\lambda}\in[0,\infty)^L:\sup_{\mathbf{q}\ge 0}\left[\mathcal{F}(\boldsymbol{\lambda},\mathbf{q})-\sum_{l=1}^{L}\frac{\lambda^{(l)}(q_l^2+2q_l)}{4}\right] \text{ is not differentiable at } \boldsymbol{\lambda}\right\}.
    \end{equation*}
  Theorem 6.7(i) in \cite{evans2018measure} implies that every multivariate convex function is locally Lipschitz on an open subset of $\mathbb{R}^L$. Using Rademacher’s Theorem (Theorem 3.1 or Theorem 6.6 in \cite{evans2018measure}), we have that any locally Lipschitz function is differentiable almost everywhere. 
  Thus the set $D\subset[0,\infty)^L$ of bad points is of Lebesgue measure zero.
    Define $\boldsymbol{\gamma}=\left(\lambda^{(1)}q_1,\ldots,\lambda^{(L)}q_L\right)$, then we obtain
    \begin{align}
        &\sup_{\mathbf{q}\ge 0}\left[\mathcal{F}(\boldsymbol{\lambda},\mathbf{q})-\sum_{l=1}^{L}\frac{\lambda^{(l)}(q_l^2+2q_l)}{4}\right]=\sup_{\boldsymbol{\gamma}\ge 0}\left[-\sum_{l=1}^L\frac{\gamma_l^2}{4\lambdal}+\bar{\mathcal{F}}(\boldsymbol{\gamma})\right],\label{eq:sup reparamtried free energy functional}\\
        \bar{\mathcal{F}}(\boldsymbol{\gamma}):=\mathbb{E}&\log \left(\sum_{x,y}p(x,y|Z)e^{\sum_{l=1}^L\gamma_l X^{(l)}x^{(l)}+\sqrt{\gamma_l}W^{(l)\prime}x^{(l)}}\right)-\sum_{l=1}^{L}\frac{\gamma_l}{2}=\mathcal{F}(\boldsymbol{\lambda},\mathbf{q})-\sum_{l=1}^{L}\frac{\lambdal q_l}{2}.\notag
    \end{align}
    It is easy to see that $\partial_{\gamma_l}\bar{\mathcal{F}}(\boldsymbol{\gamma})\le 1/2< \gamma_l/(2\lambdal)$ if $\gamma_l>\lambdal$, so we can restrict the domain $\boldsymbol{\gamma}\ge 0$ of RHS in \eqref{eq:sup reparamtried free energy functional} onto a compact subset $\boldsymbol{\gamma}\in[0,\lambda^{(1)}]\times\cdots\times[0,\lambda^{(L)}]$, thus enabling the use of envelope theorems, \cite{milgrom2002envelope}. For each $l\in[L]$, apply the last argument \cite[Corollary 4]{milgrom2002envelope} on the functional
    \begin{equation*}
        G:(\lambdal,\boldsymbol{\gamma})\mapsto -\sum_{l=1}^L\frac{\gamma_l^2}{4\lambdal}+\bar{\mathcal{F}}(\boldsymbol{\gamma}),
    \end{equation*}
    to conclude: for $\boldsymbol{\lambda}$ such that $\sup_{\boldsymbol{\gamma}}G(\lambdal,\boldsymbol{\gamma})$ is differentiable in $\lambdal$, there exists a unique maximizer $\boldsymbol{\gamma}^\ast$ to $\sup_{\boldsymbol{\gamma}}G(\lambdal,\boldsymbol{\gamma})$ and 
    \begin{equation*}
        \frac{\partial}{\partial \lambda^{(l)}}\left[\sup_{\boldsymbol{\gamma}}G(\lambdal,\boldsymbol{\gamma})\right]=\frac{\partial G}{\partial \lambda^{(l)}}(\lambdal,\boldsymbol{\gamma}^\ast)=\frac{\left(\gamma^\ast_{l}\right)^2}{4\lambda^{(l)2}}.
    \end{equation*}
    Equivalently, we conclude that for any $\boldsymbol{\lambda}\in D^c$, there exists a unique solution $\mathbf{q}^\ast$ to the maximization problem in functional \eqref{eq:limiting functional of free energy} and its derivative is given by
    \begin{equation}
        \frac{\partial}{\partial \lambda^{(l)}}\left[\sup_{\mathbf{q}\ge 0}\left[\mathcal{F}(\boldsymbol{\lambda},\mathbf{q})-\sum_{l=1}^{L}\frac{\lambda^{(l)}(q_l^2+2q_l)}{4}\right]\right]=\frac{\left(q^{\ast}_{l}\right)^2}{4}.\label{eq:envelope}
    \end{equation}

    The next display  provides a direct relation connecting the minimal mean squared error to first-order derivatives $\frac{\partial}{\partial \lambda^{(l)}}F_n(\boldsymbol{\lambda})$,
    \begin{align*}
        \MMSE\left(\bX^{(l)};\bA,\bZ\right)=&\frac{2}{n(n-1)}\sum_{i<j}\mathbb{E}\left(X_i^{(l)}X_j^{(l)}-\mathbb{E}\left[X_i^{(l)}X_j^{(l)}|\bA,\bZ\right]\right)^2\\
        =&\frac{2}{n(n-1)}\sum_{i<j}\left[1-\mathbb{E}\left\langle X_i^{(l)}X_j^{(l)}x_i^{(l)}x_j^{(l)}\right\rangle\right]\\
        =&1- \frac{2}{n(n-1)}\mathbb{E}\left[\sum_{i<j}\left\langle X_i^{(l)}X_j^{(l)}x_i^{(l)}x_j^{(l)}\right\rangle\right]\\
        =&2-\frac{4n}{n-1}\frac{\partial}{\partial \lambda^{(l)}}F_n(\boldsymbol{\lambda}).
    \end{align*}
    If $\boldsymbol{\lambda}\notin D$, the univariate function $\lambdal\mapsto \lim_{n\rightarrow\infty}F_n(\boldsymbol{\lambda})$ is differentiable at 
    $\lambdal$. 
    By convexity in $\lambdal$ and pointwise convergence of $F_n(\boldsymbol{\lambda})$, one can derive the convergence of first-order derivatives,
    \begin{equation*}
        \frac{\partial}{\partial \lambda^{(l)}} F_n(\boldsymbol{\lambda})\nconverge \frac{\partial}{\partial \lambda^{(l)}} \lim_{n\rightarrow\infty}F_n(\boldsymbol{\lambda})= \frac{\partial}{\partial \lambda^{(l)}}\left\{\frac{\sum_l\lambdal}{4}+\sup_{\mathbf{q}\ge 0}\left[\mathcal{F}(\boldsymbol{\lambda},\mathbf{q})-\sum_{l=1}^{L}\frac{\lambda^{(l)}(q_l^2+2q_l)}{4}\right]\right\}.
    \end{equation*}
    Using \eqref{eq:envelope} we conclude
    \begin{equation*}
        \MMSE\left(\bX^{(l)};\bA,\bZ\right)\nconverge 1-\left(q^{\ast}_{l}\right)^2,\quad\forall l\in[L]
    \end{equation*}
    for all $\boldsymbol{\lambda}\in D^c$.
\end{proof}
\begin{proof}[Proof of Theorem~\ref{thm:asymptotic MI} (iii)]
    We turn to optimal estimation of the implicit memberships $\bY^{(l_1)}$, $l_1\in[L_1]$ in this part. Fix $\lambda^y\ge0$ and assume that we observe an additional spiked matrix of the form
    \begin{equation}\label{eq:augmented spike mat}
        \mathbf{A}^{y}=\sqrt{\frac{\lambda^{y}}{n}}\bY^{(l_1)}\bY^{(l_1)\top}+\mathbf{W}^{y},
    \end{equation}
    where $\mathbf{W}^y$ is another independent GOE. Our framework can incorporate this additional observation, and thus for any $(\lambda^y,\boldsymbol{\lambda})\in[0,\infty)^{L+1}$, \eqref{eq:limiting conditional MI} implies the conditional mutual information converges to
    \begin{align*}
        \frac{1}{n}I(\bX,\bY;\bA,\bA^y|\bZ)\nconverge \mathcal{I}(\lambda^y,\boldsymbol{\lambda}),
    \end{align*}
    with limiting functional given by
    \begin{equation*}
        \mathcal{I}(\lambda^y,\boldsymbol{\lambda})=\frac{\lambda^y+\sum_l\lambdal}{4}-\sup_{q_y,\mathbf{q}\ge 0}\left[\mathcal{F}(\lambda^y,\boldsymbol{\lambda},q^y,\mathbf{q})-\frac{\lambda^y(q_y^2+2q_y)}{4}-\sum_{l=1}^{L}\frac{\lambda^{(l)}(q_l^2+2q_l)}{4}\right],
    \end{equation*}
    where the augmented free energy functional is defined as
    \begin{equation}
        \mathcal{F}(\lambda^y,\boldsymbol{\lambda},q^y,\mathbf{q})=\mathbb{E}\log \left(\sum_{x,y}p(x,y|Z)e^{\lambda^{y}q_yY^{(l_1)}y^{(l_1)}+\sqrt{\lambda^{y}q_y}W^{\prime}y^{(l_1)}+\sum_{l=1}^L\lambda^{(l)}q_lX^{(l)}x^{(l)}+\sqrt{\lambda^{(l)}q_l}W^{(l)\prime}x^{(l)}}\right).\label{eq:free energy augmented scalar channel}
    \end{equation}
    We fix $\boldsymbol{\lambda}\in(0,\infty)^L\backslash D$ and vary $\lambda^y$. Since $\mathcal{I}(\lambda^y,\boldsymbol{\lambda})$ is concave in $\lambda^y\in[0,\infty)$, there exists be a countable set $D^y\subset\R$ such that $\mathcal{I}(\lambda^y,\boldsymbol{\lambda})$ is continuously differentiable in $\lambda^y\in(0,\infty)\backslash D^y$. Using the  envelope theorem \cite[Corollary 4]{milgrom2002envelope} for $\lambda^y\in(0,\infty)\backslash D^y$, there exists a unique maximizer $(q_y^\ast(\lambda^y),\mathbf{q}^\ast(\lambda^y))\in\mathbb{R}^{L+1}$ and it satisfies
    \begin{equation*}
        \MMSE\left(\bY^{(l_1)};\bA,\bA^y,\bZ\right)\nconverge 1-q^{\ast}_y\left(\lambda^y\right)^2.
    \end{equation*}
    \begin{Lemma}
    Suppose $\boldsymbol{\lambda}\in[0,\infty)^L\backslash D$ is fixed and denote the associated unique maximizer $\mathbf{q}^\ast\in [0,\infty)^L$ (in the model without the augmented imaginary observation $\bA^y$ in \eqref{eq:augmented spike mat}). As $\lambda^y\in(0,\infty)\backslash D^y$ and $\lambda^y\rightarrow 0$, we have
        \begin{equation*}
            q_y^\ast(\lambda^y)\rightarrow\mathbb{E}\left[Y^{(l_1)}\mathbb{E}\left(Y^{(l_1)}|\sqrt{\boldsymbol{\lambda}\odot\mathbf{q}^\ast}\odot X+W^\prime\right)\right].
        \end{equation*}
    \end{Lemma}
    \begin{proof}
        For $\lambda^y\in(0,\infty)\backslash D^y$ and $\lambda^y\rightarrow 0$, maximizers $(q_y^\ast(\lambda^y),\mathbf{q}^\ast(\lambda^y))$ are monotonic in $\lambda^y$ and bounded, therefore having a finite limit. The last assertion of Proposition~\ref{prop:opt derivative} implies $(q_y^\ast(\lambda^y),\mathbf{q}^\ast(\lambda^y))$ satisfies a set of fixed point conditions,
        \begin{align}
        q_y^\ast(\lambda^y)&=\mathbb{E}\left[Y^{(l_1)}\mathbb{E}\left(Y^{(l_1)}\big|\sqrt{\lambda^y q_y^\ast(\lambda^y)}Y^{(l_1)}+W_y^\prime,\quad \sqrt{\boldsymbol{\lambda}\odot\mathbf{q}^\ast(\lambda^y)}\odot X+W^\prime\right)\right],\label{eq:fixed point eq augmented}\\
        q_l^\ast(\lambda^y)&=\mathbb{E}\left[X^{(l)}\mathbb{E}\left(X^{(l)}\big|\sqrt{\lambda^y q_y^\ast(\lambda^y)}Y^{(l_1)}+W_y^\prime,\quad \sqrt{\boldsymbol{\lambda}\odot\mathbf{q}^\ast(\lambda^y)}\odot X+W^\prime\right)\right],\quad l\in[L].\notag
        \end{align}
        The functional
        \begin{equation*}
            (\lambda^y q_y,\lambda^{(1)}q_1,\ldots,\lambda^{(L)}q_L)\mapsto\mathbb{E}\left[X^{(l)}\mathbb{E}\left(X^{(l)}\big|\sqrt{\lambda^y q_y^\ast(\lambda^y)}Y^{(l_1)}+W_y^\prime,\quad \sqrt{\boldsymbol{\lambda}\odot\mathbf{q}^\ast(\lambda^y)}\odot X+W^\prime\right)\right]
        \end{equation*}
        is smooth, so that the finite limit $\lim_{\lambda^y\rightarrow0}\mathbf{q}^\ast(\lambda^y)$ solves
        \begin{equation*}
            q_l=\mathbb{E}\left[X^{(l)}\mathbb{E}\left(X^{(l)}\big| \sqrt{\boldsymbol{\lambda}\odot\mathbf{q}}\odot X+W^\prime\right)\right],\quad l\in[L].
        \end{equation*}
        Recall that we set $\boldsymbol{\lambda}\in(0,\infty)^L\backslash D$, so this equation set only has at most two solutions. Using continuity, we show that $\lim_{\lambda^y\rightarrow0}\mathbf{q}^\ast(\lambda^y)$ is a maximizer of the free energy functional. Therefore, we must have $\lim_{\lambda^y\rightarrow0}\mathbf{q}^\ast(\lambda^y)=\mathbf{q}^\ast$. Lastly, take $\lambda^y\rightarrow 0$ in \eqref{eq:fixed point eq augmented} to conclude the desired result.
    \end{proof}
    Equipped with this lemma, the lower bound of \eqref{eq:limiting mmse implicit} is immediate,
    \begin{align*}
    &\quad\liminf_{n\rightarrow\infty}\MMSE\left(\bY^{(l_1)};\bA,\bZ\right)\ge\lim_{\lambda^y\notin D^y,\lambda^y\rightarrow0^+}\lim_{n\rightarrow\infty}\MMSE\left(\bY^{(l_1)};\bA,\bA^y,\bZ\right)\\
    &=\lim_{\lambda^y\notin D^y,\lambda^y\rightarrow0^+}\left[ 1-q^{\ast}_y\left(\lambda^y\right)^2\right]
    =1-\mathbb{E}\left[Y^{(l_1)}\mathbb{E}\left(Y^{(l_1)}|\sqrt{\boldsymbol{\lambda}\odot\mathbf{q}^\ast}\odot X+W^\prime\right)\right]^2.
    \end{align*}
    Next, we turn to a matching upper bound. 

    \vspace{2mm}
    \noindent\underline{Step 1.} The first step is to interpolate the additional side information $\bA^y$, similar to Section~\ref{subsec:Guerra's interpolation}, finally deriving an upper bound on $\MMSE\left(\bY^{(l_1)}\bY^{(l_1)\top};\bA,\bZ\right)$. For some $\lambda^y,q_y\ge 0$ and $t\in[0,1]$ to be determined, let's focus on an observation model as this,
    \begin{equation}\label{eq:interpolating additional implicit mat}
    \left\{\begin{aligned}
        &(X_i,Y_i,Z_i)\indsim p(\cdot), &1\le i\le n;\\
        &A_{i,j}^{(l)}=\sqrt{\frac{\lambda^{(l)}}{n}}X_i^{(l)}X_j^{(l)}+W_{i,j}^{(l)}, &1\leq i<j\leq n,1\le l\le L;\\
        &A_{i,j}^{y}=\sqrt{\frac{\lambda^{y}t}{n}}Y_i^{(l_1)}Y_j^{(l_1)\top}+W_{i,j}^{y}, &1\leq i<j\leq n;\\
        &Z^{\prime}_i=\sqrt{(1-t)\lambda^y q_y}Y_i+W^{\prime}_i,&1\leq i\leq n,1\le l\le L;
    \end{aligned}\right.
    \end{equation}
    with every $W$ being an independent standard normal. When $t=1$, \eqref{eq:interpolating additional implicit mat} corresponds to only adding additional observed matrix $\bA^y$; when $t=0$, it corresponds to only adding per-node side information $\bZ^\prime$. Abbreviate normalized mutual information by
    \begin{align*}
        &\quad I_n(\lambda^y,q_y,t,\boldsymbol{\lambda})\\
        &=\frac{1}{n}I(\bX,\bY;\bA,\bA^y,\bZ^\prime|\bZ)\\
        &=\sum_{l}\frac{(n-1)\lambdal}{2n}+\frac{(n-1)\lambda^yt}{2n}+(1-t)\lambda^yq_y-\frac{1}{n}\mathbb{E}\log\left[\sum_{\bx,\by}p(\bx,\by|\bZ)H_n(\bx,\by,t)\right],
    \end{align*}
    where the Hamiltonian $H_n(\bx,\by,t)$ is given by
    \begin{equation*}
        H_n(\bx,\by,t)=\sum_{l}\sum_{i<j}\sqrt{\frac{\lambdal}{n}}A_{i,j}^{(l)}x_i^{(l)}x_j^{(l)}+\sum_{i<j}\sqrt{\frac{\lambda^yt}{n}}A^y_{i,j}y_iy_j+\sum_{i}\sqrt{(1-t)\lambda^yq_y}y_iZ_i^\prime.
    \end{equation*}
    Guerra's interpolation technique, introduced in in Section~\ref{subsec:Guerra's interpolation}, upper bounds the derivative of the mutual information with respect to $t$ by
    \begin{equation*}
        \frac{\partial}{\partial t}I_n(\lambda^y,q_y,t,\boldsymbol{\lambda})\le \frac{\lambda^y}{4}-\frac{\lambda^yq_y}{2}+\frac{\lambda^y q_y^2}{4}.
    \end{equation*}
    Compared to the derivation in Section~\ref{subsec:Guerra's interpolation}, a $o(1)$ term is made more explicit here as we want to derive a non-asymptotic upper bound. Therefore, the following holds for any $\lambda^y,q_y\ge 0$,
    \begin{align*}
        I_n(\lambda^y,q_y,1,\boldsymbol{\lambda})&= I_n(\lambda^y,q_y,0,\boldsymbol{\lambda})+\int_0^1\frac{\partial}{\partial t}I_n(\lambda^y,q_y,t,\boldsymbol{\lambda})\mathrm{d}t\\
        &\le I_n(\lambda^y,q_y,0,\boldsymbol{\lambda})+\frac{\lambda^y}{4}-\frac{\lambda^yq_y}{2}+\frac{\lambda^y q_y^2}{4}.
    \end{align*}
    As RHS equals LHS at $\lambda^y=0$ and they are both continuously differentiable at this boundary for any finite $n$, it must hold for any $q_y\ge0$ that
    \begin{align*}
        \frac{\partial}{\partial \lambda^y}I_n(\lambda^y,q_y,1,\boldsymbol{\lambda}) \Big|_{\lambda^y=0+}&\le \frac{\partial}{\partial \lambda^y}I_n(\lambda^y,q_y,0,\boldsymbol{\lambda})\Big|_{\lambda^y=0+}+\frac{1}{4}(1-q_y)^2.
    \end{align*}
    When $t=1$, differentiating with respect to $\lambda^y$ again gives the minimal mean squared error of $\bY^{(l_1)}\bY^{(l_1)\top}$ for any $\lambda^y\ge0$,
    \begin{align*}
        \frac{\partial}{\partial \lambda^y}I_n(\lambda^y,q_y,1,\boldsymbol{\lambda})&=\frac{1}{2(n-1)^2}\sum_{i<j}\left[1-\mathbb{E}\left\langle Y_i^{(l_1)}Y_j^{(l_1)}y_i^{(l_1)}y_j^{(l_1)}\right\rangle\right]\\
        &=\frac{n}{4(n-1)}\MMSE\left(\bY^{(l_1)};\bA,\bA^y,\bZ\right).
    \end{align*}
    When $t=0$, differentiating with respect to $\lambda^y$ instead gives the minimal mean squared error of $\bY^{(l_1)}$ for any $\lambda^y\ge0$,
    \begin{align*}
        \frac{\partial}{\partial \lambda^y}I_n(\lambda^y,q_y,0,\boldsymbol{\lambda})&=\frac{q_y}{2n}\sum_{i=1}^n\left[1-\mathbb{E}\left\langle Y_i^{(l_1)}y_i^{(l_1)}\right\rangle\right].
    \end{align*}
    Therefore, take $\lambda^y=0$ to derive
    \begin{align}
        \frac{1}{4}\MMSE\left(\bY^{(l_1)}\bY^{(l_1)\top};\bA,\bZ\right)&\le \frac{q_y}{2}\MMSE\left(\bY^{(l_1)};\bA,\bZ\right)+\frac{1}{4}(1-q_y)^2\notag\\
        &\le \frac{1}{4}+\frac{q_y^2}{4}-\frac{q_y}{2}\mathbb{E}\left[Y_{n}^{(l_1)}\langle y_n^{(l_1)}\rangle\right].\label{eq:upper bound on matrix MMSE}
    \end{align}

    \vspace{2mm}
    \noindent\underline{Step 2.}
    Our second step is to conduct a rigorous cavity computation as in Sections~\ref{subsec:Aizenman-Sims-Starr} and \ref{subsec:pf mean field approx}.
Increasing $n$ to $n+1$, we can decompose the Hamiltonian in the same way as Section~\ref{subsec:Aizenman-Sims-Starr},
\begin{align*}
    H_{n+1}(\bx,x_{n+1})&=H_n(\bx)-\delta(\bx)+\sum_{l=1}^L s_l(\bx)x_{n+1}^{(l)},\\
    s_l(\bx):&=\sum_{i=1}^n\frac{\lambdal}{n+1}x_i^{(l)}X_i^{(l)}X_{n+1}^{(l)}x_{n+1}^{(l)}+\sqrt{\frac{\lambdal}{n+1}}W_{i,n+1}^{(l)}x_i^{(l)}x_{n+1}^{(l)},
\end{align*}
where $\bx$ comprises the first $n$ elements and $s_l(\bx)$ is used to characterize how $\bx$ influences $x_{n+1}$, while
\begin{equation*}
    \delta(\bx)=\sum_{l=1}^L\frac{\lambdal n}{2(n+1)}\left(\bx^{(l)}.\bX^{(l)}\right)^2-\frac{\lambdal}{2(n+1)}+\sum_{i<j}\sqrt{\frac{\lambdal}{n(n+1)}}\widetilde{W}_{i,j}^{(l)} x_i^{(l)}x_j^{(l)}
\end{equation*}
is a term which adjusts changes in the denominator from $n$ to $n+1$. We emphasize that $\langle\cdot\rangle_{n+1}$ (or $\langle\cdot\rangle_{n}$ resp.) is used to denote expectation over Gibbs distribution defined by Hamiltonian $H_{n+1}$ (or $H_n$ resp.). Continued from \eqref{eq:upper bound on matrix MMSE} in Step 1, we are interested at the following quantity
\begin{equation*}
    \left\langle y^{(l_1)}_{n+1}\right\rangle_{n+1}=\frac{\left\langle \sum_{x_{n+1},y_{n+1}}y^{(l_1)}_{n+1}p(x_{n+1},y_{n+1}|Z_{n+1})\exp\left(-\delta(\bx)+\sum_{l}s_l(\bx)x_{n+1}^{(l)}\right)\right\rangle_n}{\left\langle \sum_{x_{n+1},y_{n+1}}p(x_{n+1},y_{n+1}|Z_{n+1})\exp\left(-\delta(\bx)+\sum_{l}s_l(\bx)x_{n+1}^{(l)}\right)\right\rangle_n}:=\frac{V_1}{U_1}.
\end{equation*}

\begin{Lemma}\label{lemma:overlap concentration without perturb}
For any $\blam\in(0,\infty)^L\backslash D$, suppose $\bx,\bar{\bx}$ are drawn independently from $\langle\cdot\rangle_n$, then the following overlap concentration holds,
    \begin{equation}
    \mathbb{E}\left\langle\left[\left( \bx^{(l)}.\bar{\bx}^{(l)}\right)^2-\left(q_l^\ast\right)^2\right]^2\right\rangle_n\nconverge0,\quad l\in[L].
\end{equation}
\end{Lemma}
\begin{proof}
In the process of proving the RS prediction on normalized free energy (Proposition~\ref{prop:RS prediction on free energy}), Section~\ref{subsec:overlap concentration with perturb} shows the overlap to concentrate under a small perturbation defined in \eqref{eq:perturbation def}. In this case we start from  Proposition~\ref{prop:RS prediction on free energy} and thus we do not need to introduce the perturbation. The proof is essentially the same as Section 4 in~\cite{lelarge2019fundamental}, which is an adaptation of the proof for Ghirlanda-Guerra identities from Section 3.7 in~\cite{panchenko2013sherrington}.
\end{proof}

Setting $q_{l,n}=\langle\frac{1}{n}\sum_{i=1}^n x_i^{(l)}\bar{x}_i^{(l)}\rangle_n$ to denote the expected overlap, it follows that $\E\left|q_{l,n}-q^\ast\right|\nconverge0$ from Lemma~\ref{lemma:overlap concentration without perturb}. Building upon this intuition, $s_l(\bx)x_{n+1}^{(l)}$ will be well approximated by $s_l^0x_{n+1}^{(l)}+\frac{\lambdal-\lambdal q_{l,n}}{2}$ where
\begin{equation*}
    s_l^0:=\lambdal q_{l,n}X_{n+1}^{(l)}+\sum_{i=1}^n\sqrt{\frac{\lambdal}{n}}W_{i,n+1}^{(l)}\left\langle x_i^{(l)}\right\rangle_n.
\end{equation*}
In addition, $\delta(\bx)$ will be approximated by
\begin{equation*}
    \delta^0:=\sum_{l=1}^L \left[\frac{\lambdal q_{l,n}^2+\lambdal}{4}+\sum_{1\le i<j\le n}\frac{\sqrt{\lambdal}}{n}\widetilde{W}_{i,j}^{(l)}\left\langle x_i^{(l)}x_j^{(l)}\right\rangle_n\right].
\end{equation*}
Replacing $(\delta,s_l)$ by their respective mean-field surrogates, we define
\begin{align*}
    U_2=&\sum_{x_{n+1},y_{n+1}}p(x_{n+1},y_{n+1}|Z_{n+1})e^{-\delta^0+\sum_{l}s_l^0x_{n+1}^{(l)}+\frac{\lambdal-\lambdal q_{l,n}}{2}},\\
    V_2=&\sum_{x_{n+1},y_{n+1}}y_{n+1}^{(l_1)}p(x_{n+1},y_{n+1}|Z_{n+1})e^{-\delta^0+\sum_{l}s_l^0x_{n+1}^{(l)}+\frac{\lambdal-\lambdal q_{l,n}}{2}}.
\end{align*}
Using Lemma~\ref{lemma:overlap concentration without perturb}, we can repeat the arguments in the proof of Lemma~\ref{lemma:convergence in L2} to derive
\begin{equation*}
    \mathbb{E}\left|U_1-U_2\right|^2\nconverge0,\quad \mathbb{E}\left|V_1-V_2\right|^2\nconverge0.
\end{equation*}
In addition, repeating Lemma~\ref{lemma:bounded inverse square}, we also have $\mathbb{E}U_1^{-2}+\mathbb{E}U_2^{-2}\le C$ for some $C>0$ independent of $n$. Note that by definition, $|V_1/U_1|=\left|\left\langle y^{(l_1)}_{n+1}\right\rangle_{n+1}\right|\le 1$. Thus we have, 
\begin{align*}
    \mathbb{E}\left|\frac{V_1}{U_1}-\frac{V_2}{U_2}\right|&\le\mathbb{E}\left[\frac{V_1}{U_1}U_2^{-1}\left|U_1-U_2\right|\right]+\mathbb{E}U_2^{-1}\left|V_1-V_2\right|\\
    &\le\sqrt{\mathbb{E}U_2^{-2}} \sqrt{\mathbb{E}\left|U_1-U_2\right|^2}+\sqrt{\mathbb{E}U_2^{-2}} \sqrt{\mathbb{E}\left|V_1-V_2\right|^2}\nconverge0,
\end{align*}
which implies $\left|\mathbb{E}\left[Y_{n+1}^{(l_1)}\langle y_{n+1}^{(l_1)}\rangle\right]-\mathbb{E}Y_{n+1}^{(l_1)}V_2/U_2\right|\nconverge0$. Lastly, since $\mathbb{E}\left|q_{l,n}-q^\ast_l\right|\rightarrow0$, we conclude
\begin{equation*}
    \mathbb{E}\left[Y_{n+1}^{(l_1)}\langle y_{n+1}^{(l_1)}\rangle\right]\rightarrow\mathbb{E}\left[Y^{(l)}\mathbb{E}\left(Y^{(l_1)}|Z,\sqrt{\boldsymbol{\lambda}\odot\mathbf{q}^\ast}\odot X+W^\prime\right)\right].
\end{equation*}
Consequently, using \eqref{eq:upper bound on matrix MMSE} we have,
\begin{equation*}
    \limsup_{n\rightarrow\infty}\MMSE\left(\bY^{(l_1)};\bA,\bZ\right)\le \frac{1}{4}+\frac{q_y^2}{4}-\frac{q_y}{2}\mathbb{E}\left[Y^{(l)}\mathbb{E}\left(Y^{(l_1)}|Z,\sqrt{\boldsymbol{\lambda}\odot\mathbf{q}^\ast}\odot X+W^\prime\right)\right].
\end{equation*}
Finally, we set $q_y=\mathbb{E}\left[Y^{(l)}\mathbb{E}\left(Y^{(l_1)}|Z,\sqrt{\boldsymbol{\lambda}\odot\mathbf{q}^\ast}\odot X+W^\prime\right)\right]$ to conclude the proof.
\end{proof}

\subsection{Universality of Mutual Information}
We turn to a proof of Proposition~\ref{prop:Universality MI qualitative} in this section. Since $I(\bX,\bY;\mathbf{G},\bZ)$ (resp. $I(\bX,\bY;\mathbf{A},\bZ)$) always differs from $I(\bX,\bY;\mathbf{G}\mid\bZ)$ (resp. $I(\bX,\bY;\mathbf{A}\mid\bZ)$) with a fixed term $ni_p(x,y;z)$, it suffices to control the difference between the conditional mutual information in the graph problem (with finite $d_n$) and the mutual information in the spiked gaussian matrix problem.

\begin{Proposition}\label{thm:Universality MI quantitative}
    Under the asymptotics that (i) $\lambdal_n\rightarrow\lambdal $ for any $l\in[L]$ and (ii) $\min_l d^{(l)}(1-d^{(l)}/n)\rightarrow\infty$, the following universality result holds in terms of mutual information for the graph model~\eqref{eq:original graph model} and its spiked matrix counterpart~\eqref{eq:spiked matrix each layer},
    \begin{equation*}
        \left|\frac{1}{n}I(\bX,\bY;\mathbf{G}\mid\bZ)-\frac{1}{n}I(\bX,\bY;\bA\mid\bZ)\right|\le O\left(\sum_l|\lambdal_n-\lambdal|+\sum_{l}\frac{(\lambdal)^{3/2}}{\sqrt{d^{(l)}(1-d^{(l)}/n)}}\right).
    \end{equation*}
\end{Proposition}
Our proof strategy mirrors the one used in \cite[Proposition 4.1]{deshpande2017asymptotic} and \cite[Theorem 54]{lelarge2019fundamental}. To establish this proposition, we need two intermediate lemmas. To this end, let
\begin{equation*}
    \tau_l=\frac{a^{(l)}-b^{(l)}}{a^{(l)}+b^{(l)}}=\sqrt{\frac{\lambdal_n(1-d^{(l)}/n)}{d^{(l)}}}
\end{equation*}
measure the discrepency of in-community connectivity and cross-community connectivity for layer $l$. For any $i<j,l\in[L]$, define 
\begin{equation}\label{eq:universality MI binary noise}
    D_{i,j}^{(l)}:=\frac{\tau_l}{1-d^{(l)}/n}\left(G_{i,j}^{(l)}-\mathbb{E}\left[G_{i,j}^{(l)}\mid \bX\right]\right)=\frac{\tau_l}{1-d^{(l)}/n}\left(G_{i,j}^{(l)}-\frac{d^{(l)}}{n}-\frac{d^{(l)}\tau_l}{n}X_i^{(l)}X_j^{(l)}\right).
\end{equation}

\noindent 
 Our first lemma simplifies the definition of $I(\bX,\bY;\mathbf{G}|\bZ)$ via Taylor expansion. 

\begin{Lemma}\label{lemma:Universality MI Taylor}
    Under the same asymptotics as Proposition~\ref{thm:Universality MI quantitative}, the conditional mutual information $I(\bX,\bY;\mathbf{G}|\bZ)$ can be approximated by
    \begin{align}\label{eq:intermediate gaussian approx}
        J_1=-\mathbb{E}\log\left[\sum_{\bx,\by}p(\bx,\by|\bZ)e^{\sum_{l=1}^L\sum_{i<j}D_{i,j}^{(l)}x_i^{(l)}x_j^{(l)}+\frac{\lambdal_n}{n}X_i^{(l)}X_j^{(l)}x_i^{(l)}x_j^{(l)}-\frac{\lambdal_n}{n}}\right].
    \end{align}
    in the sense that
    \begin{equation*}
        \limsup_{n\rightarrow\infty}\frac{1}{n}\left| J_1-I(\bX,\bY;\mathbf{G}|\bZ)\right|\le O\left(\sum_l\frac{(\lambdal)^{3/2}}{\sqrt{ d^{(l)}(1-d^{(l)}/n)}}\right).
    \end{equation*}
\end{Lemma}

\noindent 
Recall now the definition of $I(\bX,\bY;\bA|\bZ)$,
\begin{equation*}
    I(\bX,\bY;\bA|\bZ)=-\mathbb{E}\log\left[\sum_{\bx,\by}p(\bx,\by|\bZ)e^{\sum_{l=1}^L\sum_{i<j}\sqrt{\frac{\lambda^{(l)}}{n}}W_{i,j}^{(l)}x_i^{(l)}x_j^{(l)}+\frac{\lambda^{(l)}}{n}X_i^{(l)}X_j^{(l)}x_i^{(l)}x_j^{(l)}-\frac{\lambdal}{n}}\right].
\end{equation*}
By direct comparison, we find that  $J_1$ differs from $I(\bX,\bY;\bA|\bZ)$ by replacing $\sqrt{\frac{\lambdal}{n}}W_{i,j}^{(l)}$ by $D_{i,j}^{(l)}$ and $\lambdal$ by $\lambdal_n$. Our next lemma uses Lindeberg swapping to show that this change is asymptotically negligible.

\begin{Lemma}\label{lem:Universality MI Lindeberg}
    Under the same asymptotics as Proposition~\ref{thm:Universality MI quantitative}, we have
    \begin{equation*}
        \limsup_{n\rightarrow\infty}\frac{1}{n}\left| J_1-I(\bX,\bY;\bA|\bZ)\right|\le O\left(\sum_l|\lambdal_n-\lambdal|+\frac{(\lambdal)^{3/2}}{\sqrt{ d^{(l)}(1-d^{(l)}/n)}}\right).
    \end{equation*}
\end{Lemma}

\noindent
Proposition~\ref{thm:Universality MI quantitative} directly follows from combining Lemmas~\ref{lemma:Universality MI Taylor} and \ref{lem:Universality MI Lindeberg}.

\paragraph{Taylor Expansion of Bernoulli p.m.f.}
We start with a Chernoff-style bound on the upper tail of the total number of edges. The proof is standard, and thus omitted. 
\begin{proof}[Proof of Lemma~\ref{lemma:Universality MI Taylor}]
    By definition, $I(\bX,\bY;\mathbf{G}|\bZ)=-\mathbb{E} \log \frac{p(\mathbf{G}|\mathbf{Z})}{p(\mathbf{G}\mid\mathbf{X},\mathbf{Y},\bZ)}$. Thus 
    \begin{equation}
        I(\bX,\bY;\mathbf{G}|\bZ)=-\mathbb{E} \log \frac{\sum_{\bx,\by} p(\mathbf{x},\mathbf{y}|\bZ)p(\mathbf{G}|\mathbf{x},\mathbf{y})}{p(\mathbf{G}\mid\mathbf{X},\mathbf{Y})}.
    \end{equation}
    Recall that $p(\mathbf{G}|\mathbf{x},\mathbf{y})=\prod_{l}p(\mathbf{G}^{(l)}|\mathbf{x}^{(l)},\mathbf{y}^{(l)})$ with
    \begin{align}
        &p(\mathbf{G}^{(l)}|\mathbf{x}^{(l)},\mathbf{y}^{(l)})=\prod_{i<j} M(x_i^{(l)}, x_j^{(l)})^{G_{i, j}^{(l)}}\left(1-M(x_i^{(l)}, x_j^{(l)})\right)^{1-G_{i, j}^{(l)}}\\
        =&\exp \left[\sum_{i<j} G_{i, j}^{(l)} \log M(x_i^{(l)}, x_j^{(l)})+\left(1-G_{i, j}^{(l)}\right) \log \left(1-M(x_i^{(l)}, x_j^{(l)})\right)\right]
    \end{align}
    where we set  $M(x_i^{(l)}, x_j^{(l)})=\frac{d^{(l)}}{n}\left(1+x_i^{(l)} x_j^{(l)} \tau_l\right)$. Equipped with this notation,
    \begin{align*}
        I(\bX,\bY;\mathbf{G}|\bZ)=-\mathbb{E}\log \Bigg[\sum_{\bx,\by}p(\bx,\by|\bZ) \exp \Bigg(\sum_{l=1}^L\sum_{i<j} &G_{i, j}^{(l)} \log \left(\frac{M(x_i^{(l)}, x_j^{(l)})}{M(X_i^{(l)}, X_j^{(l)})}\right)\\
        &+\left(1-G_{i, j}^{(l)}\right) \log \left(\frac{1-M(x_i^{(l)}, x_j^{(l)})}{1-M(X_i^{(l)}, X_j^{(l)})}\right)\Bigg)\Bigg].
    \end{align*}
    Note that $\log \left(\frac{M(x_i^{(l)}, x_j^{(l)})}{M(X_i^{(l)}, X_j^{(l)})}\right)=\log \left(\frac{1+x_i^{(l)} x_j^{(l)}\tau_l}{1+X_i^{(l)} X_j^{(l)} \tau_l}\right)$. By the Taylor-Lagrange inequality, there exists a universal constant $C>0$ such that, for $\tau_l$ small enough (i.e. for $d$ large enough),
    \begin{align*}
        \left|\log \left(\frac{1+x_i^{(l)} x_j^{(l)}\tau_l}{1+X_i^{(l)} X_j^{(l)} \tau_l}\right)-\tau_l\left(x_i^{(l)} x_j^{(l)}-X_i^{(l)} X_j^{(l)}\right)\right| \leq C \tau_l^3.
    \end{align*}
    Moreover, since $\log(c+dx)=\frac{1}{2}\log\left((c+d)(c-d)\right)+\frac{x}{2}\log\frac{c+d}{c-d}$ for any $x\in\{\pm 1\}$, we find
    \begin{align*}
        \log \left(\frac{1-M(x_i^{(l)}, x_j^{(l)})}{1-M(X_i^{(l)}, X_j^{(l)})}\right)=\frac{1}{2}\left(x_i^{(l)} x_j^{(l)}-X_i^{(l)} X_j^{(l)}\right)\log\frac{1-d^{(l)}/n-\tau_ld^{(l)}/n}{1-d^{(l)}/n+\tau_ld^{(l)}/n}.
    \end{align*}
    Taylor expansion also implies $\left|\frac{1}{2}\log\left(\frac{1+z}{1-z}\right)-z\right|\le z^3$, which further deduces
    \begin{align*}
        \left|\log \left(\frac{1-M(x_i^{(l)}, x_j^{(l)})}{1-M(X_i^{(l)}, X_j^{(l)})}\right)+\frac{\tau_ld^{(l)}/n}{1-d^{(l)}/n} \left(x_i^{(l)} x_j^{(l)}-X_i^{(l)} X_j^{(l)}\right)\right| \leq C \left(\frac{\tau_ld^{(l)}/n}{1-d^{(l)}/n}\right)^3.
    \end{align*}
    Summing over $i<j$  and using triangle inequality, we have
    \begin{align*}
        &\sum_{i<j}G_{i, j}^{(l)} \log \left(\frac{M(x_i^{(l)}, x_j^{(l)})}{M(X_i^{(l)}, X_j^{(l)})}\right)+\left(1-G_{i, j}^{(l)}\right) \log \left(\frac{1-M(x_i^{(l)}, x_j^{(l)})}{1-M(X_i^{(l)}, X_j^{(l)})}\right)\\
        =&\sum_{i<j} \left(\tau_l\left(x_i^{(l)} x_j^{(l)}-X_i^{(l)} X_j^{(l)}\right)G_{i,j}^{(l)}-\left(1-G_{i, j}^{(l)}\right)\left(x_i^{(l)} x_j^{(l)}-X_i^{(l)} X_j^{(l)}\right)\frac{\tau_ld^{(l)}/n}{1-d^{(l)}/n}\right)+\Delta_n^{(l)}\\
        =&\sum_{i<j}\frac{\tau_l}{1-d^{(l)}/n}\left(x_i^{(l)} x_j^{(l)}-X_i^{(l)} X_j^{(l)}\right)\left(G_{i, j}^{(l)}-\frac{d^{(l)}}{n}\right)+\Delta_n^{(l)}\\
        =&\sum_{i<j}\left(\left(x_i^{(l)} x_j^{(l)}-X_i^{(l)} X_j^{(l)}\right)D_{i, j}^{(l)}+\frac{\lambdal_n}{n}x_i^{(l)} x_j^{(l)}X_i^{(l)} X_j^{(l)}-\frac{\lambdal_n}{n}\right)+\Delta_n^{(l)},
    \end{align*}
    where $\Delta_n^{(l)}$ is a residual term satisfying 
    \begin{equation}\label{eq:universality MI Hamilton err bound}
        \left|\Delta_n^{(l)}\right| \leq C\left[\sum_{i<j} G_{i, j}^{(l)} \tau_l^3+(1-G_{i,j}^{(l)})\left(\frac{\tau_ld^{(l)}/n}{1-d^{(l)}/n}\right)^3\right],
    \end{equation}
    where $C>0$ is a universal constant. Consequently, we can express
    \begin{align*}
        I(\bX,\bY;\mathbf{G}|\bZ)=&-\mathbb{E}\log \Bigg[\sum_{\bx,\by}p(\bx,\by|\bZ) \exp \Bigg(\Delta_n^{(l)}\\
        &+\sum_{l}\sum_{i<j}\left(\left(x_i^{(l)} x_j^{(l)}-X_i^{(l)} X_j^{(l)}\right)D_{i, j}^{(l)}+\frac{\lambdal_n}{n}x_i^{(l)} x_j^{(l)}X_i^{(l)} X_j^{(l)}-\frac{\lambdal_n}{n}\right.\Bigg)\Bigg].
    \end{align*}
    Define $\mathsf{err}=I(\bX,\bY;\mathbf{G}|\bZ)-J_1$ with $J_1$ from \eqref{eq:intermediate gaussian approx}. For any fixed realization $(\bX,\bY,\bZ,\bG)$, error term $\left|\Delta_n^{(l)}\right|$ has a uniform upper bound as shown in \eqref{eq:universality MI Hamilton err bound} for any possible $(\bx,\by)$, so it follows that
    \begin{align*}
        |\mathsf{err}|&\leq C\E\left[\sum_{l}\sum_{i<j} G_{i, j}^{(l)} \tau_l^3+(1-G_{i,j}^{(l)})\left(\frac{\tau_ld^{(l)}/n}{1-d^{(l)}/n}\right)^3\right],\\
        &\leq O\left(\sum_{l}\tau_l^3\mathbb{E}\sum_{i<j}G_{i,j}^{(l)}\right)+O\left(\sum_{l}\left(\frac{\tau_ld^{(l)}/n}{1-d^{(l)}/n}\right)^3\mathbb{E}\sum_{i<j}(1-G_{i,j}^{(l)})\right)\\
        &\leq O\left(\sum_l\frac{n(\lambdal_n)^{3/2}(1-d^{(l)}/n)^{3/2}}{\sqrt{d^{(l)}}}\right)+O\left(\sum_l\frac{(\lambdal_n)^{3/2}(d^{(l)})^{3/2}}{n\sqrt{1-d^{(l)}/n}}\right)\\
        &\le O\left(n\sum_l\frac{(\lambdal)^{3/2}}{\sqrt{ d^{(l)}(1-d^{(l)}/n)}}\right).
    \end{align*}
    where the third inequality uses the direct fact that,
    \begin{equation*}
        \mathbb{E}\sum_{i<j}G_{i,j}^{(l)}\leq Cnd^{(l)},\quad \mathbb{E}\sum_{i<j}\left(1-G_{i,j}^{(l)}\right)\leq Cn^2(1-d^{(l)}/n);
    \end{equation*}
    while the final inequality uses $\lambdal_n\rightarrow\lambdal$, $d^{(l)}\le n$ and $1-d^{(l)}/n\le1$. Until this end, we can conclude Lemma~\ref{lemma:Universality MI Taylor}.
\end{proof}

\paragraph{Lindeberg Exchange}
The following generalization of Lindeberg's theorem is adapted from \cite[Theorem 2]{korada2011applications}.
\begin{Lemma}[\cite{korada2011applications}]\label{lem:Lindeberg}
    Let $\left(U_i\right)_{1 \leq i \leq n}$ and $\left(V_i\right)_{1 \leq i \leq n}$ be two collection of random variables with independent components and $f: \mathbb{R}^n \rightarrow \mathbb{R}$ a $\mathcal{C}^3$ function. Denote $a_i=\left|\mathbb{E} U_i-\mathbb{E} V_i\right|$ and $b_i=\left|\mathbb{E} U_i^2-\mathbb{E} V_i^2\right|$. Then
    \begin{align*}
       & |\mathbb{E} f(U)-\mathbb{E} f(V)|  \leq \sum_{i=1}^n\left(a_i \mathbb{E}\left|\partial_i f\left(U_{1: i-1}, 0, V_{i+1: n}\right)\right|+\frac{b_i}{2} \mathbb{E}\left|\partial_i^2 f\left(U_{1: i-1}, 0, V_{i+1: n}\right)\right|\right. \\
        & \left.+\frac{1}{2} \mathbb{E} \int_0^{U_i}\left|\partial_i^3 f\left(U_{1: i-1}, 0, V_{i+1: n}\right)\right|\left(U_i-s\right)^2 d s+\frac{1}{2} \mathbb{E} \int_0^{V_i}\left|\partial_i^3 f\left(U_{1: i-1}, 0, V_{i+1: n}\right)\right|\left(V_i-s\right)^2 d s\right).
    \end{align*}
\end{Lemma}

\begin{proof}[Proof of Lemma~\ref{lem:Universality MI Lindeberg}]
    Conditioned on $(\bX,\bY,\bZ)$, introduce the following function
    \begin{equation*}
        \Phi(\mathbf{u},\boldsymbol{\lambda})=-\log\left[\sum_{\bx,\by}p(\bx,\by|\bZ)e^{\sum_{l=1}^L\sum_{i<j}u_{i,j}^{(l)}x_i^{(l)}x_j^{(l)}+\frac{\lambda^{(l)}}{n}X_i^{(l)}X_j^{(l)}x_i^{(l)}x_j^{(l)}-\frac{\lambdal}{n}}\right]
    \end{equation*}
    for each aggregation $\mathbf{u}=\left\{u_{i,j}^{(l)}:l\in[L],i<j\right\}$ and $\boldsymbol{\lambda}=\left\{\lambdal:l\in[L]\right\}$.
    Function $\Phi(\mathbf{u},\boldsymbol{\lambda})$ is thus $\mathcal{C}^3$ in $\mathbf{u}$ with bounded derivatives since the support of $\mathbf{X}$ and $\mathbf{x}$ are bounded. Notice that $J_1=\mathbb{E} \Phi(\mathbf{D},\boldsymbol{\lambda}_n)$ and $I(\bX,\bY;\bA|\bZ)=\mathbb{E} \Phi\left(\sqrt{\frac{\boldsymbol{\lambda}}{n}} \mathbf{W},\boldsymbol{\lambda}\right)$. We then proceed to compute the moments of $D_{i,j}^{(l)}$ whose definition is from \eqref{eq:universality MI binary noise}, conditionally to $\mathbf{X}$. Firstly, $D_{i,j}^{(l)}$ is also of mean zero like $\sqrt{\lambdal/n}W_{i, j}^{(l)}$,
    \begin{equation*}
        \mathbb{E}\left(D_{i, j}^{(l)} \mid \mathbf{X}\right)=\mathbb{E}\left(\sqrt{\lambdal/n}W_{i, j}^{(l)} \right)=0.
    \end{equation*}
    Secondly, its second moment can be computed as
    \begin{align*}
        &\quad\mathbb{E}\left(D_{i, j}^{(l)2} \mid \mathbf{X}\right) =\frac{\tau_l^2}{(1-d^{(l)}/n)^2} \operatorname{Var}\left(G_{i, j}^{(l)} \mid \mathbf{X}\right)\\
        &=\frac{\lambdal_n/n}{d^{(l)}(1-d^{(l)}/n)}\left[a^{(l)}(1-a^{(l)}/n)\mathbbm{1}\{X_i^{(l)}=X_j^{(l)}\}+b^{(l)}(1-b^{(l)}/n)\mathbbm{1}\{X_i^{(l)}\neq X_j^{(l)}\}\right],
    \end{align*}
    while $\mathbb{E}\left(\sqrt{\lambdal/n}W_{i, j}^{(l)} \right)^2=\lambdal/n$. Under the asymptotics that (i) $\lambdal_n\rightarrow\lambdal $ for any $l\in[L]$ and (ii) $\min_l d^{(l)}(1-d^{(l)}/n)\rightarrow\infty$, we have
    \begin{equation*}
        \left|\frac{a^{(l)}(1-a^{(l)}/n)}{d^{(l)}(1-d^{(l)}/n)}-1\right|\le C\sqrt{\frac{\lambdal}{d^{(l)}(1-d^{(l)}/n)}};
    \end{equation*}
    and the same holds as well when $a^{(l)}$ is replaced by $b^{(l)}$. Consequently, differences in the second moments can be controlled by
    \begin{equation*}
        \left|\mathbb{E}\left(D_{i, j}^{(l)2} \mid \mathbf{X}\right)-\frac{\lambdal}{n}\right|\le\frac{|\lambdal_n-\lambdal|}{n}+\frac{C}{n}\frac{(\lambdal)^{3/2}}{\sqrt{d^{(l)}(1-d^{(l)}/n)}}.
    \end{equation*}
    The third moments is also bounded by $\mathbb{E}\left(D_{i, j}^{(l)3} \mid \mathbf{X}\right)=O\left(\frac{1}{n}\frac{(\lambdal)^{3/2}}{\sqrt{d^{(l)}(1-d^{(l)}/n)}}\right)$. From Lemma~\ref{lem:Lindeberg} we obtain
    \begin{align*}
        &\quad\left| J_1-\mathbb{E} \Phi\left(\sqrt{\frac{\boldsymbol{\lambda}}{n}} \mathbf{W},\boldsymbol{\lambda}_n\right)\right|=\left|\mathbb{E} \Phi\left( \mathbf{D},\boldsymbol{\lambda}_n\right)-\mathbb{E} \Phi\left(\sqrt{\frac{\boldsymbol{\lambda}}{n}} \mathbf{W},\boldsymbol{\lambda}_n\right)\right| \\
        &\leq \sum_l\sum_{i<j} O\left(\frac{|\lambdal_n-\lambdal|}{n}+\frac{1}{n}\frac{(\lambdal)^{3/2}}{\sqrt{d^{(l)}(1-d^{(l)}/n)}}\right)\\
        &=O\left(n\sum_l|\lambdal_n-\lambdal|+n\sum_{l}\frac{(\lambdal)^{3/2}}{\sqrt{d^{(l)}(1-d^{(l)}/n)}}\right).
    \end{align*}
    The last step is to note that $|\partial_{\lambdal}\Phi(\mathbf{u},\boldsymbol{\lambda})|\le Cn$ always holds for some universal constant $C>0$ since $\bX,\bx$ are bounded. Followed by
    \begin{equation*}
        \left|\mathbb{E} \Phi\left( \sqrt{\frac{\boldsymbol{\lambda}}{n}} \mathbf{W},\boldsymbol{\lambda}\right)-\mathbb{E} \Phi\left(\sqrt{\frac{\boldsymbol{\lambda}}{n}} \mathbf{W},\boldsymbol{\lambda}_n\right)\right|\le O\left(n\sum_l|\lambdal_n-\lambdal|\right),
    \end{equation*}
    Lemma~\ref{lem:Universality MI Lindeberg} is finally established.
\end{proof}

\subsection{Universality of Minimal MSE}
We establish Proposition~\ref{prop:universality mmse} in this section. 
To this end, the following lemma connects the graph models minimal MSE to the derivative of the mutual information. 
\begin{Lemma}
\label{lem:MI_derivative}
    There exists a universal constant $C>0$ such that: for any fixed $\lambda_0$, there exists $n_0(\lambda_0)$ such that for all $0<\lambdal\le \lambda_0$, and $n\ge n_0(\lambda_0)$,
    \begin{equation*}
    \left|\frac{1}{n}\frac{\partial I(\bX,\bY;\mathbf{G}\mid\bZ)}{\partial\lambdal}-\frac{1}{4}\MMSE\left(\bX^{(l)} \bX^{(l)\top};\bG,\bZ\right)\right|\le C\left(\sqrt{\frac{\lambdal}{d^{(l)}(1-d^{(l)}/n)}}\lor\frac{1}{n}\right).
    \end{equation*}    
\end{Lemma}
This lemma is in the same spirit as \cite[Lemma 7.2]{deshpande2017asymptotic} and \cite[Proposition 62]{lelarge2019fundamental}. Consequently, we omit the proof. In the following, to emphasize the dependence on $\lambda^{(l)}$, we write
\begin{align*}
    \MMSE_l^\mathsf{G}(\lambda^{(l)},n)&=\MMSE\left(\bX^{(l)}\bX^{(l)\top};\bG,\bZ\right),\\
    \Phi(\lambda^{(l)})=i_p(x,y;z)+&\frac{\sum_{l^\prime}\lambda^{(l^\prime)}}{4}-\sup_{\mathbf{q}\ge 0}\left[\mathcal{F}(\boldsymbol{\lambda},\mathbf{q})-\sum_{l^\prime=1}^{L}\frac{\lambda^{(l^\prime)}(q_{l^\prime}^2+2q_{l^\prime})}{4}\right],
\end{align*}
where the last term appears in \eqref{eq:limiting MI}. As defined in \eqref{eq:dummy mmse in general model}, we are comparing $\MMSE_l(\lambdal)$ with the best possible estimate only from per-vertex side information $\bZ$,
\begin{align*}
    \DMSE_l(p):&=1-\left\{\E\left[X^{(l)}\E\left(X^{(l)}|Z\right)\right]\right\}^2,
\end{align*}
which only depends on the prior $p$ and is irrelevant of $n$. As the free energy is convex in $\lambda^{(l)}$, using Theorem~\ref{thm:asymptotic MI}, for all but countably many $\lambdal$,
\begin{equation*}
    \Phi^\prime(\lambdal)=\lim_{n\rightarrow\infty}\frac{1}{4}\MMSE\left(\bX^{(l)}\bX^{(l)\top};\bA,\bZ\right).
\end{equation*}
Using Lemma \ref{lem:MI_derivative}, under the asymptotics that $\boldsymbol{\lambda}_n\rightarrow\boldsymbol{\lambda}$ and $d^{(l)}(1-d^{(l)}/n)\rightarrow\infty$, we have  
\begin{align}
    &\quad\limsup_{n\rightarrow\infty}\left|\frac{1}{4}\int_{0}^{\lambda^{(l)}}\MMSE_l^\mathsf{G}(\lambda^\prime,n)\mathrm{d}\lambda^\prime-\Phi(\lambda^{(l)})\right|\notag\\
    &\le \limsup_{n\rightarrow\infty}\left|\frac{1}{4}\int_{0}^{\lambda^{(l)}}\MMSE_l^\mathsf{G}(\lambda^\prime,n)\mathrm{d}\lambda^\prime-\frac{1}{n}I(\bX,\bY;\mathbf{G}\mid\bZ)\right|+\left|\frac{1}{n}I(\bX,\bY;\mathbf{G}\mid\bZ)-\Phi(\lambda^{(l)})\right|\notag\\
    &\le O\left(|\lambdal_n-\lambdal|+\frac{(\lambdal)^{3/2}}{\sqrt{d^{(l)}(1-d^{(l)}/n)}}\lor\frac{1}{n}\right),\label{eq:upper bound mmse universality integration}
\end{align}
where for last inequality we also use Proposition~\ref{prop:Universality MI qualitative} to upper bound $\left|\frac{1}{n}I(\bX,\bY;\mathbf{G}\mid\bZ)-\Phi(\lambda^{(l)})\right|$.

\begin{proof}[Proof of Proposition~\ref{prop:universality mmse}]
(i) If $\blam$ is below the weak recovery threshold for the graph model, such that 
\begin{equation*}
    \lim_{n\rightarrow\infty}\MMSE_l^\mathsf{G}(\lambda^{(l)},n)=\DMSE_l(p),
\end{equation*}
it follows from monotonicity that $\MMSE_l^\mathsf{G}(\lambda^{\prime},n)=\DMSE_l(p)$ remains constant for any $0\le\lambda^\prime\le\lambdal$. Thereafter, \eqref{eq:upper bound mmse universality integration} deduces $\Phi(\lambda^\prime)=\frac{\lambda^\prime}{4}\DMSE_l(p)$ to be linear in $0\le\lambda^\prime\le\lambdal$. It follows
\begin{equation*}
    \lim_{n\rightarrow\infty}\frac{1}{4}\MMSE\left(\bX^{(l)}\bX^{(l)\top};\bA,\bZ\right)=\Phi^\prime(\lambdal)=\DMSE_l(p).
\end{equation*}
Consequently, $\boldsymbol{\lambda}$ is also below the weak recovery threshold for the spiked matrix model. (ii) On the other hand, if $\blam$ is above the weak recovery threshold for the graph model, such that
\begin{equation}\label{eq:mmse graph smaller than dmse}
    \lim_{n\rightarrow\infty}\MMSE_l^\mathsf{G}(\lambda^{(l)},n)<\DMSE_l(p).
\end{equation}
We then proceed by contradiction. Suppose $\blam$ is below the weak recovery threshold for the spiked matrix model, i.e.
\begin{equation*}
    \limsup_{n\rightarrow\infty}\MMSE\left(\bX^{(l)}\bX^{(l)\top};\bA,\bZ\right)=\DMSE_l(p).
\end{equation*}
By definition, we must know $\Phi(\lambda^\prime)=\frac{\lambda^\prime}{4}\DMSE_l(p)$ to be linear in $0\le\lambda^\prime\le\lambdal$. But instead, \eqref{eq:mmse graph smaller than dmse} and \eqref{eq:upper bound mmse universality integration} together suggest $\Phi(\lambdal)<\frac{\lambdal}{4}\DMSE_l(p)$, resulting in a contradiction. Therefore, we must find $\blam$ above the weak recovery threshold for the spiked matrix model. Combining these two assertions yields the weak recovery thresholds to be universal from spiked matrices to random graphs, thereby proving Proposition~\ref{prop:universality mmse}.
\end{proof}

\section{Approximate Message Passing}
\label{sec:AMP_derivation}

We derive the AMP algorithm (Algorithm \ref{alg:coupled AMP}) in Section~\ref{subsec:heurestic derivation}. Subsequently, we characterize the state evolution behavior of this algorithm in Section \ref{sec:state_evolution}. Finally, we establish the universality of the AMP algorithm in Section \ref{sec:alg_universality}.   

\subsection{Heuristic Derivations}\label{subsec:heurestic derivation}
To simplify notations, we omit the dependence of the posterior on $(\bZ,\bA)$ for this subsection. Specifically, we use
\begin{equation*}
    p_i(x_i,y_i):=p(x_i,y_i|Z_i),
\end{equation*}
and denote the joint posterior $p(\bx,\by|\bZ,\bA)$ by
\begin{equation*}
    \mu(\bx,\by)\propto \prod_{i}p_i(x_i,y_i)\cdot\exp\left\{\sum_{l=1}^L\sqrt{\frac{\lambda^{(l)}}{n}}A^{(l)}_{i,j}x_i^{(l)}x_j^{(l)}\right\}.
\end{equation*}
With a factorized base measure $\otimes_{i} p_i$ on $(\bx,\by)$, the posterior is a Gibbs measure with Hamiltonian $H(\bx):=\sum_{i<j}\sum_{l=1}^L\sqrt{\frac{\lambda^{(l)}}{n}}A^{(l)}_{i,j}x_i^{(l)}x_j^{(l)}$.

We treat each $(x_i,y_i)$ jointly as a spin variable, and each marginal is denoted as $\mu_i$. For any disjoint $A,B\subset[n]$, we set $\nu_{A\rightarrow B}(\bx_A,\by_A)$ as the marginal distribution of $(\bx_A,\by_A)$ with $(\bx_B,\by_B)$ removed. Due to the separability of $H(\bx)$, we have, 
\begin{align*}
    &\mu_{i}\left(x_i,y_i\right)/p_i(x_i,y_i)\\
    \propto &\sum_{\bx_{[n]\backslash i},\by_{[n]\backslash i}}\prod_{k\neq i}p_i(x_k,y_k)\cdot\exp\left\{\sum_{k<j:k,j\in[n]\backslash i}\sum_{l=1}^L\sqrt{\frac{\lambda^{(l)}}{n}}A^{(l)}_{k,j}x_k^{(l)}x_j^{(l)}+\sum_{j\in[n]\backslash i}\sum_{l=1}^L\sqrt{\frac{\lambda^{(l)}}{n}}A^{(l)}_{i,j}x_i^{(l)}x_j^{(l)}\right\}\\
    \propto&\sum_{\bx_{[n]\backslash i},\by_{[n]\backslash i}}\nu_{[n]\backslash i\rightarrow i}(\bx_{[n]\backslash i},\by_{[n]\backslash i})\cdot\exp\left\{\sum_{j\in[n]\backslash i}\sum_{l=1}^L\sqrt{\frac{\lambda^{(l)}}{n}}A^{(l)}_{i,j}x_i^{(l)}x_j^{(l)}\right\}.
\end{align*}
We invoke the replica symmetry assumption to continue our derivation : cavity distribution $\nu_{[n]\backslash i\rightarrow i}(\bx_{[n]\backslash i},\by_{[n]\backslash i})$ is approximately independent across the spins, namely $\nu_{[n]\backslash i\rightarrow i}(\bx_{[n]\backslash i},\by_{[n]\backslash i})=\prod_{j\in[n]\backslash i}\nu_{j\rightarrow i}(x_j,y_j)+o(1)$. In turn, this implies
\begin{equation}\label{eq:Gibbs marginal decomposition}
    \mu_{i}\left(x_i,y_i\right)\propto p_i\left(x_i,y_i\right)\sum_{\bx_{[n]\backslash i},\by_{[n]\backslash i}}\prod_{j\in[n]\backslash i}\nu_{j\rightarrow i}(x_j,y_j)\exp\left\{\sum_{l=1}^L\sqrt{\frac{\lambda^{(l)}}{n}}A^{(l)}_{i,j}x_i^{(l)}x_j^{(l)}\right\}+o(1).
\end{equation}

\paragraph{Parametrization.} To proceed, we need a simple parametrization for any distribution on $(x,y)\in\{\pm 1\}^{L+L_1}$. Inspired by synchronization problems on unitary groups \cite{perry2018message} and Boolean analysis \cite{o2014analysis}, we consider the following $2^{L+L_1}$ different functions on $\{\pm 1\}^{L+L_1}$, which are indexed by each subset $I\subset[L+L_1]$,
\begin{equation*}
    \mathsf{R}_I(u^{(1)},\ldots,u^{(L+L_1)})=\prod_{l\in I}u^{(l)},\quad\forall u^{(l)}=\pm 1.
\end{equation*}
For the rest of this section, we view $\{\pm 1\}^{L+L_1}$ as a discrete Abelian group. Consequently, these functions turn out to be all the irreducible group representations. 

\vspace{2mm}
\noindent\textit{Fact 1.} Our Hamiltonian can be well encoded into this set of representations $\{\mathsf{R}_I\}$.
Specifically, we focus on $\mathsf{R}_{I}$ when $I$ is a singleton subset of $[L]$. For simplicity, denote function $\mathsf{R}_l(u)=u^{(l)}$ for each $l\in[L]$. It follows that $x_i^{(l)}x_j^{(l)}=\mathsf{R}_l\left((x_i,y_i)(x_j,y_j)^{-1}\right)$. Consequently, we have
\begin{align*}
    H\left(\bx,\by\right)=\sum_{i<j}\sum_{l=1}^L\sqrt{\frac{\lambda^{(l)}}{n}}A^{(l)}_{i,j}\mathsf{R}_l\left((x_i,y_i)(x_j,y_j)^{-1}\right).
\end{align*}

\vspace{2mm}
\noindent\textit{Fact 2.} Placing a uniform measure on $\{\pm 1\}^{L+L_1}$, $\{2^{-\frac{L+L_1}{2}}\mathsf{R}_I:I\subset[L+L_1]\}$ readily forms an orthonormal basis for the space $L^2(\{\pm 1\}^{L+L_1})$ of square integrable functions. Therefore, any log probability mass function $\log \pi(x,y)$ on $\{\pm 1\}^{L+L_1}$ can be parametrized by its coefficients when expanded onto this basis,
\begin{equation*}
    \log \pi(x,y)=\sum_{I}r_I\mathsf{R}_I(x,y).
\end{equation*}
In conclusion, any distribution on $\{\pm 1\}^{L+L_1}$, is parametrized by $2^{L+L_1}$ real numbers $\{r_I\}$ plus one normalizing constraint $\sum_{u\in\{\pm 1\}^{L+L_1}}\exp\left[\sum_{I}r_I\mathsf{R}_I(u)\right]=1$. Due to the orthonomality of these functions, each coefficient admits a simple expression
\begin{equation*}
    r_I=\frac{1}{2^{L+L_1}}\sum_{x,y}\mathsf{R}_I(x,y)\log\pi(x,y).
\end{equation*}
For each log prior distribution $\log p_i$ on $(x_i,y_i)\in\{\pm 1\}^{L+L_1}$, it is convenient to record its coefficients by $r^{p_i}_{I}=2^{-L-L_1}\sum_{x_i,y_i}\mathsf{R}_I(x_i,y_i)\log p_i(x_i,y_i)$ for each $I$. For example, 
in the multilayer inhomogeneous SBM from Example~\ref{eg:multilayer SBM}, every $p_i=p$ is the same because there is no side information $\bZ$ in this model. In this case, we find
    \begin{equation*}
        r^{p}_{\{l,L+1\}}=\frac{1}{2}\log\left(\frac{1-\rho}{\rho}\right),\quad\forall \,\, l\in[L],
    \end{equation*}
    and $r^p_I=0$ for any other  $I$.
    Similarly, for the dynamic SBM introduced in Example~\ref{eg:dynamic SBM}, every $p_i=p$. In this case, we find
    \begin{equation*}
        r^{p}_{\{l,l+1\}}=\frac{1}{2}\log\left(\frac{1-\rho}{\rho}\right),\quad\forall l\in[L-1],
    \end{equation*}
    and $r^p_I=0$ for other $I$.

\paragraph{Most Parameters are Redundant.} 
Returning to \eqref{eq:Gibbs marginal decomposition}, we parametrize the cavity distributions by $\log \mu_i=\sum_{I}r_{i,I}\mathsf{R}_I$ and $\log \nu_{i\rightarrow j}=\sum_{I}r_{i\rightarrow j,I}\mathsf{R}_I$.
For any $I\neq\emptyset$, we compute coefficient $r_{i,I}$ of $\log \mu_i$ from \eqref{eq:Gibbs marginal decomposition},
\begin{align*}
    r_{i,I}=&\sum_{x_i,y_i} \frac{\mathsf{R}_I(x_i,y_i)}{2^{L+L_1}}\log \mu_i(x_i,y_i)\\
    =&\sum_{x_i,y_i} \frac{\mathsf{R}_I(x_i,y_i)}{2^{L+L_1}}\left(\log p_i(x_i,y_i)+\sum_{j \neq i} \log\left\{\sum_{x_j,y_j}\nu_{j\rightarrow i}(x_j,y_j)\exp\left[\sum_{l=1}^L\sqrt{\frac{\lambda^{(l)}}{n}}A^{(l)}_{i,j}x_i^{(l)}x_j^{(l)}\right]\right\}\right)\\
    =&r^{p_i}_I+\sum_{j\neq i}\sum_{x_i,y_i} \frac{\mathsf{R}_I(x_i,y_i)}{2^{L+L_1}}\log\left\{1+\sum_{x_j,y_j}\nu_{j\rightarrow i}(x_j,y_j)\sum_{l=1}^L\sqrt{\frac{\lambda^{(l)}}{n}}A^{(l)}_{i,j}x_i^{(l)}x_j^{(l)}+o\Big(\frac{1}{n}\Big)\right\}\\
    =&r^{p_i}_I+\sum_{j\neq i}\sum_{x_i,y_i} \frac{\mathsf{R}_I(x_i,y_i)}{2^{L+L_1}}\sum_{x_j,y_j}\nu_{j\rightarrow i}(x_j,y_j)\sum_{l=1}^L\sqrt{\frac{\lambda^{(l)}}{n}}A^{(l)}_{i,j}x_i^{(l)}x_j^{(l)}+o(1)\\
    =&r^{p_i}_I+\sum_{j\neq i}\sum_{l=1}^L \mathbbm{1}_{I=\{l\}}\sqrt{\frac{\lambda^{(l)}}{n}}A^{(l)}_{i,j}\left[\sum_{x_j,y_j}\nu_{j\rightarrow i}(x_j,y_j)x_j^{(l)}\right]+o(1).
\end{align*}
The first equation is due to the orthonormality of $\{2^{-\frac{L+L_1}{2}}\mathsf{R}_I:I\subset\{0,1,\ldots,L\}\}$ under uniform measure. The second equation follows by plugging in \eqref{eq:Gibbs marginal decomposition}. Since we set $I\neq\emptyset$,  $\sum_{x_i,y_i}\mathsf{R}_I(x_i,y_i)=0$ and thus the  normalizing constant in \eqref{eq:Gibbs marginal decomposition} can be ignored. The third display follows from $\exp(x)=1+x+o(x^2)$, while the forth display uses $\log(1+x)=x+o(1)$. The last equation also uses orthogonality,
\begin{equation*}
    \sum_{x_i,y_i} \frac{\mathsf{R}_I(x_i,y_i)}{2^{L+1}}x_i^{(l)}=\sum_{x_i,y_i} \frac{\mathsf{R}_I(x_i,y_i)R_l(x_i,y_i)}{2^{L+1}}=\mathbbm{1}_{I=\{l\}}.
\end{equation*}

As a result, once $I\neq \emptyset,\{1\},\ldots,\{L\}$, corresponding coefficient $r_{i,I}=r^{p_i}_I+o(1)$ is approximately constant. On the other hand,  for $I=\{l\}$ we have
\begin{equation}
    r_{i,l}=r^{p_i}_{l}+\sum_{j\neq i}\sqrt{\frac{\lambda^{(l)}}{n}}A^{(l)}_{i,j}\left[\sum_{x_j,y_j}\nu_{j\rightarrow i}(x_j,y_j)x_j^{(l)}\right]+o(1).
\end{equation}
Lastly, $r_{i,\emptyset}$ adjusts automatically to normalize the distribution.

\paragraph{Linearized Belief Propagation.} Leaving out one more spin in previous equations, one can derive for any $i \neq j$ that
\begin{equation*}
    r_{i\rightarrow j,l}=r^{p_i}_{l}+\sum_{k\neq i,j}\sqrt{\frac{\lambda^{(l)}}{n}}A^{(l)}_{i,k}\left[\sum_{x_k,y_k}\nu_{k\rightarrow i}(x_k,y_k)x_k^{(l)}\right]+o(1),
\end{equation*}
and $r_{i\rightarrow j,I}= r^{p_i}_{I}+o(1)$ with any $I\neq \emptyset,\{1\},\ldots,\{L\}$. Therefore, to implement our message-passing algorithm, we assume $r_{i\rightarrow j,I}= r^{p_i}_{I}$ exactly for any $I\neq \emptyset,\{1\},\ldots,\{L\}$ and only iteratively update $r_{i\rightarrow j,l}$ by
\begin{equation}\label{eq:LBP prototype}
    r_{i\rightarrow j,l}^{t+1}=r^{p_i}_{l}+\sum_{k\neq i,j}\sqrt{\frac{\lambda^{(l)}}{n}}A^{(l)}_{i,k}\left[\sum_{x_k,y_k}\nu_{k\rightarrow i}^t(x_k,y_k)x_k^{(l)}\right],\forall i\neq j, l\in[L],
\end{equation}
where $\nu_{k\rightarrow i}^t$ is determined by $r^t_{k\rightarrow i}$.
After omitting the $o(1)$ term, every cavity field is given by 
\begin{equation*}
    \nu_{k\rightarrow i}^t(x_k,y_k)=\frac{1}{Z_{k\rightarrow i}^t}\exp\left[\sum_{l\in[L]} r_{k\rightarrow i,l}^tx_k^{(l)}+\sum_{|I|>1}r^{p_k}_I\mathsf{R}_I(x_k,y_k)\right],
\end{equation*}
so $\sum_{x_k,y_k}\nu_{k\rightarrow i}^t(x_k,y_k)x_k^{(l)}$ only depends on $\{r_{i\rightarrow j,l}^t:i\neq j\}$ across iterations. This fact motivates us to record it as a non-linear denoiser defined by
\begin{equation*}
    \mathcal{E}^{(l)}(r_{1:L};p)=\frac{\sum_{x,y}x^{(l)}\exp\left[\sum_{l^\prime\in[L]} r_{k\rightarrow i,l^\prime}^tx^{(l^\prime)}+\sum_{|I|>1}r^{p}_I\mathsf{R}_I(x,y)\right]}{\sum_{x,y}\exp\left[\sum_{l^\prime\in[L]} r_{k\rightarrow i,l^\prime}^tx^{(l^\prime)}+\sum_{|I|>1}r^{p}_I\mathsf{R}_I(x,y)\right]}
\end{equation*}
for any $r\in\R^L$ and distribution $p$ on $\{\pm 1\}^{L+L_1}$. Thereafter, recursions \eqref{eq:LBP prototype} are simplified to
\begin{equation*}
    r_{i\rightarrow j,l}^{t+1}=r^{p_i}_{l}+\sum_{k\neq i,j}\sqrt{\frac{\lambda^{(l)}}{n}}A^{(l)}_{i,k}\mathcal{E}^{(l)}(r_{k\rightarrow i}^t;p_k),\forall i\neq j, l\in[L].
\end{equation*}
For the three representative models we start from, $\mathcal{E}$ admits much more concise formalization.
\begin{itemize}
    \item In the Example~\ref{eg:multilayer SBM} of inhomogeneous multilayer SBM, there is no node-wise side information, so every node admits the same prior $p_k=p_\ML$ and we suppress the dependency of $\mathcal{E}^{(l)}(r;p)$ on prior distribution $p$. Lemma~\ref{lemma:ML denoiser} presents a closed form for $\mathcal{E}(r)$
    \begin{equation*}
        \mathcal{E}_{\ML,\rho}^{(l)}(r_{1:L})=\frac{\tanh\left(r_l+\bar{\rho}\right)\prod_{l^\prime}\cosh\left(r_{l^\prime}+\bar{\rho}\right)+\tanh\left(r_l-\bar{\rho}\right)\prod_{l^\prime}\cosh\left(r_{l^\prime}-\bar{\rho}\right)}{\prod_{l^\prime}\cosh\left(r_{l^\prime}+\bar{\rho}\right)+\prod_{l^\prime}\cosh\left(r_{l^\prime}-\bar{\rho}\right)}.
    \end{equation*}
    \item In the Example~\ref{eg:dynamic SBM} of dynamical SBM, there is no node-wise side information as well, so every node admits the same prior $p_k=p_\Dyn$ and we suppress the dependency of $\mathcal{E}^{(l)}(r;p)$ on prior distribution $p$. Lemma~\ref{lemma:dyn denoiser} provides an efficient Algorithm~\ref{alg:denoiser dyn} for $\mathcal{E}_{\Dyn,\rho}^{(l)}(r_{1:L})$ which runs in $O(L)$ time.
    \item In the Example~\ref{eg:semi sbm} of a graph model with partially observed labels, each $p_k$ is specified to $Z_k$ so $\mathcal{E}^{(l)}$ has to depend on $Z_k$. Equation~\eqref{eq:semi denoiser} provides a closed form,
    \begin{equation*}
        \mathcal{E}_{\Semi}(r,z)=z\mathbbm{1}\{z\neq\ast\}+\frac{\tanh(r)+\delta}{1+\delta\tanh(r)}\mathbbm{1}\{z=\ast\}.
    \end{equation*}
\end{itemize}

\paragraph{Approximate Message Passing.} Now we complete the derivation of coupled AMP algorithms by replacing the non-backtracking nature with an Onsager term, following \cite{bayati2011dynamics}. Write
\begin{equation*}
    r_{i\rightarrow j,l}^{t+1}=r^{p_i}_{l}+\sum_{k\neq i,j}\sqrt{\frac{\lambda^{(l)}}{n}}A^{(l)}_{i,k}\mathcal{E}^{(l)}(r^t_{k\rightarrow i})=r_{i\rightarrow j,l}^{t+1}+\delta_{i\rightarrow j,l}^{t+1},
\end{equation*}
where $\delta_{i\rightarrow j,l}^{t+1}=\sqrt{\frac{\lambda^{(l)}}{n}}A^{(l)}_{i,j}\mathcal{E}^{(l)}(r^t_{j\rightarrow i})$. Now use $r^t_{k\rightarrow i}=r^t_{k}-\delta^t_{k\rightarrow i}$ and expand $\mathcal{E}$ to its first derivative, to find
\begin{align*}
    r_{i\rightarrow j,l}^{t+1}=&r^{p_i}_{l}+\sum_{k\neq i}\sqrt{\frac{\lambda^{(l)}}{n}}A^{(l)}_{i,k}\mathcal{E}^{(l)}(r^t_{k\rightarrow i})\\
    \approx&r^{p_i}_{l}+\sum_{k\neq i}\sqrt{\frac{\lambda^{(l)}}{n}}A^{(l)}_{i,k}\mathcal{E}^{(l)}(r^t_{k})-\sqrt{\frac{\lambda^{(l)}}{n}}A^{(l)}_{i,k}\partial_l\mathcal{E}^{(l)}(r^t_{k})\delta^t_{k\rightarrow i,l}\\
    =&r^{p_i}_{l}+\sum_{k\neq i}\sqrt{\frac{\lambda^{(l)}}{n}}A^{(l)}_{i,k}\mathcal{E}^{(l)}(r^t_{k})-\frac{\lambda^{(l)}}{n}A^{(l)2}_{i,k}\partial_l\mathcal{E}^{(l)}(r^t_{k})\mathcal{E}^{(l)}(r^{t-1}_{i\rightarrow k}).
\end{align*}
Now we approximate $r^{t-1}_{i\rightarrow k}\approx r^{t-1}_{i}$ to get
\begin{equation*}
    r_{i,l}^{t+1}\approx r^{p_i}_{l}+\sum_{k\neq i}\sqrt{\frac{\lambda^{(l)}}{n}}A^{(l)}_{i,k}\mathcal{E}^{(l)}(r^t_{k})-\lambda^{(l)}\left[\sum_{k\neq i}\frac{A^{(l)2}_{i,k}}{n}\partial_l\mathcal{E}^{(l)}(r^t_{k})\right]\mathcal{E}^{(l)}(r^{t-1}_{i}).
\end{equation*}
We assume that $A^{(l)2}_{i,k}$ well concentrates around $1$ in the previous equation. Then the subtracted term becomes the empirical Onsager term
\begin{equation*}
    \sum_{k\neq i}\frac{A^{(l)2}_{i,k}}{n}\partial_l\mathcal{E}^{(l)}(r^t_{k})\approx \frac{1}{n}\sum_{k\neq i}\partial_l\mathcal{E}^{(l)}(r^t_{k}).
\end{equation*} 
To generalize our iterative scheme, we allow the use of  time dependent nonlinear mappings $\mathcal{E}_t$. We put these elements together to conclude
\begin{align}
    r_{i,l}^{t+1}=& r^{p_i}_{l}+\sum_{k\neq i}\sqrt{\frac{\lambda^{(l)}}{n}}A^{(l)}_{i,k}\mathcal{E}^{(l)}_t(r^t_{k})-\lambda^{(l)}\mathsf{d}_t^{(l)}\mathcal{E}^{(l)}_{t-1}(r^{t-1}_{i}),\\
    \mathsf{d}_t^{(l)}=&\frac{1}{n}\sum_{k\neq i}\partial_l\mathcal{E}^{(l)}_t(r^t_{k}).
\end{align}
Finally, we choose $m^{t+1}_{i,l}:=r_{i,l}^{t+1}-r^{p_i}_{l}$ so that the algorithm only depends on $p_i$ through the non-linearity $\mathcal{E}_t$. Moreover, diagonal terms $A^{(l)}_{i,i}$ are independent of the signal under detection, and their effect is negligible compared to off-diagonal terms. This completes the heuristic derivation of Algorithm~\ref{alg:coupled AMP}.

\subsection{State Evolution}
\label{sec:state_evolution}
The core strategy in establishing state evolution is to reduce  Algorithm~\ref{alg:coupled AMP} to an AMP algorithm with a single sensing matrix and non-separable denoising functions.

\paragraph{Reduction to Single Symmetric AMP.} Let us consider an abstract class of AMP recursions, and make a remark about \cite{gerbelot2021graph} before moving on.
\begin{align}
    \bM^{t+1} &= \bar{\bA}\hbX^{t}-\hbX^{t-1}(\bb^{t})^{\top} && \in \R^{nL\times L} \, , \label{eq:sym-amp-iteration-1} \\
    \hbX^{t} &=f^{t}(\bM^{t}) && \in \R^{nL\times L} \, , \\ 
    \bb_t &= \frac{1}{nL} \sum_{i=1}^{nL} \frac{\partial f^t_i}{\partial \bM_i}(\bM^t) && \in 
    \R^{L\times L}\, . \label{eq:sym-amp-iteration-2}
\end{align}
 
\begin{Remark}\label{remark:comparison to gerbelot-berthier}
    As an intermediate step, \cite{gerbelot2021graph} adapt the results in \cite{berthier2020state} and provide Lemma 13 to address the state evolution of recursions \eqref{eq:sym-amp-iteration-1}-\eqref{eq:sym-amp-iteration-2}, which is restated as Lemma~\ref{thm:symmetric-biased} in our manuscript later. Both Algorithm~\ref{alg:coupled AMP} and their algorithm can be embedded onto the abstract class of recursions \eqref{eq:sym-amp-iteration-1}-\eqref{eq:sym-amp-iteration-2}, so their intermediate result (Lemma~\ref{thm:symmetric-biased}) greatly simplifies our derivation of state evolution, Theorem~\ref{thm:coupled-AMP-SE}. But it is noteworthy that Algorithm~\ref{alg:coupled AMP} is not included in the framework proposed by \cite{gerbelot2021graph}, since we only get to observe symmetric sensing matrices and prior knowledge connect all these observations.
\end{Remark}

For our setting, we define
\begin{equation*}
    \bar{\bA} = 
    \begin{pmatrix}
    \sqrt{\frac{1}{nL}}\bA^{(1)} & &\ast\\
    & \ddots & \\
    \ast & & \sqrt{\frac{1}{nL}}\bA^{(L)}
    \end{pmatrix}\in\R^{nL\times nL},
\end{equation*}
where $\ast$ denotes additional independent standard normal variables. Later, we will realize $\bar{\bA}$ as a rescaled GOE plus a low rank signal component. Initialization $\bM^0$ is then given by
\begin{equation}\label{eq:reduction-sym-init}
    \bM^0 = \begin{pmatrix}
        m^0_{\cdot,1} & & & \star\\
        & m^0_{\cdot,2} & &\\
        & & \ddots & \\
        \star & & & m^0_{\cdot,L}
    \end{pmatrix} \in \R^{nL \times L},\quad m^0_{\cdot,l}=\begin{pmatrix}
        m^0_{1,l}\\
        m^0_{2,l}\\
        \cdots\\
        m^0_{n,l}\\
    \end{pmatrix}\in\R^{n\times 1},
\end{equation}
where $\star$ indicates entries whose values do not influence the output of the algorithm. Finally, we absorb different SNRs across layers into the definitions of non-linear denoisers,
\begin{equation}\label{eq:reduction-sym-non-linear}
    f^{t}\begin{pmatrix}
        m_{1,1} & & \\
        m_{2,1} & & \\
        \cdots & & \star \\
        m_{n,1} & & \\
        & \ddots & \\
        & & m_{1,L} \\
        & & m_{2,L} \\
        \star & & \cdots  \\
        & & m_{n,L} \\
    \end{pmatrix}=\begin{pmatrix}
        \sqrt{\lambda^{(1)}L}\mathcal{E}_t^{(1)}(m_1,Z_1) & & \\
        \sqrt{\lambda^{(1)}L}\mathcal{E}_t^{(1)}(m_2,Z_2) & & \\
        \cdots & & 0 \\
        \sqrt{\lambda^{(1)}L}\mathcal{E}_t^{(1)}(m_n,Z_n) & & \\
        & \ddots & \\
        & & \sqrt{\lambda^{(L)}L}\mathcal{E}_t^{(L)}(m_1,Z_1) \\
        & & \sqrt{\lambda^{(L)}L}\mathcal{E}_t^{(L)}(m_2,Z_2) \\
        0 & & \cdots  \\
        & & \sqrt{\lambda^{(L)}L}\mathcal{E}_t^{(L)}(m_n,Z_n) \\
    \end{pmatrix} \in \R^{nL \times L},
\end{equation}
where $m_i=(m_{i,1},\ldots,m_{i,L})$ indicates all variables corresponding to node $i$ across all layers. The following lemma formalizes this reduction. 

\begin{Lemma}\label{lemma:reduction to sym AMP}
    With initialization and non-linearities specified as above,
    iterations \eqref{eq:sym-amp-iteration-1}-\eqref{eq:sym-amp-iteration-2} produce the same sequence as our original multilayer AMP Algorithm~\ref{alg:coupled AMP}, i.e.
    \begin{equation}\label{eq:reduction-sym-X}
        \bM^t = \begin{pmatrix}
            m^t_{\cdot,1} & & & \star\\
            & m^t_{\cdot,2} & &\\
            & & \ddots & \\
            \star & & & m^t_{\cdot,L}
        \end{pmatrix} \in \R^{nL \times L},\quad m^t_{\cdot,l}=\begin{pmatrix}
            m^t_{1,l}\\
            m^t_{2,l}\\
            \cdots\\
            m^t_{n,l}\\
        \end{pmatrix}\in\R^{n\times 1}.
    \end{equation}
\end{Lemma}

\begin{proof}
    This reduction can be verified via induction. The conclusion holds automatically for $t=0$. Suppose the reduction holds until iterate $t\in\mathbb{N}$.
    
    We need to verify that the Onsager terms arising from the two algorithms (\eqref{eq:coupled AMP Onsager} and \eqref{eq:sym-amp-iteration-2}) are, in fact, the same. By definition $\bb_t = \frac{1}{nL} \sum_{i=1}^{nL} \frac{\partial f^t_i}{\partial \bM_i}(\bM^t)  \in 
    \R^{L\times L}$, each $f^t_i$ and $\bM_i$ are row vectors of dimension $1\times L$. Therefore, for $l_1,l_2\in[L]$,
    \begin{align*}
        \bb_t[l_1,l_2]=\frac{1}{nL}\sum_{i=1}^n\sum_{l^\prime=1}^L \frac{\partial f^t_{i-n+nl^\prime}[l_1]}{\partial \bM_{i-n+nl^\prime}[l_2]}(\bM^t).
    \end{align*}
    If $l^\prime =l_1$, we have $f^t_{i-n+nl^\prime}[l_1](\bM^t)=\sqrt{\lambda^{(l_1)}L}\mathcal{E}_t^{(l_1)}(m_i^t,Z_i)=\sqrt{\lambda^{(l_1)}L}\mathcal{E}_t^{(l_1)}(m_{i,1}^t,\ldots,m_{i,L}^t,Z_i)$; otherwise $f^t_{i-n+nl^\prime}[l_1]=0$. As for the second coordinate, if $l^\prime=l_2$, we have $\partial \bM_{i-n+nl^\prime}^t[l_2]=m_{i,l^\prime}^t$; otherwise, $\partial \bM_{i-n+nl^\prime}^t[l_2]$ would be denoted as $\star$, implying that it wouldn't be used in later iterations and thus irrelevant. As a result, for any $l_1\neq l_2$, we directly have $\bb_t[l_1,l_2]=0$. While for diagonal terms, it holds
    \begin{equation*}
        \bb_t[l,l]=\frac{1}{nL}\sum_{i=1}^n\sum_{l=1}^L \frac{\partial f^t_{i-n+nl}[l]}{\partial \bM_{i-n+nl}[l]}(\bM^t)=\frac{\sqrt{\lambda^{(l)}}}{n\sqrt{L}}\sum_{i=1}^n\partial_l\mathcal{E}^{(l)}_t(m^t_{i},Z_i)= \sqrt{\frac{\lambdal}{L}}\mathsf{d}_t^{(l)}
    \end{equation*}
    matching exactly with \eqref{eq:coupled AMP Onsager}.

    Now we will further employ the block structure in \eqref{eq:sym-amp-iteration-1}-\eqref{eq:sym-amp-iteration-2},
    \begin{equation*}
        \bar{\bA}\hbX^{t}=\bar{\bA}f^{t}(\bM^{t})=\begin{pmatrix}
            \sqrt{\frac{\lambda^{(1)}}{n}}\bA^{(1)}\mathcal{E}_t^{(1)}(\mathbf{m}^t,\bZ)  & & \star\\
            &  \ddots & \\
            \star & &  \sqrt{\frac{\lambdal}{n}}\bA^{(L)}\mathcal{E}_t^{(L)}(\mathbf{m}^t,\bZ)
        \end{pmatrix}
    \end{equation*}
    and
    \begin{equation*}
        \hbX^{t-1}(\bb^{t})^{\top}=\begin{pmatrix}
            \lambda^{(1)}\mathsf{d}_t^{(1)}\mathcal{E}^{(1)}_{t-1}(\mathbf{m}^{t-1},\bZ) & &  \star\\
            & \ddots & \\
            \star & & \lambda^{(L)}\mathsf{d}_t^{(L)}\mathcal{E}^{(L)}_{t-1}(\mathbf{m}^{t-1},\bZ)
        \end{pmatrix}.
    \end{equation*}
    By subtracting these two matrices, we find \eqref{eq:reduction-sym-X} to hold for $t+1$.
\end{proof}

\paragraph{SE with Non-separable Denoisers.} 
Arrange the signal matrix to be detected by
\begin{equation}\label{eq:symmetric-low-rank-component}
    \bV_0 = \begin{pmatrix}
        \sqrt[4]{L\lambda^{(1)}}\bX^{(1)} & & & 0\\
        & \sqrt[4]{L\lambda^{(2)}}\bX^{(2)} & &\\
        & & \ddots & \\
        0 & & & \sqrt[4]{L\lambda^{(L)}}\bX^{(L)}
    \end{pmatrix} \in \R^{nL \times L},\quad \bX^{(l)}=\begin{pmatrix}
        X^{(l)}_1\\
        X^{(l)}_2\\
        \cdots\\
        X^{(l)}_n\\
    \end{pmatrix}\in\R^{n\times 1}.
\end{equation}
Decompose the sensing matrices into $\bA=\frac{1}{nL}\bV_0\bV_0^{\top}+\bar{\bW}$, where $\bar{\bW}$ is defined by
\begin{equation*}
    \bar{\bW} = 
    \begin{pmatrix}
    \sqrt{\frac{1}{nL}}\bW^{(1)} & &\ast\\
    & \ddots & \\
    \ast & & \sqrt{\frac{1}{nL}}\bW^{(L)}
    \end{pmatrix}\in\R^{nL\times nL},
\end{equation*}
with $\bW^{(1)},\ldots,\bW^{(L)}$ introduced in \eqref{eq:spiked matrix each layer} as GOEs. Therefore, $\bar{\bW}$ is itself a rescaled GOE of dimension $nL$.

\begin{Definition}[state evolution iterates]\label{def:se-sym}
    For a biased sensing matrix, the state evolution iterates are composed of two components: one infinite-dimensional array $\{\bar{\boldsymbol{\mu}}^t\in\R^{L\times L}\}_{t>0}$ denoting the \textbf{bias coefficient}; and an infinite-dimensional two-way array $\{\bar{\boldsymbol{\kappa}}^{s,r}\in\R^{L\times L}\}_{r,s>0}$ denoting \textbf{covariances}. These two arrays are generated as follows. 
    Define the first state evolution iterate
    \begin{align}
        \bar{\boldsymbol{\mu}}^{1} =& \lim_{n \to \infty} \frac{1}{nL}  \mathbb{E}\left[\bV_0^\top f^{0}(\bM^{0})\right]\in\R^{L\times L},\label{eq:sym-bias-se-mu-init}\\
        \bar{\boldsymbol{\kappa}}^{1,1} =& \lim_{n \to \infty} \frac{1}{nL} \mathbb{E}\left[ f^{0}(\bM^0)^\top f^{0}(\bM^0)\right]\in\R^{L\times L}.\label{eq:sym-bias-se-kappa-init}
    \end{align}
    Recursively, once $\{\bar{\boldsymbol{\kappa}}^{s,r}\}_{s,r \leq t}$ and $\{\bar{\boldsymbol{\mu}}^{s}\}_{s\leq t}$ are defined for some $t \geq 1$, take $(\bN^1,\dots, \bN^t) \in (\R^{nL \times L})^t$ a centered Gaussian vector of covariance $(\bar{\boldsymbol{\kappa}}^{s,r})_{s,r\leqslant t}\otimes \mathbf{I}_{nL}$ and $\bN^0 = \bM^0$. We then define new state evolution iterates by
    \begin{align}
        \bar{\boldsymbol{\mu}}^{t+1} = &\lim_{n \to \infty} \frac{1}{nL} \E\left[ \bV_0^\top f^{t}(\bV_0\bar{\boldsymbol{\mu}}^t+\bN^{t}) \right]  ,\label{eq:sym-bias-se-mu}\\
        \bar{\boldsymbol{\kappa}}^{t+1, s+1} = \bar{\boldsymbol{\kappa}}^{s+1, t+1} = &\lim_{n \to \infty} \frac{1}{nL} \E\left[ f^{s}(\bV_0\bar{\boldsymbol{\mu}}^s+\bN^{s})^\top f^{t}(\bV_0\bar{\boldsymbol{\mu}}^t+\bN^{t}) \right] \, , \qquad s \in \{ 0, \dots, t \} \, .\label{eq:sym-bias-se-kappa}
    \end{align}
\end{Definition}

Given all these definitions, \cite[Lemma 13]{gerbelot2021graph} can be stated as below.


\begin{Lemma}\label{thm:symmetric-biased}
    Assume the same conditions as \cite[Lemma 13]{gerbelot2021graph} which mainly address the regularity of adopted denoisers. Define, as above,  $(\bN^1,\dots, \bN^t) \in (\R^{nL \times L})^t$ a centered Gaussian vector of covariance $(\bar{\boldsymbol{\kappa}}^{s,r})_{s,r\leqslant t}\otimes \mathbf{I}_{nL}$ and $\bN^0 = \bM^0$. Then for any sequence $\phi_n : (\R^{nL \times L})^{t+1} \to \R$ of pseudo-Lipschitz functions, 
    \begin{equation*}
        \phi_{n}\left(\bM^{0}, \bM^{1}, ..., \bM^{t}\right) - \mathbb{E}\left[\phi_{n}\left(\bN^{0},\mathbf{V}_{0}\bar{\boldsymbol{\mu}}^{1}+\bN^{1},...,\mathbf{V}_{0}\bar{\boldsymbol{\mu}}^{t}+\bN^{t}\right)\right]\overset{\text{P}}{\rightarrow}0,
    \end{equation*}
\end{Lemma}

\begin{proof}[Proof Outline.]
This theorem starts from a symmetric AMP algorithm with a mean-zero sensing matrix $\bar{\bW}$. Specifically, define a bias-corrected sequence by
\begin{align}
    \mathbf{S}^{t+1} &= \bar{\bW}\mathbf{T}^t-\mathbf{T}^{t-1}(\widetilde{\bb}^{t})^{\top} && \in \R^{nL\times L} \, , \label{eq:sym-amp-iteration-unbiased-1} \\
    \mathbf{T}^t &=\hat{f}^{t}(\bV_0\bar{\boldsymbol{\mu}}^t+\mathbf{S}^{t}) && \in \R^{nL\times L} \, , \\ 
    \widetilde{\bb}_t &= \frac{1}{nL} \sum_{i=1}^{nL} \frac{\partial f^t_i}{\partial \bM_i}(\bV_0\bar{\boldsymbol{\mu}}^t+\mathbf{S}^{t}) && \in 
    \R^{L\times L}\, . \label{eq:sym-amp-iteration-unbiased-2}
\end{align}
Subsequently, \cite[Theorem 2]{gerbelot2021graph} would yield that
\begin{equation*}
    \phi_n\left(\mathbf{S}^0, \mathbf{S}^1, \dots, \mathbf{S}^t\right) - \E\left[ \phi_n\left(\bN^0, \bN^1, \dots, \bN^t\right) \right]\overset{\text{P}}{\rightarrow}0.
\end{equation*} 
Lastly, following the techniques developed in \cite[Section D.1]{gerbelot2021graph}, \cite[Section B.4]{deshpande2017asymptotic} and \cite[Section 6.8]{feng2022unifying}, we should choose a specific $\hat{f}^t$ in \eqref{eq:sym-amp-iteration-unbiased-1}-\eqref{eq:sym-amp-iteration-unbiased-2} based on those denoisers $f^t$ used in \eqref{eq:sym-amp-iteration-1}-\eqref{eq:sym-amp-iteration-2}. In this way, we are able to upper bound $\|\bV_0\bar{\boldsymbol{\mu}}^t+\mathbf{S}^t-\bM^t\|$ to conclude the whole theorem. 
\end{proof}

\paragraph{Back to Coupled AMP.} Primarily, we need to make sure that state evolution iterates are identical under two different settings.
\begin{Lemma}\label{lemma:connect SE}
    Comparing Definitions~\ref{def:SE iterates} and \ref{def:se-sym}, it holds that $\bar{\boldsymbol{\mu}}^t_{l,l}=\frac{\lambda^{(l)\frac{3}{4}}}{L^{\frac{1}{4}}} \mu^t_l$ and $\bar{\boldsymbol{\kappa}}^{t,t}_{l,l}=\lambdal\kappa^t_l$.
\end{Lemma}
\begin{proof}
First, since $f^t$ and $\bV_0$ are both constructed into a block diagonal form as \eqref{eq:reduction-sym-non-linear} and \eqref{eq:symmetric-low-rank-component} respectively, one would conclude $\bV_0^\top f^{t}(\bV_0\bar{\boldsymbol{\mu}}^t+\bN^{t})$ to be diagonal. It then follows from \eqref{eq:sym-bias-se-mu-init} and \eqref{eq:sym-bias-se-mu} that all matrices $\bar{\boldsymbol{\mu}}^t,t\ge 1$ are diagonal. Therefore, $\bV_0\bar{\boldsymbol{\mu}}^t$ can be written as
    \begin{equation*}
        \bV_0\bar{\boldsymbol{\mu}}^t=\diag\left(\bar{\boldsymbol{\mu}}^t_{1,1}\left(\lambda^{(1)}L\right)^{\frac{1}{4}}\bX^{(1)},\ldots,\bar{\boldsymbol{\mu}}^t_{L,L}\left(\lambda^{(L)}L\right)^{\frac{1}{4}}\bX^{(L)}\right).
    \end{equation*}

    Secondly, implied by the construction of $f^t$ in \eqref{eq:reduction-sym-non-linear}, later iterations in \eqref{eq:sym-bias-se-mu} and \eqref{eq:sym-bias-se-kappa} only depend on the joint asymptotic distribution of certain entries of $\bV_0\bar{\boldsymbol{\mu}}^t+\bN^{t}$, not all. Specifically, for each row (indexed by $i-n+nl$ with $i\in[n],l\in[L]$) in $\bV_0\bar{\boldsymbol{\mu}}^t+\bN^{t}$, only the entry in $l$-th column really matters and its distribution is $\mathcal{N}\left(\bar{\boldsymbol{\mu}}^t_{l,l}\left(\lambda^{(l)}L\right)^{\frac{1}{4}}X^{(l)}_i,\bar{\boldsymbol{\kappa}}^{t,t}_{l,l}\right)$. In conclusion, $(\diag(\bar{\boldsymbol{\kappa}}^{s,r}))_{s,r>0}$ suffices to describe the selected entries.
    
    Lastly, since we do not care about the correlation between iterations, it suffices to only study $(\diag(\bar{\boldsymbol{\kappa}}^{t,t}))_{t>0}$. The whole state evolution iterates boil down to only diagonal entries:
    \begin{align}
        \bar{\boldsymbol{\mu}}^{t+1}_{l,l} = &\lim_{n \to \infty} \frac{\lambda^{(l)\frac{3}{4}}}{nL^{\frac{1}{4}}} \E\left[ \sum_{i=1}^n X^{(l)}_i \mathcal{E}_t^{(l)}\left(\widetilde{m}^t_i,Z_i\right)\right],\\
        \bar{\boldsymbol{\kappa}}^{t+1, t+1}_{l,l} = &\lim_{n \to \infty} \frac{\lambda^{(l)}}{n} \E\left[ \sum_{i=1}^n\mathcal{E}_t^{(l)}\left(\widetilde{m}^t_i,Z_i\right)^\top \mathcal{E}_t^{(l)}\left(\widetilde{m}^t_i,Z_i\right) \right].
    \end{align}
    where $\widetilde{m}^t_1,\ldots,\widetilde{m}^t_n$ are independent duplications of $\widetilde{m}^t$ from Definition~\ref{def:SE iterates}. Since $(\widetilde{m}^t_i,Z_i)$ are all i.i.d., we can get rid of the summation over $i\in[n]$ and conclude this lemma.
\end{proof}
\begin{proof}[Proof of Theorem~\ref{thm:coupled-AMP-SE}]
    This theorem follows from using Lemmas~\ref{lemma:reduction to sym AMP} and \ref{lemma:connect SE}, and choosing specific separable test function $\phi_n$ in Lemma~\ref{thm:symmetric-biased}.
\end{proof}

\subsection{Algorithmic Universality}
\label{sec:alg_universality}
Recall the low-rank component extracted from adjacency matrices,
\begin{equation*}
    \frac{\mathbf{G}^{(l)}-d^{(l)}/n}{\sqrt{d^{(l)}(1-d^{(l)}/n)/n}}=\sqrt{\frac{\lambda^{(l)}}{n}}\bX^{(l)}\bX^{(l)\top}+\mathbf{H}^{(l)},
\end{equation*}
with $\lambda^{(l)}_n=\frac{n(a^{(l)}-b^{(l)})^2}{4d^{(l)}(1-d^{(l)}/n)}$ being effective SNR in layer $l$. Conditioned on each $\bX^{(l)}$, the noise matrix $\mathbf{H}^{(l)}$ has independent off-diagonal entries satisfying that for each $1\le i<j\le n$,
\begin{align*}
    &\quad\E\left[H_{i,j}^{(l)}\right]=0, \qquad\left|H_{i,j}^{(l)}\right| \leq \frac{1}{\sqrt{d^{(l)}(1-d^{(l)}/n)/n}},\\
    &\E\left[H_{i,j}^{(l)2}\right]\in\left\{\frac{a^{(l)}(1-a^{(l)}/n)}{d^{(l)}(1-d^{(l)}/n)},\frac{b^{(l)}(1-b^{(l)}/n)}{d^{(l)}(1-d^{(l)}/n)}\right\}.
\end{align*}
Under the assumption that $d^{(l)}(1-d^{(l)}/n)\rightarrow\infty$, effective SNR $\lambdal_n\rightarrow\lambdal$ and $S_{i,j}^{(l)}=\E\left[H_{i,j}^{(l)2}\right]\rightarrow1$. This implies every $\mathbf{H}^{(l)}/\sqrt{n}$ is a generalized Wigner matrix in the sense of Definition 2.3 of \cite{wang2022universality}. Moreover since $\min_l d^{(l)}(1-d^{(l)}/n)\ge C\log n/n$ for some $C>0$, we can show using Theorem
2.7 of \cite{benaych2020spectral} and (2.4) of \cite{wang2022universality} that $\|\mathbf{H}^{(l)}\|_{op}\le C\sqrt{n}$ almost surely for large $n,p$.
Hence, the formal proof of Proposition~\ref{prop:algorithmic universality} is a combination of the proof techniques of Theorem 2.4 in \cite{wang2022universality} and those of Theorem~\ref{thm:coupled-AMP-SE}.

\paragraph{Reformulation of AMP Iterates and Universality Class.} We begin with defining an auxiliary bias-corrected AMP iterates with mean-zero sensing matrices $\mathbf{H}$ given in \eqref{eq:graph low rank component}
\begin{align}
u_{k,l}^{t+1}&=\mathcal{U}^{(l)}_t(s^t_{k},X_k,Y_k,Z_k),\label{eq:coupled AMP update reformulated 1}\\
s_{i,l}^{t+1}&=\sum_{k}\sqrt{\frac{\lambda^{(l)}}{n}}H^{(l)}_{i,k}u_{k,l}^{t+1}-\lambda^{(l)}\mathsf{b}_t^{(l)}u_{i,l}^{t},\label{eq:coupled AMP update reformulated 2}
\end{align}
where we take the non-linear denoisers as
\begin{equation*}
\mathcal{U}^{(l)}_t(s^t,X,Y,Z)=\mathcal{E}^{(l)}_t(s^t+\boldsymbol{\lambda}\odot\boldsymbol{\mu}^t\odot X,Z).
\end{equation*}
Different from Algorithm~\ref{alg:coupled AMP}, we adopt a constant coefficient before the Onsager terms, which is given jointly with the state evolution iterates in the following recursive way
\begin{align}
\mathsf{b}_t^{(l)}&=\mathbb{E}\left[\partial_l\mathcal{U}^{(l)}_t(n^t,X,Y,Z)\right],\label{eq:coupled AMP Onsager reformulated}\\
\Sigma_{t+1}^{(l)}&=\mathbb{E}\left[\mathcal{U}^{(l)}_{t}(n^{t},X,Y,Z)^2\right],
\end{align}
where $n^t\sim\mathcal{N}\left(0,\diag(\Sigma_{t}^{(1)},\ldots,\Sigma_{t}^{(L)} )\right)$ is an $L$-dimensional Gaussian vector independent of the law of $(X,Y,Z)$. This recursion starts from $\Sigma_{0}^{(l)}=\mathbb{E}\left[\left(S_l^0\right)^2\right]$, where $S^0$ is the limit of initialization $\mathbf{s}^0$.

When examining the Onsager coefficients at time $t$, the empirical mean in \eqref{eq:coupled AMP Onsager} concentrates near its population mean when conditioning on previous iterates $\mathbf{m}^{0:t}$ when $n\rightarrow\infty$. As suggested by Corollary 2 in \cite{berthier2020state} and Remark 2.9 in \cite{wang2022universality}, as long as the adopted Onsager coefficient is a consistent estimator of this constant version, the state evolution results would not change. So we adopt an empirical mean as Onsager coefficients in presenting Algorithm~\ref{alg:coupled AMP} as it is easier to implement, but turn to use its expectation for proving universality.

By slightly adapting techniques in Section B.4 of \cite{deshpande2017asymptotic}, we find this sequence of bias-corrected AMP iterates is close to the output of Algorithm~\ref{alg:coupled AMP} in the sense that $\bar{m}_l^t\approx \lambdal \mu_l^t X_l+s_l^t$. Therefore, it suffices to establish the universality of sequence \eqref{eq:coupled AMP update reformulated 1}-\eqref{eq:coupled AMP update reformulated 2} in which $(X,Y,Z)$ is treated as some side information. Some necessary conditions in this regard are collected below. 

\begin{Assumption}\label{assump:ufconvergence}
As $n \to \infty$, the empirical distribution of $(\mathbf{s}^0,\mathbf{X},\mathbf{Y},\mathbf{Z})$ converges in Wasserstein distance to some distribution $(S^0,X,Y,Z)$ with finite moments of all orders. This means, for any polynomially bounded continuous function $f$,
\begin{equation}
\frac{1}{n}\sum_{i=1}^n f(s^0_i,X_i,Y_i,Z_i) \to \E\left[f(S^0,X,Y,Z)\right].
\end{equation}
Furthermore, multivariate polynomials are dense in the real $L^2$-space of functions $f:\R^{L}\times\mathcal{V} \to \R$ with inner-product
\begin{equation*}
(f,g) \mapsto \E[f(S^0,X,Y,Z)g(S^0,X,Y,Z)].
\end{equation*}
\end{Assumption}

\paragraph{Tensor Network and its Universality.}
For every $t\in\mathbb{N}$, suppose $\widetilde{\mathcal{U}}_t:\mathbb{R}^{L}\times \mathcal{Z}\rightarrow\mathbb{R}^L$ is a reasonable polynomial approximation of $\mathcal{U}_t$. Then we replace all $\mathcal{U}_t$ in \eqref{eq:coupled AMP update reformulated 1}-\eqref{eq:coupled AMP update reformulated 2} with $\widetilde{\mathcal{U}}_t$ and output another approximating sequence below
\begin{align}
\tilde{u}_{k,l}^{t+1}=\widetilde{\mathcal{U}}^{(l)}_t(\tilde{s}^t_{k},X_k,Y_k,Z_k), \qquad \tilde{s}_{i,l}^{t+1}=\sum_{k}\sqrt{\frac{\lambda^{(l)}}{n}}H^{(l)}_{i,k}\tilde{u}_{k,l}^{t+1}-\lambda^{(l)}\widetilde{\mathsf{b}}_t^{(l)}\tilde{u}_{i,l}^{t}.\label{eq:coupled AMP update polyn}
\end{align}
For sequence $\{\tilde{\mathbf{u}}^1,\tilde{\mathbf{s}}^1,\ldots,\tilde{\mathbf{u}}^t,\tilde{\mathbf{s}}^t\}$, the Onsager coefficients and state evolution iterates are respectively given as
\begin{equation}
\widetilde{\Sigma}_{t+1}^{(l)}=\mathbb{E}\left[\widetilde{\mathcal{U}}^{(l)}_{t}(n^{t},X,Y,Z)^2\right],\qquad \tilde{\mathsf{b}}_t^{(l)}=\mathbb{E}\left[\partial_l\widetilde{\mathcal{U}}^{(l)}_t(n^t,X,Y,Z)\right],\label{eq:coupled AMP Onsager polyn}
\end{equation}
in which $n^t\sim\mathcal{N}\left(0,\diag(\widetilde{\Sigma}_{t}^{(1:L)})\right)$. In the following, we introduce the notion of tensor networks.

\begin{Definition}
A \emph{diagonal tensor network} $T=(\mathcal{V},\mathcal{E},\mathcal{L},\{q_v\}_{v \in \mathcal{V}})$ in $k$ variables is an undirected tree graph with vertices $\mathcal{V}$ and edges $\mathcal{E} \subset \mathcal{V} \times \mathcal{V}$ labeled by $\mathcal{L}[u,v]\in[L]$ for any $(u,v)\in\mathcal{E}$. Moreover, each vertex $v \in \mathcal{V}$ is also labeled by a polynomial function $q_v:\R^k \to \R$. The \emph{value} of $T$ on $L$ different symmetric matrices $\left\{\mathbf{H}^{(l)} \in \mathbb{R}^{n \times n}:l\in[L]\right\}$ and vectors $\mathbf{x}_1,\ldots,\mathbf{x}_k \in \mathbb{R}^n$ is
\begin{equation*}
\val_T(\mathbf{H};\mathbf{x}_1,\ldots,\mathbf{x}_k)=\frac{1}{n} \sum_{\mathbf{i} \in [n]^\mathcal{V}}q_{\mathbf{i}|T} \cdot  H_{\mathbf{i}|T}
\end{equation*}
where, for each index tuple $\mathbf{i}=(i_v:v \in \mathcal{V}) \in [n]^\mathcal{V}$, we set
\begin{equation*}
    q_{\mathbf{i}|T}=\prod_{v \in \mathcal{V}} q_v(x_1[i_v],\ldots,x_k[i_v]), \qquad H_{\mathbf{i}|T}=\prod_{(u,v) \in \mathcal{E}} H_{i_u,i_v}^{(\mathcal{L}[u,v])}=\prod_{l}\prod_{(u,v) \in \mathcal{E}:\mathcal{L}[u,v]=l} H_{i_u,i_v}^{(l)}.
\end{equation*}
\end{Definition}

\begin{Remark}
Compared to the original notion in Definition 2.11 in \cite{wang2022universality}, our definition extends to a case where we have multiple sensing matrices.
Tensor networks can be understood as iteratively contracting all tensor-matrix-tensor products represented by the edges of the tree.
\end{Remark}

\begin{Lemma}\label{lemma:polynomial_tensor_network_decomposition}
Fix any $t \geq 1$. Let $(\tilde{\mathbf{s}}^0, \tilde{\mathbf{u}}^1,\tilde{\mathbf{s}}^1, \ldots,\tilde{\mathbf{s}}^t, \tilde{\mathbf{u}}^t)$ be the iterates given in \eqref{eq:coupled AMP update polyn}, where Onsager coefficients $\widetilde{\mathsf{b}}_{1:t}^{(1:L)}$ are scalar constants given in \eqref{eq:coupled AMP Onsager polyn}.
Suppose all non-linear denoisers $\widetilde{\mathcal{U}}_{t+1}$ are polynomial functions applied row-wise.
Then for any polynomial test function $p: \mathbb{R}^{(2t+1)L}\times\{\pm 1\}^{L+L_0}\times\mathcal{Z} \to \mathbb{R}$, there holds
\begin{align*}
\frac{1}{n}\sum_{i=1}^n p(\tilde{s}_i^0, \ldots, \tilde{s}_i^t,\tilde{u}_i^1, \ldots, \tilde{u}_i^t,X_i,Y_i,Z_i)= \sum_{T \in \mathcal{F}} \val_T(\mathbf{H}/\sqrt{n}; \tilde{\mathbf{s}}^0, \mathbf{X},\mathbf{Y},\mathbf{Z}). 
\end{align*}
for some finite set $\mathcal{F}$ of
diagonal tensor networks in $k+1$ variables.
\end{Lemma}

\begin{proof}
Begin with exploiting $p$ being a polynomial, 
\begin{align*}
\frac{1}{n}\sum_{i=1}^n p(\tilde{s}_i^{0:t},\tilde{u}_i^{1:t},Z_i)=\val_T\left(\mathbf{H}/\sqrt{n};\tilde{\mathbf{s}}^{0:t},\tilde{\mathbf{u}}^{1:t},\mathbf{Z}\right).
\end{align*}
where $T$ is a tensor network with only one vertex $v$ whose associated polynomial is $q_v = p$. We claim that given any tensor network $T=(\mathcal{V},\mathcal{E},\mathcal{L},\{q_v\}_{v \in \mathcal{V}})$ in the variables $(\tilde{s}^{0:t},\tilde{u}^{1:t},X,Y,Z)$, we can decompose
\begin{equation}\label{eq:sym_reduction_1}
\val_T\left(\mathbf{H}/\sqrt{n};\tilde{\mathbf{s}}^{0:t},\tilde{\mathbf{u}}^{1:t},\mathbf{X},\mathbf{Y},\mathbf{Z}\right) = \sum_{T^\prime \in \mathcal{F}} \val_{T^\prime}\left(\mathbf{H}/\sqrt{n};\tilde{\mathbf{s}}^{0:(t-1)},\tilde{\mathbf{u}}^{1:t},\mathbf{X},\mathbf{Y},\mathbf{Z}\right)
\end{equation}
where $\mathcal{F}$ is a finite set of tensor networks in the variables $(\tilde{s}^{0:(t-1)}, \tilde{u}^{1:t},X,Y,Z)$. To show this, recall that
\begin{align*}
\val_T\left(\mathbf{H}/\sqrt{n};\tilde{\mathbf{s}}^{0:t},\tilde{\mathbf{u}}^{1:t},\mathbf{X},\mathbf{Y},\mathbf{Z}\right)&= \frac{1}{n}\sum_{\mathbf{i} \in [n]^\mathcal{V}} \prod_{v \in \mathcal{V}} q_v(\tilde{s}^{0:t}_{i_v}, \tilde{u}^{1:t}_{i_v}, Z_{i_v}) \cdot  H_{\mathbf{i}|T}.
\end{align*}
Subsequently expand each $q_v$ in terms of a polynomial in $\tilde{\mathbf{s}}^t$, where the dependence on 
$(\tilde{\mathbf{s}}^{0:(t-1)}, \tilde{\mathbf{u}}^{1:t}, \mathbf{X},\mathbf{Y},\mathbf{Z})$ is absorbed in coefficients. Then plug in  $\tilde{\mathbf{s}}^t_{\cdot,l} = \sqrt{\lambdal/n}\mathbf{H}^{(l)} \tilde{\mathbf{u}}^t_{\cdot,l}-\lambdal\widetilde{\mathsf{b}}_{t-1}^{(l)} \tilde{\mathbf{u}}_{\cdot,l}^{t-1}$ to rearrange the previous expansion in terms of every $\mathbf{H}^{(l)}\tilde{\mathbf{u}}_{\cdot,l}^t$,
\begin{align*}
q_v(\tilde{s}^{1:t}_{i_v}, \tilde{u}^{1:t}_{i_v}, Z_{i_v}) = \sum_{\theta\in\Theta_v} q_{v,\theta}(\tilde{s}^{0:(t-1)}_{i_v}, \tilde{u}^{1:t}_{i_v}, Z_{i_v}) \cdot \prod_{l=1}^L\bigg(\sum_{j=1}^n H_{i_v, j}^{(l)} \tilde{u}^t_j /\sqrt{n}\bigg)^{\theta_l}
\end{align*}
where $\Theta_v\subset\mathbb{N}^L$ is a finite set of index tuples corresponding to $q_v$ in $\widetilde m^t$, and $\{q_{v,\theta}:\theta\in\Theta_v \}$ are all polynomials that depend on $q_v$ and $\widetilde{\mathsf{d}}$. Therefore
\begin{align*}
&\val_T\left(\mathbf{H}/\sqrt{n};\tilde{\mathbf{s}}^{0:t},\tilde{\mathbf{u}}^{1:t},\mathbf{X},\mathbf{Y},\mathbf{Z}\right) \\
=&\frac{1}{n} \sum_{\boldsymbol{\theta} \in \prod_{v \in \mathcal{V}} \Theta_v} \sum_{\mathbf{i} \in [n]^\mathcal{V}} \prod_{v \in \mathcal{V}} q_{v,\theta_v}(\tilde{s}^{0:(t-1)}_{i_v}, \tilde{u}^{1:t}_{i_v}, Z_{i_v})\cdot \prod_{l=1}^L\bigg(\sum_{j=1}^n \frac{H_{i_v, j}^{(l)}}{\sqrt{n}} \tilde{u}^t_{j,l} \bigg)^{\theta_{v,l}} \left(\frac{H}{\sqrt{n}}\right)_{\mathbf{i}|T}\\
=&\frac{1}{n} \sum_{\boldsymbol{\theta} \in \prod_{v \in \mathcal{V}} \Theta_v} \sum_{\mathbf{i} \in [n]^\mathcal{V}} \prod_{v \in \mathcal{V}} q_{v,\theta_v}(\tilde{s}^{0:(t-1)}_{i_v}, \tilde{u}^{1:t}_{i_v}, Z_{i_v})\cdot \prod_{l=1}^L\left[\bigg(\sum_{j=1}^n \frac{H_{i_v, j}^{(l)}}{\sqrt{n}} \tilde{u}^t_{j,l} \bigg)^{\theta_{v,l}} \prod_{(u,v) \in \mathcal{E}:\mathcal{L}[u,v]=l} \frac{H_{i_u,i_v}^{(l)}}{\sqrt{n}}\right].
\end{align*}
For each $\boldsymbol{\theta} \in \prod_{v\in\mathcal{V}}\Theta_v$, we define a new tensor network $T_{\boldsymbol{\theta}}$ from $T$ as follows: (1) for each $v \in \mathcal{V}$, replace the associated polynomial $q_v$ by $q_{v, \theta_v}$; (2) for each $v \in \mathcal{V}$ and $l\in[L]$, connect $v$ with $\theta_{v,l}$ new vertices, where the associated polynomial for each new vertex is $q(\tilde{s}^{0:(t-1)}, \tilde{u}^{1:t}, Z) = \tilde u^t_{l}$ and newly added edges are labeled with $l$. Then the above is precisely
\begin{align*}
\val_T\left(\mathbf{H}/\sqrt{n};\tilde{\mathbf{s}}^{0:t},\tilde{\mathbf{u}}^{1:t},\mathbf{X},\mathbf{Y},\mathbf{Z}\right) = \sum_{\boldsymbol{\theta} \in \prod_{v\in\mathcal{V}}\Theta_v} \val_{T_{\boldsymbol{\theta}}}\left(\mathbf{H}/\sqrt{n};\tilde{\mathbf{s}}^{0:(t-1)},\tilde{\mathbf{u}}^{1:t},\mathbf{X},\mathbf{Y},\mathbf{Z}\right)
\end{align*}
which shows the claim \eqref{eq:sym_reduction_1}. We next claim that for any tensor network $T$ in the variables $(\tilde{s}^{1:(t-1)}, \tilde{u}^{1:t}, X,Y,Z)$, we have
\begin{equation}\label{eq:sym_reduction_2}
\val_{T}\left(\mathbf{H}/\sqrt{n};\tilde{\mathbf{s}}^{0:(t-1)},\tilde{\mathbf{u}}^{1:t},\mathbf{X},\mathbf{Y},\mathbf{Z}\right)=\val_{T^\prime}\left(\mathbf{H}/\sqrt{n};\tilde{\mathbf{s}}^{0:(t-1)},\tilde{\mathbf{u}}^{0:(t-1)},\mathbf{X},\mathbf{Y},\mathbf{Z}\right)
\end{equation}
for a tensor network $T^\prime$ in the variables $(\tilde{s}^{0:(t-1)}, \tilde{u}^{1:(t-1)},X,Y,Z)$. This holds because $\tilde{\mathbf{u}}^t = \widetilde{\mathcal{U}}_t\left(\tilde{\mathbf{s}}^{t-1}, \mathbf{X},\mathbf{Y},\mathbf{Z}\right)$ is itself a polynomial of $(\tilde{\mathbf{s}}^{0:(t-1)},\tilde{\mathbf{u}}^{0:(t-1)},\mathbf{X},\mathbf{Y},\mathbf{Z})$, so for each vertex $v$ of $T$, we may write
\begin{equation*}
q_v\left(\tilde{\mathbf{s}}^{0:(t-1)},\tilde{\mathbf{u}}^{0:(t-1)}, \widetilde{\mathcal{U}}_t\left(\tilde{\mathbf{s}}^{t-1}, \mathbf{X},\mathbf{Y},\mathbf{Z}\right), \mathbf{X},\mathbf{Y},\mathbf{Z}\right) = \tilde q_v\left(\tilde{\mathbf{s}}^{0:(t-1)},\tilde{\mathbf{u}}^{0:(t-1)},\mathbf{X},\mathbf{Y},\mathbf{Z}\right)
\end{equation*}
for some polynomial $\tilde q_v$. Then we can define $T^\prime$ by replacing each polynomial $q_v$ with $\tilde q_v$ and preserving all other structures of $T$.

Having shown the reductions \eqref{eq:sym_reduction_1} and \eqref{eq:sym_reduction_2}, the proof is completed by recursively applying
these reductions for $t,t-1,t-2,\ldots,1$.
\end{proof}
We introduce another lemma below which  addresses universality for tensor networks. Since the notion of tensor networks are generalized to multiple independent sensing matrices, this lemma is a trivial extension of Lemma 2.13 in \cite{wang2022universality}. Interested readers should consult Section 3.1 in \cite{wang2022universality} for a self-contained proof.

\begin{Lemma}\label{lemma:Wignermoments}
Let $\mathbf{x}_1,\ldots,\mathbf{x}_k \in \mathbb{R}^n$ be (random or deterministic) vectors and let $(X_1,\ldots,X_k)$ have finite moments of all orders, such that almost surely as $n \to \infty$,
\begin{equation}\label{eq:WignerMomentsWpConvergence}
(\mathbf{x}_1,\ldots,\mathbf{x}_k) \toW (X_1,\ldots,X_k).
\end{equation}
Suppose that for some constant $C>0$, there holds $\min_{l}d^{(l)}\ge C\log n$.
Then for any diagonal tensor network $T$ in $k$ variables, there is a deterministic value $\limval_T(X_1,\ldots,X_k)$ depending only on $T$ and the joint law of $(X_1,\ldots,X_k)$ such that almost surely,
\begin{align*}
    &\quad\lim_{n \to \infty} \val_T(\mathbf{H}/\sqrt{n};\mathbf{x}_1,\ldots,\mathbf{x}_k)\\
    &=\lim_{n \to \infty} \val_T(\mathbf{W}/\sqrt{n};\mathbf{x}_1,\ldots,\mathbf{x}_k)\\
    &=\limval_T(X_1,\ldots,X_k)
\end{align*}
where $\mathbf{W}$ is the collection $L$ independent $\text{GOE}(n)$ matrices in \eqref{eq:spiked matrix each layer}.
\end{Lemma}

\paragraph{Universality of AMP Iterates.}
\begin{proof}[Proof of Proposition~\ref{prop:algorithmic universality}]

Recall the sequence $\left(\mathbf{s}^0,\ldots,\mathbf{s}^t\right)$ defined in \eqref{eq:coupled AMP update reformulated 1}-\eqref{eq:coupled AMP update reformulated 2}, with respect to a collection of mean-zero sensing matrices $\mathbf{H}$.
When replacing $\mathbf{H}$ with $\mathbf{W}$ in 
\eqref{eq:coupled AMP update reformulated 1}-\eqref{eq:coupled AMP update reformulated 2}, we can define another auxiliary AMP iterates by
\begin{align}
\check{u}_{k,l}^{t+1}=\mathcal{U}^{(l)}_t(\check{s}^t_{k},X_k,Y_k,Z_k), \qquad \check{s}_{i,l}^{t+1}=\sum_{k}\sqrt{\frac{\lambda^{(l)}}{n}}W^{(l)}_{i,k}\check{u}_{k,l}^{t+1}-\lambda^{(l)}\mathsf{b}_t^{(l)}\check{u}_{i,l}^{t},\label{eq:coupled AMP update GOE}
\end{align}
where Onsager coefficients and state evolution iterates are both unchanged. Finally, by setting $\min_{l}d^{(l)}\ge C\log n$, \cite{benaych2020spectral} implies $\max_l\|\mathbf{H}^{(l)}/\sqrt{n}\|\le C$ with high probability. Therefore, equipped with the polynomial approximation technique detailed in Section 3.3 of \cite{wang2022universality} enabled by the polynomial growth condition \eqref{eq:polygrowth}, we are able to apply Lemmas~\ref{lemma:polynomial_tensor_network_decomposition} and \ref{lemma:Wignermoments} to fully conclude the universality of AMP iterates,
\begin{equation*}
    \frac{1}{n}\sum_{i=1}^n\psi(\check{s}^{0:t}_i,X_i,Y_i,Z_i)-\frac{1}{n}\sum_{i=1}^n\psi(s^{0:t}_i,X_i,Y_i,Z_i)\overset{\text{P}}{\rightarrow}0,
\end{equation*}
for any pseudo-Lipschitz test function $\psi$.

By slightly adopting the techniques in Section B.4 of \cite{deshpande2017asymptotic} or Theorem 7.2 of \cite{ma2023community}, we find this sequence of bias-corrected AMP iterates is close to the output of Algorithm~\ref{alg:coupled AMP} in the sense that $\bar{m}_l^t\approx \lambdal \mu_l^t X_l+s_l^t$. Similarly, we get $m_l^t\approx \lambdal \mu_l^t X_l+\check{s}_l^t$.
Therefore, this universality result holds as well for the original sequences with biased sensing matrices, i.e. \eqref{eq:algorithmic universality} holds.
\end{proof}

\section{Regularity of the limiting free energy, and analysis of examples}
\label{sec:misc_proofs}

We prove Propositions~\ref{prop:opt derivative} and \ref{prop:all-in-all-out} in Section~\ref{subsec:opt derivative}. The remaining subsections establish the remaining results for the specific examples.

\subsection{Proof of Propositions~\ref{prop:opt derivative} and \ref{prop:all-in-all-out}}\label{subsec:opt derivative}

\begin{proof}[Proof of Proposition~\ref{prop:opt derivative}]
    (i) Since $\mathcal{F}(\boldsymbol{\lambda},\mathbf{q})$ depends on $(\boldsymbol{\lambda},\mathbf{q})$ only jointly through $\boldsymbol{\lambda}\odot\mathbf{q}$, it is reasonable and convenient to use a reparametrization $\boldsymbol{\gamma}=\left(\lambda^{(1)}q_1,\ldots,\lambda^{(L)}q_L\right)$. Denoting
    \begin{equation*}
        \bar{\mathcal{F}}(\boldsymbol{\gamma}):=\mathbb{E}\log \left(\sum_{x,y}p(x,y|Z)e^{\sum_{l=1}^L\gamma_l X^{(l)}x^{(l)}+\sqrt{\gamma_l}W^{(l)\prime}x^{(l)}}\right)-\sum_{l=1}^{L}\frac{\gamma_l}{2}=\mathcal{F}(\boldsymbol{\lambda},\mathbf{q})-\sum_{l=1}^{L}\frac{\lambdal q_l}{2},
    \end{equation*}
    the optimization objective is then restated as
    \begin{align*}
        \sup_{\mathbf{q}\ge 0}\left[\mathcal{F}(\boldsymbol{\lambda},\mathbf{q})-\sum_{l=1}^{L}\frac{\lambda^{(l)}(q_l^2+2q_l)}{4}\right]=\sup_{\boldsymbol{\gamma}\ge 0}\left[\bar{\mathcal{F}}(\boldsymbol{\gamma})-\sum_{l=1}^L\frac{\gamma_l^2}{4\lambdal}\right].
    \end{align*}
    In the definition of $\bar{\mathcal{F}}$, every $(x,y)$ is weighted by $e^{\sum_{l=1}^L\gamma_lX^{(l)}x^{(l)}+\sqrt{\gamma_l}W^{(l)\prime}x^{(l)}}$, so we introduce the notation
    \begin{equation*}
        \langle f(x,y)\rangle=\frac{\sum_{x,y}p(x,y|Z)e^{\sum_{l=1}^L\gamma_lX^{(l)}x^{(l)}+\sqrt{\gamma_l}W^{(l)\prime}x^{(l)}}f(x,y)}{\sum_{x,y}p(x,y|Z)e^{\sum_{l=1}^L\gamma_lX^{(l)}x^{(l)}+\sqrt{\gamma_l}W^{(l)\prime}x^{(l)}}}.
    \end{equation*}
    Consequently, it yields the first-order derivative of $\bar{\mathcal{F}}$,
    \begin{align*}
        \partial_l \bar{\mathcal{F}}(\boldsymbol{\gamma})=\mathbb{E}\left[\left\langle X^{(l)}x^{(l)}+\frac{1}{2}\sqrt{\frac{1}{\gamma_l}}W^{(l)\prime}x^{(l)}\right\rangle\right]-\frac{1}{2}.
    \end{align*}
    We use Gaussian integration by parts to find that
    \begin{equation*}
        \mathbb{E}\left[W^{(l)\prime}\left\langle x^{(l)}\right\rangle\right]=\mathbb{E}\left[\frac{\partial}{\partial W^{(l)\prime}}\left\langle x^{(l)}\right\rangle\right]=\sqrt{\gamma_l}\left(1-\mathbb{E}\left[\left\langle x^{(l)}\right\rangle^2\right]\right).
    \end{equation*}
    Since Nishimori identity (tower property) implies $\mathbb{E}\left[X^{(l)}\left\langle x^{(l)}\right\rangle\right]=\mathbb{E}\left[\left\langle x^{(l)}\right\rangle^2\right]$, we further derive
    \begin{equation*}
        \partial_l \bar{\mathcal{F}}(\boldsymbol{\gamma})=\frac{1}{2}\mathbb{E}\left[X^{(l)}\left\langle x^{(l)}\right\rangle\right].
    \end{equation*}
    Lastly, since $G(\mathbf{q})=\bar{\mathcal{F}}(\boldsymbol{\lambda}\odot\mathbf{q})-\sum_{l=1}^L\frac{\lambda^{(l)2}q_l^2}{4\lambdal}$, we know
    \begin{equation*}
        \partial_l G(\mathbf{q})=\frac{\lambdal}{2}\mathbb{E}\left[\left\langle x^{(l)}\right\rangle^2\right]-\frac{\lambdal q_l}{2}.
    \end{equation*}
    Then we continue to calculate the second-order derivatives,
    \begin{equation*}
        \partial_{l_2}\partial_{l_1}\bar{\mathcal{F}}(\boldsymbol{\gamma})=\frac{1}{2}\mathbb{E}\left[X^{(l_1)}\frac{\partial}{\partial \gamma_{l_2}}\left\langle x^{(l_1)}\right\rangle\right].
    \end{equation*}
    First it follows from chain rule that 
    \begin{equation*}
        \frac{\partial}{\partial \gamma_{l_2}}\left\langle x^{(l_1)}\right\rangle=\left(\left\langle x^{(l_1)}x^{(l_2)}\right\rangle-\left\langle x^{(l_1)}\right\rangle\left\langle x^{(l_2)}\right\rangle\right)\left(X^{(l_2)}+\frac{1}{2}\sqrt{\frac{1}{\gamma_{l_2}}}W^{(l_2)\prime}\right).
    \end{equation*}
    Again by Gaussian integration by parts, we have
    \begin{align*}
        &\mathbb{E}\left[W^{(l_2)\prime}X^{(l_1)}\left(\left\langle x^{(l_1)}x^{(l_2)}\right\rangle-\left\langle x^{(l_1)}\right\rangle\left\langle x^{(l_2)}\right\rangle\right)\right]\\
        =&\mathbb{E}\left[X^{(l_1)}\frac{\partial}{\partial W^{(l_2)\prime}}\left(\left\langle x^{(l_1)}x^{(l_2)}\right\rangle-\left\langle x^{(l_1)}\right\rangle\left\langle x^{(l_2)}\right\rangle\right)\right]\\
        =&2\sqrt{\gamma_{l_2}}\cdot\mathbb{E}\left[X^{(l_1)}\left(\left\langle x^{(l_1)}\right\rangle\left\langle x^{(l_2)}\right\rangle^2-\left\langle x^{(l_1)}x^{(l_2)}\right\rangle\left\langle x^{(l_2)}\right\rangle\right)\right].
    \end{align*}
    Putting together, we would finally arrive at
    \begin{equation*}
        \partial_{l_2}\partial_{l_1}\bar{\mathcal{F}}(\boldsymbol{\gamma})=\frac{1}{2}\mathbb{E}\left[\left(\left\langle x^{(l_1)}x^{(l_2)}\right\rangle-\left\langle x^{(l_1)}\right\rangle\left\langle x^{(l_2)}\right\rangle\right)^2\right].
    \end{equation*}
    It follows from rescaling that
    \begin{equation*}
        \partial_{l_2}\partial_{l_1}G(\boldsymbol{\gamma})=\frac{1}{2}\mathbb{E}\left[\left(\left\langle x^{(l_1)}x^{(l_2)}\right\rangle-\left\langle x^{(l_1)}\right\rangle\left\langle x^{(l_2)}\right\rangle\right)^2\right]-\frac{1}{2\lambda^{(l_1)}}\mathbbm{1}\{l_1=l_2\}.
    \end{equation*}
    (ii) If $\boldsymbol{\gamma}\in(0,\infty)^L$ is not on the boundary, the conclusion follows naturally. Since $\mathbf{0}$ is also a solution to $T(\boldsymbol{\gamma})=\boldsymbol{\gamma}$, it suffices to deal with the case that there is a proper subset $S\subset[L]$ such that $\gamma_S=0,\gamma_{S^c}>0$.
    
    For any $j\notin S$, by the Lagrange condition, we must have $\partial_j G(\boldsymbol{\gamma})=0$ leading to $T_j(\boldsymbol{\gamma})=\gamma_j>0$. 
    For any $j\in S$, the optimality condition at boundaries imply $\partial_j G(\boldsymbol{\gamma})\le 0$. But by its definition, $\partial_j G(\boldsymbol{\gamma})=T_j(\boldsymbol{\gamma})-\gamma_j=T_j(\boldsymbol{\gamma})\ge 0$. So we also have $T_j(\boldsymbol{\gamma})=\gamma_j=0$ for each $j\in S$.
\end{proof}

\begin{proof}[Proof of Proposition~\ref{prop:all-in-all-out}]
    For simplicity of proof presentation, we only consider the spiked matrices model \eqref{eq:spiked matrix each layer}, as our universality results Propositions~\ref{prop:Universality MI qualitative} and \ref{prop:universality mmse} allow us to translate the results to the  original graph model \eqref{eq:original graph model}. 

    Since $\boldsymbol{\lambda}\notin D$, denote $\mathbf{q}^\ast$ as the unique maximizer. The last assertion of Proposition~\ref{prop:opt derivative} establishes that $\mathbf{q}^\ast$ to solves a fixed point  equation $\mathbf{q}^\ast = T(\boldsymbol{\lambda}\odot\mathbf{q}^\ast)$. Now suppose a certain layer $\bX^{(l)}$ is weakly recoverable. Then Theorem~\ref{thm:asymptotic MI} and Lemma~\ref{lemma:MSE to weak recovery} together imply $q_l^\ast>0$. As a result, in the effective scalar channel \eqref{eq:scalar channel}, only $A^{(l)\prime}$ suffices to weakly recover $X^{(l)}$, since
    \begin{equation*}
        \mathbb{E}\left[X^{(l)}\E\left(X^{(l)}|A^{(l)\prime}\right)\right]=\E\left[\tanh\left(\lambdal \mathbf{q}^\ast_l+\sqrt{\lambdal \mathbf{q}^\ast_l}\mathcal{N}(0,1)\right)\right]>0,
    \end{equation*}
    as long as $\lambdal \mathbf{q}^\ast_l>0$. For any other $l^\prime\in[L]$, since $X^{(l)}$ and $X^{(l^\prime)}$ are correlated under $p$,
    \begin{align*}
    q^\ast_{l^\prime}&=\mathbb{E}\left[X^{(l^\prime)}\E\left(X^{(l^\prime)}|A^{\prime}\right)\right]\\
    &\ge \left|\mathbb{E}\left[X^{(l^\prime)}\E\left(X^{(l)}|A^{(l)\prime}\right)\right]\right|\\
    &=\left|\mathbb{E}\left[X^{(l)}X^{(l^\prime)}\right]\right|\cdot \mathbb{E}\left[X^{(l)}\E\left(X^{(l)}|A^{(l)\prime}\right)\right]>0.
    \end{align*}
    This guarantees the weak recovery feasibility of $\bX^{(l^\prime)}$ by Theorem~\ref{thm:asymptotic MI}(ii). In the same way, for any $l_1\in[L_1]$, we also know
    \begin{equation*}
        \mathbb{E}\left[Y^{(l_1)}\E\left(X^{(l^\prime)}|A^{\prime}\right)\right]>0
    \end{equation*}
    which guarantees the weak recovery feasibility of $\bY^{(l_1)}$ by Theorem~\ref{thm:asymptotic MI}(iii). To complete the proof, it suffices to show nothing can be weakly recovered if $\mathbf{q}^\ast=0$. And this is immediately true due to Theorem~\ref{thm:asymptotic MI}(ii),(iii).
\end{proof}

\subsection{Proofs for Multilayer SBM}\label{subsec:pf multilayer}

\begin{proof}[Proof of Corollary~\ref{cor:global and local MI multilayer}]
    (i) By the definition of conditional mutual information, we know
    \begin{equation*}
        I(\bY;\bG)=I(\bX,\bY;\bG)-I(\bX;\bG|\bY)=I(\bX,\bY;\bG)-\sum_{l}I(\bX^{(l)};\bG^{(l)}|\bY),
    \end{equation*}
    where the first term $I(\bX,\bY;\bG)$ is already derived in  Theorem~\ref{thm:multilayer joint MI}.  The second term $I(\bX;\bG|\bY)$ is easily decoupled by conditional independence. Now we can treat $\bY$ as node-wise independent side information in the general framework presented in Section~\ref{subsec:general model setup}, and we are only given one graph $\bG^{(l)}$ to make inference on $\bX^{(l)}$. The free energy functional of this channel is given in \eqref{eq:free energy multilayer conditional global}. Again, Theorem~\ref{thm:asymptotic MI} and Proposition~\ref{prop:Universality MI qualitative} together imply
    \begin{align*}
        &\quad\lim_{d^{(l)}\rightarrow\infty}\lim_{n\rightarrow\infty}\frac{1}{n}I(\bX^{(l)};\bG^{(l)}|\bY)=\lim_{n\rightarrow\infty}\frac{1}{n}I(\bX^{(l)};\bA^{(l)}|\bY)\\
        &=\frac{\lambdal}{4}-\sup_{q_l\ge 0}\left[\mathcal{F}^{(l)}_{\ML,\rho}(\lambdal,q_l)-\frac{\lambda^{(l)}(q_l^2+2q_l)}{4}\right].
    \end{align*}
    Finally \eqref{eq:global MI multilayer} follows from plugging these expressions.

    (ii) As for each individualized layer $\bX^{(l)}$, it holds that
    \begin{equation*}
        I(\bX^{(l)};\bG)=I(\bX,\bY;\bG)-I(\bX^{(-l)},\bY;\bG^{(-l)}|\bX^{(l)}).
    \end{equation*}
    To deal with the second term, we treat $\bX^{(l)}$ as node-wise independent side information in the general framework presented in Section~\ref{subsec:general model setup}. Conditioned on $\bX^{(l)}$, we want to infer $(\bX^{(-l)},\bY)$ from $\bG^{(-l)}$. This input prior corresponds to energy functional defined in \eqref{eq:free energy multilayer conditional local}. Theorem~\ref{thm:asymptotic MI} and Proposition~\ref{prop:Universality MI qualitative} together imply
    \begin{align*}
        &\quad\lim_{d^{(l)}\rightarrow\infty}\lim_{n\rightarrow\infty}\frac{1}{n}I(\bX^{(-l)},\bY;\bG^{(-l)}|\bX^{(l)})=\lim_{n\rightarrow\infty}\frac{1}{n}I(\bX^{(-l)},\bY;\bA^{(-l)}|\bX^{(l)})\\
        &=\sum_{l_1\neq l}\frac{\lambda^{(l_1)}}{4}-\sup_{\mathbf{q}_{-l}\ge 0}\left[\bar{\mathcal{F}}^{(l)}_{\ML,\rho}(\boldsymbol{\lambda}^{(-l)},\mathbf{q}_{-l})-\sum_{l_1\neq l}\frac{\lambda^{(l_1)}(q_{l_1}^2+2q_{l_1})}{4}\right].
    \end{align*}
    We plug in this expression to finally get \eqref{eq:local MI multilayer}.
\end{proof}

\begin{proof}[Proof of Lemma~\ref{lemma:ML denoiser}]
By reparametrizing $\rho$ into $\bar{\rho}$, we know $p_{\ML}(x|y)\propto\prod_l\exp\left(\bar{\rho}yx^{(l)}\right)$. The likelihood also satisfies a similar form $p(m_l|x^{(l)})\propto\exp(m_lx^{(l)})$, mainly due to the fact that $x^{(l)2}=1$. Therefore, after observing $\sqrt{\lambdal q_l^t}A^{(l)\prime}=m_l,\forall l\in[L]$, the joint posterior can be denoted as
    \begin{equation*}
        p(x,y|m)\propto p(m|x)p_{\ML}(x|y)\propto \exp\left(\sum_{l}m_lx^{(l)}+\bar{\rho}yx^{(l)}\right).
    \end{equation*}
    Consequently, we can compute the marginals by
    \begin{align*}
        p(x^{(l)}=1|m)\propto&\sum_{y,x^{(-l)}}\exp\left(\sum_{l}m_lx^{(l)}+\bar{\rho}yx^{(l)}\right)\\
        \propto&\sum_{y}\exp\left(m_l+\bar{\rho}y\right)\prod_{l^\prime\neq l}\cosh(m_{l^\prime}+\bar{\rho}y).
    \end{align*}
    Similarly, we would have $ p(x^{(l)}=-1|m)\propto\sum_{y}\exp\left(-m_l-\bar{\rho}y\right)\prod_{l^\prime\neq l}\cosh(m_{l^\prime}+\bar{\rho}y)$. As a result,
    \begin{align*}
        \mathcal{E}^{(l)}_t(m_{1:L})&=p(x^{(l)}=1|m)-p(x^{(l)}=-1|m)\\
        &=\frac{\sum_{y}\exp\left(m_l+\bar{\rho}y\right)\prod_{l^\prime\neq l}\cosh(m_{l^\prime}+\bar{\rho}y)-\sum_{y}\exp\left(-m_l-\bar{\rho}y\right)\prod_{l^\prime\neq l}\cosh(m_{l^\prime}+\bar{\rho}y)}{\sum_{y}\exp\left(m_l+\bar{\rho}y\right)\prod_{l^\prime\neq l}\cosh(m_{l^\prime}+\bar{\rho}y)+\sum_{y}\exp\left(-m_l-\bar{\rho}y\right)\prod_{l^\prime\neq l}\cosh(m_{l^\prime}+\bar{\rho}y)}\\
        &=\frac{\prod_{l^\prime\neq l}\cosh\left(m_{l^\prime}+\bar{\rho}\right)\sinh\left(m_l+\bar{\rho}\right)+\prod_{l^\prime\neq l}\cosh\left(m_{l^\prime}-\bar{\rho}\right)\sinh\left(m_l-\bar{\rho}\right)}{\prod_{l^\prime}\cosh\left(m_{l^\prime}+\bar{\rho}\right)+\prod_{l^\prime}\cosh\left(m_{l^\prime}-\bar{\rho}\right)}.
    \end{align*}
    coinciding with \eqref{eq:ML denoiser}.
\end{proof}

\begin{proof}[Proof of Proposition \ref{prop:multilayer-threshold}]
Note that for almost all $\boldsymbol{\lambda}$ except for a set of measure zero, the maximizer $\mathbf{q}^\ast$ of RHS in \eqref{eq:limiting MI} is unique. 

Note that the Hessian matrix $\nabla^2 G(\mathbf{0})$ is
\begin{equation*}
    \frac{1}{2}\begin{pmatrix}
        \lambda^{(1)2}-\lambda^{(1)} & \lambda^{(1)}\lambda^{(2)}(1-2\rho)^4 & & \lambda^{(1)}\lambda^{(L)}(1-2\rho)^4\\
        \lambda^{(1)}\lambda^{(2)}(1-2\rho)^4 & \lambda^{(2)2}-\lambda^{(2)} & & \lambda^{(2)}\lambda^{(L)}(1-2\rho)^4\\
        & & \ddots & \\
        \lambda^{(1)}\lambda^{(L)}(1-2\rho)^4 & & & \lambda^{(L)2}-\lambda^{(L)}
    \end{pmatrix},
\end{equation*}
which can be rewritten as:
\begin{equation*}
    \frac{(1-2\rho)^4}{2}\begin{pmatrix}
           \lambda^{(1)} \\
           \lambda^{(2)} \\
           \vdots \\
           \lambda^{(L)}\end{pmatrix}^\intercal\begin{pmatrix}
           \lambda^{(1)} &
           \lambda^{(2)} &
           \ldots &
           \lambda^{(L)}\end{pmatrix}-\frac{1}{2}\text{diag}(\lambda^{(l)}-(1-(1-2\rho)^4)\lambda^{(l)2})
\end{equation*}
Denote the maximum eigenvector of $\nabla^2 G(\mathbf{0})$ as $\mathbf{v} = (v_1,\ldots, v_L)^\intercal$. We first note that all entries of $v$ must either be non-negative or non-positive. 
%
To see this, we note 
\begin{equation*}
    v^\intercal\nabla^2G(\mathbf{0})v = \frac{(1-2\rho)^4}{2}\left(\sum_{l = 1}^L \lambda^{(l)}v_l\right)^2-\frac{1}{2}\sum_{l = 1}^L (\lambda^{(l)}-(1-(1-2\rho)^4)\lambda^{(l)})v_l^2.
\end{equation*}
Without loss of generality, we assume $\sum_{l = 1}^L\lambda^{(l)}v_l>0$, then we can flip the sign of any negative $v_l$ to increase $\sum_{l = 1}^L\lambda^{(l)}v_l$ while keeping $\sum_{l = 1}^L (\lambda^{(l)}-(1-(1-2\rho)^4)\lambda^{(l)})v_l^2$ and $\sum_{l = 1}^L v_l^2$ fixed, thus increasing $\mathbf{v}^\intercal\nabla^2G(0)\mathbf{v}$. This contradicts the fact that $\mathbf{v}$ is the leading eigenvector and establishes the claim. 

By direct computation, 
\begin{equation}
\label{eq:determinant-G0-inhomo-SBM}
    |\nabla^2 G(\mathbf{0})| = \left(-\frac{1}{2}\right)^L\left(1-\sum_{l = 1}^L \frac{(1-2\rho)^4\lambda^{(l)}}{1-(1-(1-2\rho)^4)\lambda^{(l)}}\right).
\end{equation}
Thus if \eqref{eq:threshold-inhomo-SBM} holds, the leading eigenvalue for $\nabla^2 G(\mathbf{0})$ must be larger than $0$. The previous observation indicates that the leading eigenvector must point to the first quadrant, which suggests that  $\mathbf{0}$ is only a saddle point and that $\mathbf{q}^\ast\neq \mathbf{0}$. Note that $\mathbf{q}^\ast$ enjoys an ``all-or-nothing'' behavior in the sense that $q_l^\ast$'s are either all equal to zero or none of them equal to zero. Thus if \eqref{eq:threshold-inhomo-SBM} holds, we achieve simultaneous weak recovery for all $\mathbf{X}^{(l)}$'s as well as $\mathbf{Y}$.

We now turn to the converse assertion: if 
\begin{equation}
\label{eq:inverse-threshold-inhomo-SBM}
        \sum_{l = 1}^L \frac{(1-2\rho)^4\lambda^{(l)}}{1-(1-(1-2\rho)^4)\lambda^{(l)}}<1,
\end{equation}
$\mathbf{0}$ is a local optimizer. Note that $\nabla^2 G(\mathbf{0})$ is a negative definite matrix perturbed by a rank one matrix, thus it admits at most one positive eigenvalue. Thus if $\sum_{l = 1}^L \frac{(1-2\rho)^4\lambda^{(l)}}{1-(1-(1-2\rho)^4)\lambda^{(l)}}<1$, $|\nabla^2G(\mathbf{0})|<0$ when $L$ is odd and $|\nabla^2G(\mathbf{0})|>0$ when $L$ is even. Then the largest eigenvalue of $\nabla^2 G(\mathbf{0})$ is negative, the proof is thus completed.

\end{proof}

\begin{proof}[Proof of Proposition~\ref{prop:ML conjectured}]
Recall that we are solving $\mathbf{q}=T^\ML(\blam\odot\mathbf{q})$ to find the maximizer of free energy functional. 
All we need here for $T^\ML$ is the monotonicity and the conjectured weakened notion of strict concavity. In fact, this type of mappings actually has a name ``standard interference mapping", and is already widely studied in the context of communication systems \cite{yates1995framework, cavalcante2015elementary, cavalcante2019connections, piotrowski2022fixed}.
\begin{itemize}
    \item[(i)] Suppose $\mathbf{q}\mapsto T^\ML(\blam\odot\mathbf{q})$ has two distinct non-zero fixed points $\mathbf{q},\mathbf{q}^\prime\in[0,+\infty)^L$. As suggested by Proposition~\ref{prop:all-in-all-out}, all coordinates of $\mathbf{q},\mathbf{q}^\prime$ should be positive. Without loss of generality, there exists some $j_0$ such that $q_{j_0} <q_{j_0}^\prime$. Hence, there exists a $t>1$ such that $t\mathbf{q}\ge\mathbf{q}^\prime$ and $tq_j=q_j^\prime$ for some $j$. Then monotonicity and concavity together imply
    \begin{equation*}
        q_j^\prime = T^\ML_j(\blam\odot\mathbf{q}^\prime)\le T^\ML_j(t\blam\odot\mathbf{q})< tT^\ML_j(\blam\odot\mathbf{q})=tq_j,
    \end{equation*}
    contradicting to $t\mathbf{q}_j=\mathbf{q}_j^\prime$. Therefore, $\mathbf{q}=T^\ML(\blam\odot\mathbf{q})$ only has at most one non-zero solution. When there is no non-zero solution, then origin $\mathbf{0}$ becomes unique maximizer of \eqref{eq:multilayer joint MI}. If there is indeed a non-zero fixed point $\mathbf{q}^\ast$, there must be
    \begin{equation*}
        T^\ML_j(\alpha\blam\odot\mathbf{q}^\ast)>\alpha T^\ML_j(\blam\odot\mathbf{q}^\ast)=\alpha\mathbf{q}_j^\ast,\quad\forall j\in[L],0<\alpha<1.
    \end{equation*}
    Since $T^\ML(\blam\odot\mathbf{q})-\mathbf{q}$ is the gradient of the free energy functional, the functional must have a strict bigger value at $\mathbf{q}^\ast$ than $\mathbf{0}$. So $\mathbf{q}^\ast$ would be the unique maximizer.

    \item[(ii)] If condition \eqref{eq:threshold-inhomo-SBM} fails, $\mathbf{0}$ must be a \textit{local} maximizer to \eqref{eq:multilayer joint MI}. Continued from the previous argument, if $\mathbf{q}\mapsto T^\ML(\blam\odot\mathbf{q})$ has one additional non-zero fixed point $\mathbf{q}^\ast$, $T^\ML(\blam\odot\mathbf{q})<\mathbf{q}$ must hold along the segment connecting $\mathbf{0}$ and $\mathbf{q}^\ast$. Then $\mathbf{0}$ cannot be a local maximizer, contradiction! So $\mathbf{0}$ is the only fixed point, thus the global maximizer.

    \item[(iii)] If condition \eqref{eq:threshold-inhomo-SBM} holds, as established in \cite{piotrowski2022fixed}, state evolution iterates converge to $\mathbf{q}^\ast$ geometrically, as long as the initialization $\mathbf{q}^0\neq \mathbf{0}$.
\end{itemize}

\end{proof}

\subsection{Proofs for Dynamic SBM}
\begin{proof}[Proof of Corollary~\ref{cor:local dyn}]
    By the definition of conditional mutual information, we know
    \begin{align*}
        I(\bX^{(l)};\bG)&=I(\bX;\bG)-I(\bX^{(<l)},\bX^{(>l)};\bG|\bX^{(l)})\\
        &=I(\bX,\bY;\bG)-I(\bX^{(<l)};\bG^{(<l)}|\bX^{(l)})-I(\bX^{(>l)};\bG^{(>l)}|\bX^{(l)}),
    \end{align*}
    where the first term $I(\bX;\bG)$ is already derived in  Theorem~\ref{thm:dyn joint MI}.  The last two terms can both be put into our general model easily. Now for $I(\bX^{(>l)};\bG^{(>l)}|\bX^{(l)})$, we can treat $\bX^{(l)}$ as node-wise independent known side information in the general framework presented in Section~\ref{subsec:general model setup}, and we are only given graphs $\bG^{(>l)}$ to make inference on $\bX^{(>l)}$. Corresponding free energy functionals are given in the corollary statement. Again, Theorem~\ref{thm:asymptotic MI} and Proposition~\ref{prop:Universality MI qualitative} together imply
    \begin{align*}
        &\quad\lim_{n\rightarrow\infty}\frac{1}{n}I(\bX^{(>l)};\bG^{(>l)}|\bX^{(l)})=\lim_{n\rightarrow\infty}\frac{1}{n}I(\bX^{(>l)};\bA^{(>l)}|\bX^{(l)})\\
        &=\frac{\sum_{l^\prime>l}\lambda^{l^\prime}}{4}-\sup_{\mathbf{q}_{>l}\ge 0}\left[\mathcal{F}_{\Dyn,\rho}(\boldsymbol{\lambda}^{(>l)},\mathbf{q}_{>l})-\sum_{l^\prime>l}\frac{\lambda^{(l^\prime)}(q_{l^\prime}^2+2q_{l^\prime})}{4}\right].
    \end{align*}
    Similar limiting formula holds for $I(\bX^{(<l)};\bG^{(<l)}|\bX^{(l)})/n$ as well. Finally \eqref{eq:local dyn} follows from plugging these expressions. 
\end{proof}

\begin{proof}[Proof of Lemma~\ref{lemma:dyn denoiser}]
Focus on the scalar channel \eqref{eq:scalar channel dyn} with input distribution $p_\Dyn$. And input $m_{1:L}$ of this algorithm satisfies $m_{l}=\sqrt{\lambdal q_l^t}A^{(l)\prime}$. Then the likelihood function corresponding to each $m_l$ is
\begin{equation*}
p\left(m_l|X^{(l)}=x^{(l)}\right)\propto\exp(m_lx^{(l)})=:q_l\left(x^{(l)}\right).
\end{equation*}
Therefore, the forward recursion is initialized by
\begin{align*}
    g_1(1)=p(X^{(1)}=1)q_1(1)&\propto p(X^{(1)}=1|m_1),\\
    g_1(-1)=p(X^{(1)}=-1)q_1(-1)&\propto p(X^{(1)}=-1|m_1).
\end{align*}
and proceeds by computing for each $x^{(l)}\in\{\pm 1\}$,
\begin{align*}
g_l(x^{(l)})&=\sum_{x^{(l-1)}}p(X^{(l)}=x^{(l)}|X^{(l-1)}=x^{(l-1)})\exp(m_lx^{(l-1)})g_{l-1}\left(x^{(l-1)}\right)\\
&=\sum_{x^{(1)}}p(X^{(1)}=x^{(1)})q_1(x^{(1)})\left(\sum_{x^{(2)}}p(X^{(2)}=x^{(2)}|X^{(1)}=x^{(1)})q_2(x^{(2)})\left(\sum_{x^{(3)}}\cdots\right)\right)\\
&\propto p(X^{(l)}=x^{(l)}|m_{1:l}).
\end{align*}
for each $l=2,\ldots,L$. Henceforth, the backward recursion starts with
\begin{align*}
    \mathcal{E}^{(L)}&=\frac{g_L(1)-g_L(-1)}{g_L(1)+g_L(-1)}=p(X^{(L)}=1|m_{1:L})-p(X^{(L)}=-1|m_{1:L})\\
    &=\mathbb{E}\left[X^{(L)}\bigg| \sqrt{\lambdal q_l^t}A^{(l)\prime}=m_l,\quad \forall l\in[L]\right]
\end{align*}
and then proceeds by
\begin{align*}
\mathcal{E}^{(l)}&=\sum_{x^{(l)},x^{(l+1)}}x^{(l)}\frac{1+x^{(l+1)}\mathcal{E}^{(l+1)}}{2}\cdot\frac{g_{l}(x^{(l)})\left[(1-\rho)\mathbbm{1}(x^{(l+1)}=x^{(l)})+\rho\mathbbm{1}(x^{(l+1)}\neq x^{(l)})\right]\exp\left(m_{l+1}x^{(l+1)}\right)}{g_{l+1}(x^{(l+1)})}\\
&=\sum_{x^{(l)},x^{(l+1)}}x^{(l)}p(X^{(l+1)}=x^{(l+1)}|m_{1:L})\frac{p(X^{(l)}=x^{(l)},X^{(l+1)}=x^{(l+1)}|m_{1:(l+1)})}{p(X^{(l+1)}=x^{(l+1)}|m_{1:(l+1)})}\\
&=\sum_{x^{(l)},x^{(l+1)}}x^{(l)}p(X^{(l+1)}=x^{(l+1)}|m_{1:L})p(X^{(l)}=x^{(l)}|X^{(l+1)}=x^{(l+1)},m_{1:(l+1)})\\
&=\sum_{x^{(l)}}x^{(l)}p(X^{(l)}=x^{(l)}|m_{1:L})\\
&=\mathbb{E}\left[X^{(L)}\bigg| \sqrt{\lambdal q_l^t}A^{(l)\prime}=m_l,\quad \forall l\in[L]\right],
\end{align*}
for $l$ going from $L-1$ down to $1$. Therefore, Algorithm~\ref{alg:denoiser dyn} outputs an exact posterior mean of the scalar channel. Consisting of two recursions, this algorithm runs efficiently in $O(L)$ time.
\end{proof}

\begin{proof}[Proof of Proposition~\ref{prop:dynamic-threshold}]
    Note that if  $\lambda^{(1)} = \lambda^{(2)} = \cdots = \lambda^{(L)} = \lambda$,  the Hessian matrix $\nabla^2G(\mathbf{0})$ is
\begin{equation*}
    \frac{1}{2}\begin{pmatrix}
        \lambda^{2}-\lambda & \lambda^2(1-2\rho)^2 & & \lambda^2(1-2\rho)^{2(L-1)}\\
        \lambda^{2}(1-2\rho)^2 & \lambda^{2}-\lambda & & \lambda^{2}(1-2\rho)^{2(L-2)}\\
        & & \ddots & \\
        \lambda^{2}(1-2\rho)^{2(L-1)} & & & \lambda^{2}-\lambda
    \end{pmatrix}.
\end{equation*}
Define the Toeplitz matrix $K$ as 
\begin{equation*}
    K = \left((1-2\rho)^{2|r-s|}\right)_{r,s = 1}^L,
\end{equation*}
then $\nabla^2G(\mathbf{0})$ can be rewritten as 
\begin{equation*}
    \frac{\lambda^2}{2}K-\frac{\lambda}{2}I,
\end{equation*}
where $I$ is the identity matrix.

Note that $K$ belongs to the family of Kac-Murdock-Szego matrices, and thus its largest eigenvalue is \citep{trench2010spectral,grenander1958toeplitz}
\begin{equation*}
    \lambda_{\max}(K) = \frac{1-(1-2\rho)^4}{1-2(1-2\rho)^2\cos\theta_\ast+(1-2\rho)^4},
\end{equation*}
where $\theta_\ast\in(0,\pi)$ is the minimum solution of equation
\begin{equation*}
        0 = \sin[(L+1)\theta_\ast]-2(1-2\rho)^2\sin[ L\theta_\ast]+(1-2\rho)^4\sin[(L-1)\theta_\ast].
    \end{equation*}
Furthermore, setting the largest eigenvector of $K$ as $v = (v_1,\ldots, v_L)^\intercal$, 
\begin{equation*}
    v_r = \cos{\frac{(L-2r+1)\theta_\ast}{2}},\quad 1\leq r\leq L.
\end{equation*}
It is proved in \cite{trench2010spectral} and \cite{grenander1958toeplitz} that $0<\theta_\ast<\pi/(L+1)$, thus $v_r>0$ for $1\leq r\leq L$. Note that $v$ is also the largest eigenvector of $\nabla^2G(\mathbf{0})$. Thus if $\lambda^2\lambda_{\max}(K)-\lambda>0$, or equivalently $\lambda>1/\lambda_{\max}(K)$, there exists $\mathbf{v}\in(0,\infty)^L$, such that $\mathbf{v}^\intercal \nabla^2G(\mathbf{0})\mathbf{v}>0$. This completes the proof. 
\end{proof}

\begin{proof}[Proof of Proposition~\ref{prop:dyn conjectured}]
    This proof is essentially the same as that of Proposition~\ref{prop:ML conjectured}, and is thus omitted. 
\end{proof}

\vspace{2mm}
\noindent
\textbf{Acknowledgments:} SS was partially supported by NSF DMS-CAREER 2239234, ONR N00014-23-1-2489 and AFOSR FA9950-23-1-0429.

\bibliography{sample}
\bibliographystyle{plainnat}

\end{document}